\documentclass[12pt,english]{article}

\usepackage[latin1]{inputenc}
\usepackage[english]{babel}
\usepackage{amsfonts}
\usepackage{latexsym}
\usepackage{amsmath}
\usepackage{mathrsfs}
\usepackage{amssymb}
\usepackage{graphpap}
\usepackage{graphicx}
\usepackage{epsfig}

\def\Z{{\mathbb Z}}

\def\R{{\mathbb R}}

\def\IP{{\mathbb P}}
\def\IE{{\mathbb E}}
\def\e{{\mathfrak E}}

\def\Po{{\mathtt P}_{\omega}}
\def\Eo{{\mathtt E}_{\omega}}
\newcommand{\eps}{\varepsilon}
\newcommand{\gam}{\gamma}
\newcommand{\vphi}{\varphi}
\newcommand{\eqlaw}{\stackrel{\text{\tiny law}}{=}}
\newcommand{\argmax}{\mathop{\mathrm{arg\,max}}}
\newcommand{\argmin}{\mathop{\mathrm{arg\,min}}}
\newtheorem{theo}{Theorem}[section]
\newtheorem{lm}{Lemma}[section]
\newtheorem{df}{Definition}[section]

\allowdisplaybreaks

\title{On the moments of the meeting time of independent random walks in random environment}
\author{Christophe~Gallesco$^{\dagger}$}

\begin{document}

\bibliographystyle{plain}

\maketitle
{\footnotesize 

\noindent $^{\dagger}$ Instituto de Matem\'atica e Estat\'istica, Universidade de
S\~ao Paulo, rua do Mat\~ao 1010, CEP 05508--090, S\~ao Paulo, SP, Brazil,
\noindent e-mail: \texttt{gallesco@ime.usp.br}
}

\begin{abstract}
We consider, in the continuous time version, $\gamma$ independent random walks on $\mathbb{Z_+}$ in random environment in the Sinai's regime. Let $T_\gam$ be the first meeting time of one pair of the $\gamma$ random walks starting at different positions. We first show that the tail of the quenched distribution of $T_\gamma$, after a suitable rescaling, converges in probability, to some functional of the Brownian motion. Then we compute the law of this functional. Eventually, we obtain results about the moments of this meeting time. Being $\Eo$ the quenched expectation, we show that, for almost all environments $\omega$, $\Eo[T_\gamma^{c}]$ is finite for $c<\gamma(\gamma-1)/2$ and infinite for $c>\gamma(\gamma-1)/2$.
\\[.3cm] {\bf Keywords:} random walk in random environment, Sinai's regime,
 $t$-stable point, meeting time, coalescing time\\
 {\bf AMS 2000 Mathematics Subject Classification:} 60K37
\end{abstract}

\section{Introduction and results}
\subsection{Preliminary facts}
Let $\omega=(\omega^+_x,\omega^-_x)_{x \in \mathbb{Z_+}}$ be a sequence of positive vectors. Then, consider $\gamma$ independent
continuous-time Markov chains
$\xi_i=(\xi_i(t))_{t \geq 0}$, $1\leq i\leq \gamma$, on~$\mathbb{Z_+}$, which jump from~$y$ to~$y+1$
with rate~$\omega^+_y$, and to~$y-1$ with rate~$\omega^-_y$ if $y\geq1$ and jump from~$0$ to~$1$ with rate~$\omega^+_0$. Now, we suppose that $\omega$ is a fixed realization of an i.i.d. sequence of positive random variables.
We refer to $\omega$ as {\it the environment},
and to the $\xi_i$-s, $1\leq i\leq \gamma$, as the {\it random walks in
the random environment} $\omega$. We will denote by
$\IP, \IE$ the probability and expectation with respect to
$\omega$, and by  $\Po$, $\Eo$  the (so-called ``quenched")
probability and expectation for random walks in the fixed
environment $\omega$. Observe that the $\xi_i$-s are independent under the quenched law $\Po$ but not in general under the annealed law $\IP \times \Po$. As usual, we will use the notation $\Po^{\vec{x}}$ with $\vec{x}=(x_1,\dots, x_\gam)$ for the quenched law of $(\xi_1,\dots, \xi_\gam)$ starting from $\xi_i(0)=x_i$,  with $x_i<x_{i+1}$ for $1\leq i< \gamma$. Nevertheless, for the sake of brevity, we will omit the superscript ${\vec{x}}$ whenever the initial positions are $\xi_{i}(0)=i$, for $1\leq i\leq \gamma$ ( for the sake of simplicity, we decided to choose these intitial positions, but our results do not depend on the initial positions of the $\xi_i$-s). We also formally define the meeting time of one pair of the random walks as
\[T_{\gamma} =\inf\{s>0;\xi_i(s)=\xi_{j}(s), \mbox{for some $1\leq i\neq j \leq \gam$}\}.\]
In this paper we study only the case of {\it Sinai's regime}, which
means that the following condition is satisfied:

\medskip
\noindent
{\bf Condition S.} We have
\[
\IE\ln \frac{\omega^+_0}{\omega^-_0} = 0,\quad
\sigma^2 := \IE  \ln^2 \frac{\omega_0^+}{\omega_0^-}
\in (0, +\infty).
\]

The one-dimensional RWRE has been a model of constant interest in the last three decades. The first important results are due to Solomon~\cite{solom}, Sinai~\cite{Sinai} and Kesten, Kozlov and Spitzer~\cite{KKS}. Solomon showed that the first part of condition S is equivalent to recurrence for a random walk $\xi$ on $\Z$. Sinai proved that abnormal diffusion takes place in the recurrent case:~$\xi(t)$ is of order $\ln^2 t$; Kesten, Kozlov and Spitzer computed the limit law of the random walk in the transient case.
Since then, we know that this model presents many interesting properties. For example, the large deviation properties were studied by Greven and den Hollander~\cite{GdH}. The moderate deviation properties were studied by Comets and Popov in~\cite{CP} and \cite{CP2}. Moreover, the model presents many variations of interest, for example, in~\cite{CMP}, Comets, Menshikov and Popov use Lyapunov functions method to study random strings and in~\cite{HuShidiff} Hu and Shi use stochastic calculus to obtain results on moderate deviations for the one-dimensional diffusion in a Brownian potential, which is the continuous-space analogue of the RWRE.

 More recently, in the sub-ballistic transient case, Fribergh, Gantert and Popov~\cite{PGF} worked on moderate deviations and Enriquez, Sabot and Zindy~\cite{ESZ} refined, in particular, the results of \cite{KKS}. In the Sinai's regime, Dembo, Gantert, Peres and Shi~\cite{DGPS} and Gantert, Peres and Shi~\cite{GPS} studied the properties of the local times. For surveys on the subject, the reader is referred to the lecture notes of a course by Zeitouni~\cite{Z}, to the book of Hughes~\cite{hughes} and to the stochastic calculus approach by Shi~\cite{shi}. In addition, we mention that the results we obtain in Theorem \ref{theo2} contrast with those of the analagous deterministic model which is the simple random walk on $\Z_+$. In the case of 2 independent random walks, it is well known that $E[\tau^{c}]<\infty$ for $c<1/2$ and $E[\tau^{c}]=\infty$ for $c\geq 1/2$. In the case of 3 independent random walks, we have that $E[\tau^{c}]<\infty$ for $c<3/2$ and $E[\tau^{c}]=\infty$ for $c>3/2$ (see for example \cite{BFMP}).
\medskip

For technical reasons we also need  the following uniform ellipticity
condition:

\medskip
\noindent
{\bf Condition B.} There exists
a positive number~$\kappa>1$ such that
\[
\kappa^{-1} \leq \omega^+_0 \leq \kappa, \qquad
\kappa^{-1} \leq \omega^-_0 \leq \kappa \quad  \mbox{ $\IP$-a.s.}
\]

\medskip
Given the realization of the random environment~$\omega$, define the potential function
for $x\in\mathbb{R_+}$, by
\[
V(x) = \left\{
    \begin{array}{ll}
     \sum_{i=0}^{[x]-1}\ln\frac{\omega^-_i}{\omega^+_i}, & x>0,\\
     0, & x=0
    \end{array}
           \right.
\]
where $[x]$ is the integer which is the closest to $x$.\\
We will focus on estimating 
the tail of the quenched distribution of $T_{\gam}$ for typical configurations $\omega$ of the environment. This enables us to
approximate~$V$ by Brownian motion, but it is most
convenient to use the well-known
Koml{\'o}s-Major-Tusn{\'a}dy~\cite{kmt} strong approximation
theorem.
Indeed, it allows us to relate the features of
long time behavior for the walk to Brownian functionals directly
built on the model, simplifying much the proof of limit properties:
in a possibly enlarged probability space
there exist a
%coupling
version of our
environment~ $\omega$ and a Brownian
motion $(W(x), x \in \mathbb{R_+})$ with diffusion constant  $\sigma$ (i.e.,
$\mathop {\mathrm {Var}}(W(x))=\sigma^2 x$), such that for some $\hat K > 0$
\begin{equation}
\label{kmth}
\IP \Big[ \limsup_{x \to +\infty} \frac{|V(x)-W(x)|}{ \ln x}
  \leq \hat K
\Big] =1.
\end{equation}
A useful consequence of the KMT strong approximation theorem is that if $x$ is not too far from the origin, then $V(x)$ and $W(x)$ are rather close for the vast majority of environments. Hence, it
is convenient to introduce the following set of ``good'' environments,
and to restrict our forthcoming
computations to this set.
Fix an arbitrary $M_0>0$; for any~$t>e$,
let
\begin{equation}
\label{blizko}
\Gamma_t = \Big\{\omega : |V(x)-W(x)|\leq K_0\ln\ln t\;, \;
%%\mbox{ for all }
 x\in [0,\ln^{M_0}t]\Big\},
\end{equation}
where we can choose $K_0 \in (0, \infty)$ in such a way that for $\IP$-almost all~$\omega$, it holds that $\omega\in\Gamma_t$
for all~$t$ large enough (cf. e.g. (13) of \cite{CP}).
Also, we will systematically use the following notation,
\begin{equation}
\tau_{A}(\xi_i) =\inf\{t>0;\xi_i(t)\in A\} \nonumber\\
\end{equation}
with $1\leq i\leq \gam$ and $A\subset \mathbb{Z_+}$.
\medskip

The next definition is central in our construction. Although it applies to more general numerical functions defined on $\Z_+$ or $\R_+$, we will use it only in the case of the Brownian motion $W$.
\begin{df}
\label{def_t_stable1}
Let $t>1$ and $m^* \in \R_+$. Denote $r=r(t,m^*)$ the smallest $x > m^*$ such that $W(x) \geq W(m^*)+ \ln t$
(observe that as $W$ is a Brownian motion $r$ exists $\IP$-a.s). We will say that $m^*$ is the first $t$-stable point to the right of the origin if $W(m^*)=\min_{x\in [0,r]}W(x)$ and $m^*$ is the smallest point on $\R_+$ verifying this condition. Then, we define the other $t$-stable points as follows. For $m\in \R_+$, denote by  $l = l(t,m)$
the largest $x<m$ such that $W(x) \geq W(m)+ \ln t$ (if such an $x$ exists). We say that $m$ is a $t$-stable point (different from the first $t$-stable point $m^*$) if $m\neq m^*$, $l(m,t)$ and $r(m,t)$ exist and
$m$ satisfies $W(m)=\min_{x\in [l,r]}W(x)$.
\end{df}

\noindent
{\it Warning:} Observe that, according to Definition~\ref{def_t_stable1},
the point~$m$ is not generally integer. Thus, throughout
this paper, the statement ``the random walk hits a $t$-stable point''
means that it hits the site~$x\in \mathbb{Z_+}$ which is closest to the
$t$-stable point. As a rule, real points $x \in \mathbb{R_+}$
will be replaced, if the context requires, with the closest integer, that
we may still denote by the same symbol $x$, if no confusion can occur.

%%%%%%%%%%%%%%%%%%%%%%%%%%%%%%%%%%%
\subsection{Results} 
Let ${\hat m}_1(t), {\hat m}_2(t),\dots, {\hat m}_{\gam}(t)$ be the first $\gam$ $t$-stable points to the right of the origin and for $1\leq i \leq \gam-1$
\[{\hat h}_i(t)=\argmax_{x\in ({\hat m}_i(t),{\hat m}_{i+1}(t))}W(x).
\] 
\\
We define the process $\zeta_\gam=(\zeta_\gam(t))_{t>1}$ by
\begin{equation}
\label{Procze}
\zeta_{\gam}(t):=\sum_{i=1}^{\gam-1}(\gam-i) \frac{W({\hat h}_i(t))-W({\hat m}_i(t))}{\ln t}.
\end{equation}
We obtain the following results.
\begin{theo}
\label{theo1}
With the process $\zeta_{\gam}$ defined above it holds that
\[
\lim_{t \rightarrow \infty}  \Big|\frac{\ln \Po[T_{\gam}>t]}{\ln t}+\zeta_{\gam}(t)\Big|=0 \phantom{**}\mbox{in $\IP$-probability}.
\] 
\end{theo}
Then we set out results about the distribution of $\zeta_{\gam}(t)$.
\begin{theo}
\label{theo3}
The distribution of $\zeta_{\gam}(t)$ does not depend on $t$. Moreover, for each $t>1$, the random variable $\zeta_{\gam}(t)-\frac{\gam(\gam-1)}{2}$ has density $f_{\gam}$ given by
\begin{equation}
\label{densityK}
f_\gam(x)=\sum_{i=1}^{\gam-1}\frac{(-1)^{\gam-1-i}i^{\gam-2}}{i!(\gam-1-i)!}e^{-x/i}
\end{equation}
for $x\geq 0$.
\end{theo}
We finish this section by formulating the result about the quenched moments of the meeting time.
\begin{theo}
\label{theo2}
Let $c$ be a positive number, we have $\IP$-a.s.
\[
 \Eo[T_\gam^{c}]<\infty  \quad \mbox{if $c<\frac{\gam(\gam-1)}{2}$} 
\]  and 
\[
 \Eo[T_\gam^{c}]=\infty  \quad \mbox{if $c >\frac{\gam(\gam-1)}{2}$}
 . \]
\end{theo}
\medskip

\begin{figure}[!htb]
\begin{center}
\includegraphics[scale= 0.9]{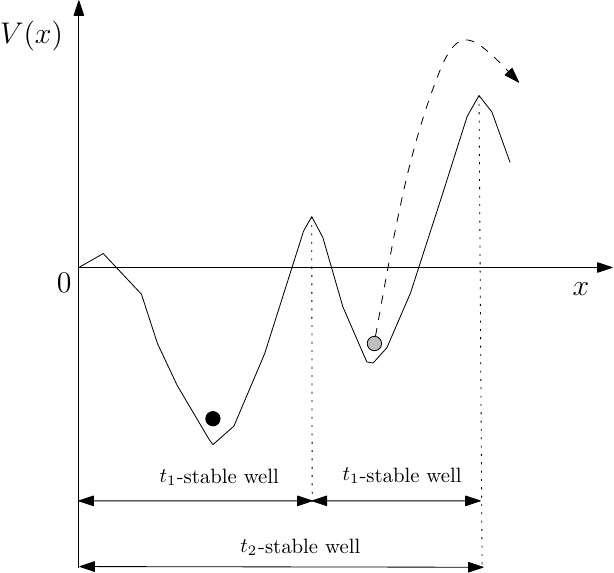}
\caption{Example of a strategy for two random walks on $\Z_+$.}
\label{figintro1}
\end{center}
\end{figure}

At this point let us explain why we consider reflected random walks on $\Z_+$ and not random walks on $\Z$. Despite the main ideas of our proofs should work for random walks on $\Z$, it is much more difficult to write them because the number of strategies for the random walks not to meet until time $t$ is much greater in the case of $\Z$ than in the case of $\Z_+$. We call a strategy the typical behavior of the random walks until time $t$. For example, on Figure \ref{figintro1} we represented two random walks on $\Z_+$ by black and grey particles.  Imagine that at some order of time smaller than $t_1$, the particles are in two different $t_1$-stable wells but in the same $t_2$-stable well with $t_1<t_2<t$.
\medskip

\begin{figure}[!htb]
\begin{center}
\includegraphics[scale= 0.9]{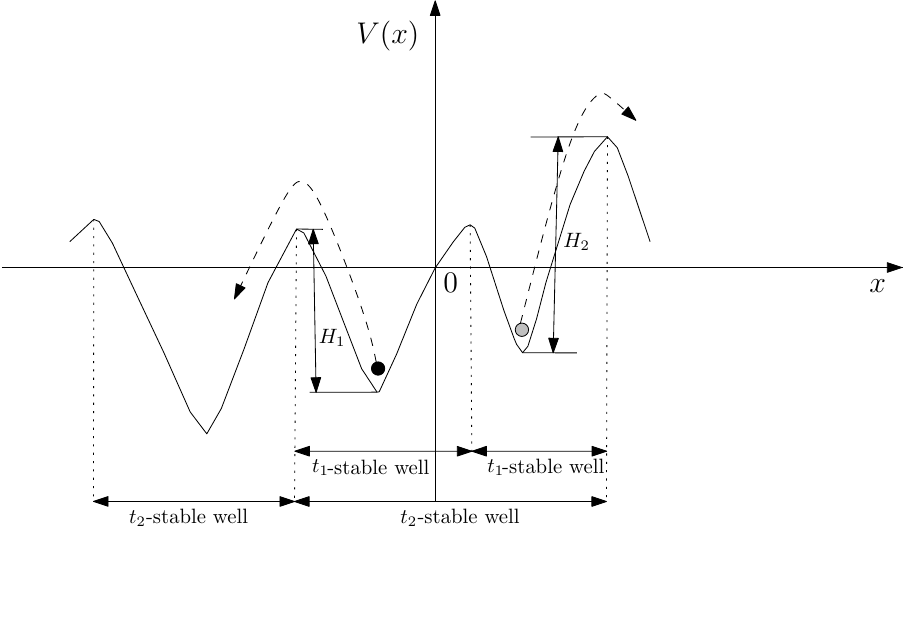}
\caption{Example of two possible strategies for two random walks on $\Z$. $H_1<H_2$.}
\label{figintro2}
\end{center}
\end{figure}

\noindent
For the two particles not to meet until time $t_2$, we need to expel the grey one to the right into an other $t_2$-stable well. More generally, in the reflected case, we will see in the proof of Theorem \ref{theo2} that at each order of time smaller than $t$, we have to expel the grey particle to the right into a sufficiently deep well until the two particles are in the first two $t$-stable wells (indeed, these two $t$-stable wells are sufficiently deep to hold the particles until time $t$ with high probability). Therefore, in some sense the strategy  does not depend on the realization $\omega$ of our random environment since we always force the rightmost particle to go to the right. On $\Z$, this not the case anymore. The possible strategies for the random walks not to meet until time $t$ will highly depend on the realization $\omega$ of our random environment. On Figure \ref{figintro2}, we have two particles on $\Z$ which are in two different $t_1$-stable wells but in the same $t_2$-stable well. In this case, we could expel the black particle to the left into an other $t_2$-stable well instead of expelling the grey one to the right since the barrier of potential $H_1$ the black particle has to jump over is smaller than the barrier $H_2$ the grey particle has to jump over. Note, however, that we may want to force the grey particle out to the right because this would lead to a situation when on next steps the ``cost'' decreases. Thus, we have to know the exact topography of the potential $V$ (that is the sequence of stable wells at each order of time smaller than $t$) to determine the strategy which has the cheapest cost. This leads to a great number of strategies to reach time $t$, which leaves the problem on $\Z$ much more complicated.

In the proof of Theorem \ref{theo1}, another possible approach would be using the Karlin-McGregor formula (see \cite{Karlin}) for determining $\Po[T_{\gam}>t]$.  However, to apply this formula, one would have to obtain estimates on transition probabilities which are much finer then those obtained in \cite{CP}; with the existing estimates the leading term in the determinant in the Karlin-McGregor formula equals zero, which yields no nontrivial results. 

The rest of the paper is organized as follows.
In section \ref{goodenv}, we will caracterize ``good environments'' and prove Lemma \ref{l_good_omega} that will be essential in order to prove Theorem \ref{theo1}.
Section \ref{SEC2} is dedicated to the proofs of Theorems \ref{theo1} and \ref{theo2} in the case $\gam=2$. As we intend to use an inductive argument to prove Theorems \ref{theo1} and \ref{theo2} in the general case, the preliminar treatment of the case $\gam=2$ is necessary. Section \ref{SECK} is dedicated to the proofs of Theorems \ref{theo1} and \ref{theo2} in the general case. Finally, section \ref{SECtheo3} is dedicated to the proof of Theorem \ref{theo3}. 

We mention here that all the constants that will be defined do not depend on the realization $\omega$ of the random environment but they may depend on the law of the environment and $\gam$. Furthermore, we will denote by $K'_1$, $K'_2$,~$\dots$ the ``local'' constants, that is, those that are used only in a small neighbourhood of the place where they appear for the first time. In this case we restart the numeration at the beginning of each section.

\section{Definition of $t$-good trajectories}
\label{goodenv}
First of all, we need to give more definitions.
Observe that between two successive $t$-stable points $m$ and $m'$ there is always a point
\[
h=\argmax_{x\in (m,m')}W(x)
\]
 which is a maximum of $W$. Thus it is possible to separate two successive $t$-stable points of $W$. Denoting by ${\mathcal S}_t$ the set of all the $t$-stable points we define 
\[
{\mathcal H}_t:=\{h\geq 0 \phantom{*} \mbox{\it{such that there exist} } m, m'\in {\mathcal S}_t: h =\argmax_{x\in (m,m')}W(x)\}\cup\{0\}.
\]
\begin{df}
\label{def_t_stable2}
If $m$ is a $t$-stable point then we define the $t$-stable well of $m$ as
${\mathcal D}_t(m):=[\max {\mathcal H}_t\cap [0,m), \min {\mathcal H}_t\cap (m,\infty))$.
\end{df}
Furthermore, we have the property that $\IP$-a.s., for all $t>1$, $\cup_{m\in {\mathcal S}_t} {\mathcal D}_t(m)=\mathbb{R_+}$, where the union is disjoint.
We end this section by giving the definition of the elevation of a finite interval $[a,b]$.
\begin{df}
 The {\it elevation\/} of $I=[a,b]$
is defined as the Brownian functional
\begin{equation*}
\e(I) =  \max_{ x,y \in I} \max_{z\in [x,y]}
   \{W(z)-W(x)-W(y)+\min_{v\in I} W(v)\} .
\end{equation*}
where $[x,y]$ denotes the interval with endpoints $x,
y$ regardless of $x<y$ or $x>y$.
\end{df}
It can be seen that in the definition of~$\e(I)$
one may assume that $y$ is the global minimum of $W$ on $I$, $x$
is one of local minima,
and~$z$ is one of local maxima
of $W$ in~$I$. 
\medskip

In order to prove Theorems \ref{theo1} and \ref{theo2}, it will be convenient to exclude some particular trajectories of the Brownian motion $W$; that is why we will give in this section a formal definition of what we call a "$t$-good trajectory" of the Brownian motion $W$. We will see that our definition of a $t$-good trajectory is not very restrictive since for almost every trajectory, there exists $t_0$ (which depends on the trajectory of $W$) such that for all $t>t_0$ the trajectory is $t$-good. For the sake of brevity, we will not indicate in our notations the dependence on $t$ and $\omega$ when no confusion can occur.
\medskip

The following construction is rather technical, its usefulness will appear more clearly in sections \ref{SEC2} and \ref{SECK}.
First take a positive decreasing function $\alpha(t)$ such that $\lim_{t \rightarrow \infty } \alpha(t)=0$ (for the moment it is not necessary to explicit $\alpha$, but in Lemma \ref{l_good_omega} we will take $\alpha(t)= \ln^{-5/6}t$). Fix $t>1$ and consider the first $\gam$ $t^{\alpha(t)}$-stable points to the right of the origin. Let us denote by $m_1(1),\dots, m_{\gam-1}(1)$ and $m'_{\gam}(1)$ these $\gam$ $t^{\alpha(t)}$-stable points. Then let us define  for $1 \leq j \leq \gam-2$,
 \[h_j(1)=\arg\max_{x\in(m_j(1),m_{j+1}(1))}W(x)\] 
 and
 \[h_{\gam-1}(1)=\arg\max_{x\in(m_{\gam-1}(1),m'_{\gam}(1))}W(x).\]
Then, for $1\leq j\leq \gam-1$, define
\begin{equation*}
r_j(1)= \frac{W(h_j(1))-W(m_j(1))}{\ln t}
\end{equation*}
and for $1\leq j \leq \gam-2$,
\begin{equation*}
l_j(1)= \frac{W(h_j(1))-W(m_{j+1}(1))}{\ln t}
\end{equation*}
and
\begin{equation*}
l_{\gam-1}(1)= \frac{W(h_{\gam-1}(1))-W(m'_{\gam}(1))}{\ln t}.
\end{equation*}
Finally let,
\[a_1=\min_{j<\gam}\{r_j(1)\}\wedge \min_{j<\gam}\{ l_j(1)\}.
\]
We will use the following rule for $n\geq 1$.
\medskip

\noindent
{\sl{Case 1:} $a_n=l_{\gam-1}(n)$.}
\medskip

\begin{figure}[!htb]
\begin{center}
\includegraphics[scale= 0.7]{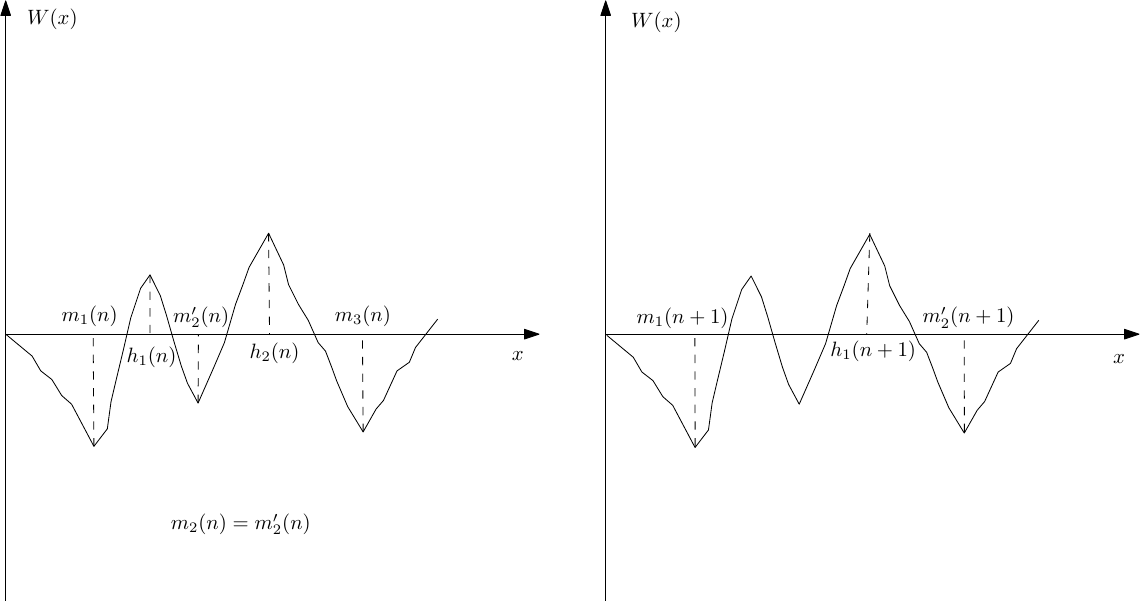}
\caption{Case $\gam=2$, $a_n=l_1(n)$, $m'_2(n)$ is $t^{a_n}$-stable.}
\label{fig1}
\end{center}
\end{figure}

If the point $m'_{\gam}(n)$ is $t^{a_n}$-stable then do the following. Rename $m'_{\gam}(n)$ as $m_{\gam}(n)$. Then consider the first $t^{a_n}$-stable point immediately after $m_\gam(n)$. Let us call it $m_{\gam+1}(n)$. See Figure \ref{fig1} (left picture).
\medskip

\noindent
Then define 
\[h_\gam(n)=\arg\max_{x\in(m_{\gam}(n),m_{\gam+1}(n))}W(x),\]
\[r_\gam(n)=\frac{W(h_\gam(n))-W(m_\gam(n))}{\ln t}\]
and
\[l_\gam(n)= \frac{W(h_\gam(n))-W(m_{\gam+1}(n))}{\ln t}.\]
Then, for $1\leq j\leq \gam-1$, define
\begin{equation*}
m_j(n+1)=m_j(n)
\end{equation*}
and
\begin{equation*}
m'_{\gam}(n+1)=m_{\gam+1}(n),
\end{equation*}
For $1\leq j\leq \gam-2$, define
\[h_j(n+1)=\arg\max_{x\in(m_j(n+1),m_{j+1}(n+1))}W(x)\]
and
\[h_{\gam-1}(n+1)=\arg\max_{x\in(m_{\gam-1}(n+1),m'_{\gam}(n+1))}W(x).\] 
See Figure \ref{fig1}.
For $1\leq j \leq \gam-1$, we define
\begin{equation*}
r_j(n+1)= \frac{W(h_j(n+1))-W(m_j(n+1))}{\ln t}
\end{equation*}
and for $1\leq j \leq \gam-2$
\begin{equation*}
l_j(n+1)= \frac{W(h_j(n+1))-W(m_{j+1}(n+1))}{\ln t}
\end{equation*}
and
\begin{equation*}
l_{\gam-1}(n+1)= \frac{W(h_{\gam-1}(n+1))-W(m'_{\gam}(n+1))}{\ln t}.
\end{equation*}
Note that
for $1\leq j\leq \gam-2$,
\begin{align*}
r_j(n+1)&=r_j(n)\nonumber\\
l_j(n+1)&=l_j(n)
\end{align*}
and
\begin{align*}
r_{\gam-1}(n+1)&=r_{\gam-1}(n)-l_{\gam-1}(n)+r_\gam(n)\nonumber\\
l_{\gam-1}(n+1)&=l_\gam(n).
\end{align*}
Finally define

\[a_{n+1}=\min_{j<\gam}\{r_j(n+1)\}\wedge \min_{j<\gam}\{ l_j(n+1)\}. 
\]

\begin{figure}[!htb]
\begin{center}
\includegraphics[scale= 0.7]{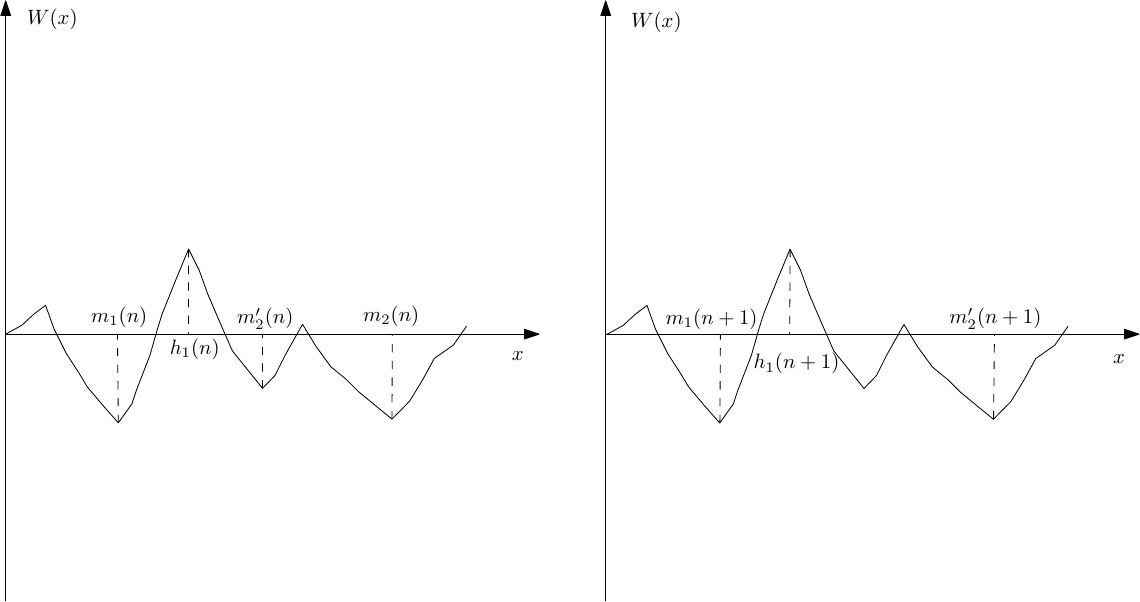}
\caption{Case $\gam=2$, $a_n=l_1(n)$, $m'_2(n)$ is not $t^{a_n}$-stable.}
\label{fig2}
\end{center}
\end{figure}

If the point $m'_\gam(n)$ is not $t^{a_n}$-stable, nevertheless it belongs to a $t^{a_n}$-stable well.
Then there exists a $t^{a_n}$-stable point $x_r$ such that $x_r>m'_\gam(n)$ and such that $m'_\gam(n) \in {\mathcal D}_{t^{a_n}} (x_r)$. Then rename $x_r$ as $m_{\gam}(n)$. See Figure \ref{fig2} (left picture).
In this case we define
\begin{equation*}
h^*(n)=\arg\max_{x\in(m'_{\gam}(n),m_{\gam}(n))}W(x)
\end{equation*}
and
\begin{equation*}
r_{\gam}(n)=\frac{W(h^*(n))-W(m'_{\gam}(n))}{\ln t}.
\end{equation*}
\medskip

\medskip

\noindent
Next, for $1\leq j\leq \gam-1$, we define
\begin{equation*}
m_j(n+1)=m_j(n)
\end{equation*}
and
\begin{equation*}
m'_{\gam}(n+1)=m_{\gam}(n).
\end{equation*}
Define for $1\leq j\leq \gam-2$,
\begin{equation}
\label{Cons1}
h_j(n+1)=\arg\max_{x\in(m_j(n+1),m_{j+1}(n+1))}W(x)
\end{equation}
and
\begin{equation}
\label{Cons2}
h_{\gam-1}(n+1)=\arg\max_{x\in(m_{\gam-1}(n+1),m'_{\gam}(n+1))}W(x).
\end{equation}
See Figure \ref{fig2}.
For $1\leq j \leq \gam-1$,
\begin{equation}
\label{Cons3}
r_j(n+1)= \frac{W(h_j(n+1))-W(m_j(n+1))}{\ln t}
\end{equation}
and for $1\leq j \leq \gam-2$,
\begin{equation}
\label{Cons4}
l_j(n+1)= \frac{W(h_j(n+1))-W(m_{j+1}(n+1))}{\ln t}
\end{equation}
and
\begin{equation}
\label{Cons5}
l_{\gam-1}(n+1)= \frac{W(h_{\gam-1}(n+1))-W(m'_{\gam}(n+1))}{\ln t}.
\end{equation}
Note that in this case
for $j\leq \gam-2$,
\begin{align*}
r_j(n+1)&=r_j(n)\nonumber\\
l_j(n+1)&=l_j(n)
\end{align*}
and
\begin{align*}
r_{\gam-1}(n+1)&=r_{\gam-1}(n)\nonumber\\
l_{\gam-1}(n+1)&=\frac{W(h_{\gam-1}(n))-W(m_\gam(n))}{\ln t}.
\end{align*}
Finally define

\[a_{n+1}=\min_{j<\gam}\{r_j(n+1)\}\wedge \min_{j<\gam}\{ l_j(n+1)\}.
\]
\begin{figure}[!htb]
\begin{center}
\includegraphics[scale= 0.7]{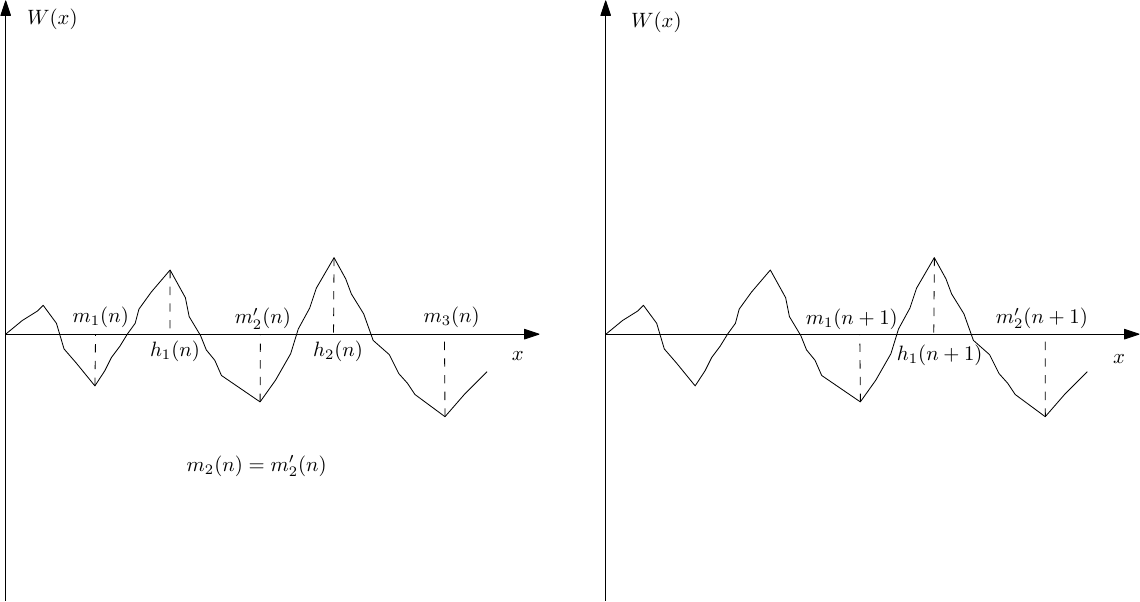}
\caption{Case $\gam=2$, $a_n=r_1(n)$, $m'_2(n)$ is $t^{a_n}$-stable.}
\label{fig3}
\end{center}
\end{figure}
\medskip

\medskip

\medskip

\medskip

\medskip

\medskip

\medskip

\medskip

\medskip

\medskip

\medskip

\medskip

\medskip

\medskip

\medskip

\medskip

\medskip

\noindent
{\sl{Case 2:} $a_n\neq l_{\gam-1}(n)$.}
\medskip

\begin{figure}[!htb]
\begin{center}
\includegraphics[scale= 0.7]{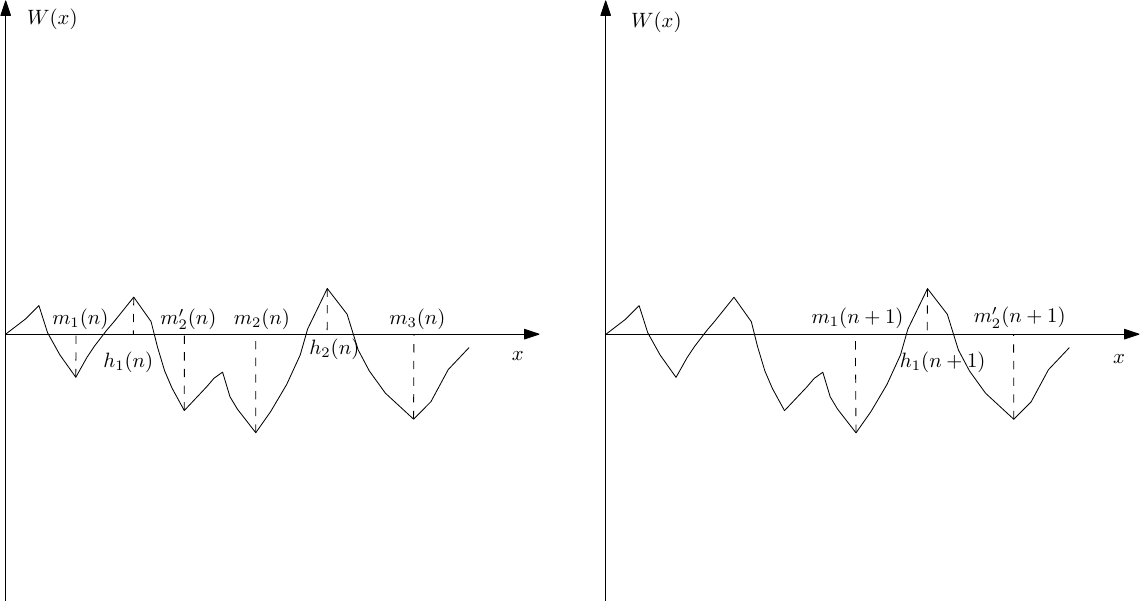}
\caption{Case $\gam=2$, $a_n=r_1(n)$, $m'_2(n)$ is not $t^{a_n}$-stable.}
\label{fig4}
\end{center}
\end{figure}

If the point $m'_{\gam}(n)$ is $t^{a_n}$-stable then do the following. Rename $m'_{\gam}(n)$ as $m_{\gam}(n)$. Then consider the first $t^{a_n}$-stable point immediately after $m_\gam(n)$. Let us call it $m_{\gam+1}(n)$. See Figure \ref{fig3} (left picture).
Then define 
\[h_\gam(n)=\arg\max_{x\in(m_{\gam}(n),m_{\gam+1}(n))}W(x),\]
\[r_\gam(n)=\frac{W(h_\gam(n))-W(m_\gam(n))}{\ln t}\]
and
\[l_\gam(n)= \frac{W(h_\gam(n))-W(m_{\gam+1}(n))}{\ln t}.\]
If the point $m'_\gam(n)$ is not $t^{a_n}$-stable, nevertheless it belongs to a $t^{a_n}$-stable well.
Hence there exists a $t^{a_n}$-stable point $x_r$ such that $x_r>m'_\gam(n)$ and such that $m'_\gam(n) \in {\mathcal D}_{t^{a_n}} (x_r)$ then consider the first  $t^{a_n}$-stable point immediately after $x_r$. Call it $m_{\gam+1}(n)$. Let us also rename $x_r$ as $m_\gam(n)$. See Figure \ref{fig4} (left picture).
Then let us define 
\[h_\gam(n):=\arg\max_{x\in(m_\gam(n),m_{\gam+1}(n))}W(x)\]
 and
\[r_\gam(n):= \frac{W(h_{\gam}(n))-W(m_{\gam}(n))}{\ln t}.\]
Next, we have to distinguish two cases.

If $a_n=r_i(n)$ for some $1\leq i\leq \gam-1$, we define the following.\\
For $j<i$,
\begin{align*}
m_j(n+1)=m_{j}(n).
\end{align*}
For $i\leq j \leq \gam-1$
\begin{align*}
m_j(n+1)=m_{j+1}(n)
\end{align*}
and 
\begin{align*}
m'_\gam(n+1)=m_{\gam+1}(n).
\end{align*}
Then for $1\leq j \leq \gam-1$, define $h_j(n+1)$, $ r_j(n+1)$ and $ l_j(n+1)$ as in (\ref{Cons1}), (\ref{Cons2}), (\ref{Cons3}), (\ref{Cons4}) and (\ref{Cons5}). See Figure \ref{fig3}.\\
Observe that in this case for $j<i-1$,
\begin{align*}
l_j(n+1)=l_{j}(n),
\end{align*}
\begin{align*}
l_{i-1}(n+1)=l_{i-1}(n)-r_i(n)+l_i(n),
\end{align*}
for $i\leq j \leq \gam-3$,
\begin{align*}
l_j(n+1)=l_{j+1}(n)
\end{align*}
and 
\begin{align*}
l_{\gam-2}(n+1)=\frac{W(h_{\gam-1}(n))-W(m_\gam(n))}{\ln t},\nonumber\\
l_{\gam-1}(n+1)=\frac{W(h_{\gam}(n))-W(m_{\gam+1}(n))}{\ln t}.
\end{align*}
Besides for $j<i$,
\begin{equation*}
r_j(n+1)=r_{j}(n)
\end{equation*}
and for $i \leq j \leq \gam-1$
\begin{equation*}
r_j(n+1)=r_{j+1}(n).
\end{equation*}
\medskip

If $a_n=l_i(n)$ for some $1\leq i\leq \gam-2$, we define the following quantities.
For $j\leq i$,
\begin{align*}
m_j(n+1)=m_{j}(n)
\end{align*}
and for $i< j \leq \gam-1$,
\begin{align*}
m_j(n+1)=m_{j+1}(n)
\end{align*}
and 
\begin{align*}
m'_\gam(n+1)=m_{\gam+1}(n).
\end{align*}
Then for $1\leq j\leq \gam-1$, define $h_j(n+1)$, $ r_j(n+1)$ and $ l_j(n+1)$ as in (\ref{Cons1}), (\ref{Cons2}), (\ref{Cons3}), (\ref{Cons4}) and (\ref{Cons5}). See Figure \ref{fig4}.\\
Observe that in this case, for $j< i$,
\begin{equation*}
r_j(n+1)=r_j(n),
\end{equation*}
\begin{equation*}
r_i(n+1)=r_i(n)-l_i(n)+r_{i+1}(n),
\end{equation*}
and for $i<j\leq \gam-1$,
\begin{equation*}
r_j(n+1)=r_{j+1}(n).
\end{equation*}
For $j< i$,
\begin{equation*}
l_j(n+1)=l_{j}(n),
\end{equation*}
for $i\leq j\leq \gam-3$,
\begin{equation*}
l_j(n+1)=l_{j+1}(n)
\end{equation*}
and
\begin{align*}
l_{\gam-2}(n+1)=\frac{W(h_{\gam-1}(n))-W(m_\gam(n))}{\ln t},\nonumber\\
l_{\gam-1}(n+1)=\frac{W(h_{\gam}(n))-W(m_{\gam+1}(n))}{\ln t}.
\end{align*}
Finally define
\[a_{n+1}=\min_{j<\gam}\{r_j(n+1)\}\wedge \min_{j<\gam}\{ l_j(n+1)\}.
\]

\medskip

\medskip

Thus by construction, we obtain an increasing sequence $(a_n)_{n\geq1}$, $\IP$-a.s. In the sequel, we will even see that the sequence reaches 1 very quickly. Hence, heuristically, the sequence $(a_n)_{n\geq1}$ will allow us to approximate time $t$ by the sequence $(t^{a_n})_{n\geq1}$. For each order of time $t^{a_n}$ with $a_n\leq1$, we will define a strategy for the random walks not to meet during this order of time. This argument will be developed in sections \ref{SEC2} for 2 random walks and \ref{SECK} for the general case.
Furthermore, observe that by construction in all the cases we have that
\begin{align}
\label{ELEVPROP}
\lefteqn{\max_{1\leq j \leq \gam-1}\{\e[m_j(n),h_j(n)] \vee \e[h_j(n),m_{j+1}(n)]\}}\phantom{*****}\nonumber\\
&\vee \{\e[m_{\gam-1}(n),h_{\gam-1}(n)] \vee \e[h_{\gam-1}(n),m'_{\gam}(n)]\} \leq  a_{n-1}\ln t,
\end{align}
for $n\geq 1$, with the convention $a_0=\alpha$. This last property will turn out to be important in the proofs of Theorems \ref{theo1} and \ref{theo2}.
Now, let us introduce $N=N(\omega,t)=\min \{n\geq 1; a_{n}\geq 1 \}$. As, by construction, the points $m_j(N)$, $1\leq j\leq \gam$ are the first $t$-stable points to the right of the origin, we can give a new definition of our process $\zeta_\gam$ of Theorem \ref{theo1} in function of the $r_i(N)$-s with $1\leq i\leq \gam-1$. Observe that we have,
\begin{equation}
\label{PROC}
\zeta_\gam(t)=\sum_{i=1}^{\gam-1}(\gam-i)r_i(N) .
\end{equation}
Now take $0<\eps(t)<\alpha(t)$ and for $1\leq k\leq \gam$, let $J^1_k(n)$ the number of $t^{a_n-\eps}$-stable wells in the interval $[m_k(n),h_k(n)]$ and  $J^2_k(n)$ the number of  $t^{a_n-\eps}$-stable wells in $[h_k(n),m_{k+1}(n)]$ (if $h_{\gam}(n)$ is not defined we pose $J^1_\gam(n)=0$ and $J^2_\gam(n)=0$). We define the set of $t$-good trajectories of $W$ in the following way. For the sake of brevity, we will use the notations $\ln_2 t:=\ln\ln t$ and $\ln_3 t:=\ln\ln\ln t$. 
\medskip

\begin{df}
\label{goodenviron}
Fix $t>e^e$. Let $A$ and $B$ be two constants (they do not depend on $\omega$). A realization $\omega$ of the Brownian motion $W$ is called $t$-good, if 
\begin{itemize}
\item[(i)] $N=N(\omega,t)\leq A\ln_2 t$;
\item[(ii)] $m'_\gam(N)\leq B \ln_3 t \ln^2 t$;
\item[(iii)] $m'_\gam(1)\leq \ln^{1/2} t$;
\item[(iv)]  $\sum_{k=1}^{\gam}J^1_k(i)+J^2_k(i) \leq \ln_2 t$ for $1\leq i \leq N$;
\item[(v)] $\max\{|W(x)|; x \in [0, B \ln_3 t \ln^2 t]\}  \leq \ln_3 t \ln t$.
\end{itemize}
\end{df}

The following lemma plays a crucial role.

\begin{lm}
\label{l_good_omega}
Take $\alpha(t)= \ln^{-5/6}t$ and $\eps=\eps(t)=\ln^{-11/12}t$ and let $\Lambda_{t}$ be the set of $\omega$ which are $t$-good. We can choose $A$ and $B$ sufficiently large such that $\IP[\lim_{t \rightarrow \infty}\Lambda_{t}]=1$. 
\end{lm}
\noindent
{\it Proof.} We start by showing that we can choose a constant $A$ sufficiently large such that there exists a set $\Omega'$ of probability 1 such that for all $\omega \in \Omega'$ there exists $t_0(\omega)$ such that for all $t>t_0(\omega)$ the item (i) of Definition \ref{goodenviron} is true.  First take $\delta>0$ and define the sequence of times $t_n=e^{(1+\delta)^n}$ for $n\geq 0$. Then we define the family of events
\begin{equation}
G_{n}=\{\omega; N(\omega,t_n)>A_1\ln_2 t_n\} \nonumber\\
\end{equation}\\
for $A_1$ a positive constant. We now argue by contradiction that
\[\IP[\limsup_{n\rightarrow \infty}G_{n}]=0
 \]
for  $A_1$ sufficiently large. Suppose that $\IP[\limsup_{n\rightarrow \infty}G_{n}]>0$. Then there exists $\omega \in \limsup_{n \rightarrow \infty}G_{n}$ such  that there exists a subsequence of times $(t_{n_l})_{l\geq1}$ such that 
\begin{equation}
\label{UF}
N(\omega,t_{n_l})> A_1 \ln_2 t_{n_l}
\end{equation}
for $l\geq1$.
\medskip

By construction of our $t$-good environments we know that the sequence $(a_i)_{i\geq 1}$ is increasing. As $W$ is a Brownian motion, we can even say that this sequence is strictly increasing $\IP$-a.s.. Nevertheless it is not clear how fast it grows. After a deeper look at the construction and using Lemma 6.1 of \cite{CP}, we can check that the sequence $(u_i)_{i\geq 1}=(a_{i\gam})_{i\geq 1}$ verifies the following properties. For each $t>1$ there exists a family of positive continuous i.i.d.\ random variables $(\vphi_i(t))_{i\geq1}$ such that
\begin{equation}
\label{Blub}
u_{i+1}(t)\geq (1+\vphi_i(t)) u_{i}(t)
\end{equation}
for $i\geq 1$.
Using item~(ii) of Lemma 6.1 of \cite{CP}, we can even show that the random variables  $(\vphi_i(t))_{i\geq1}$ have exponential moments.
Moreover, observe that the law of $(\vphi_i(t))_{i\geq1}$ does not depend on $t$.
Iterating (\ref{Blub}) and taking the logarithm of both sides, we obtain
\begin{equation}
\ln u_{i+1}(t)\geq \sum_{j=1}^i\ln (1+\vphi_j(t)) + \ln u_1(t)\nonumber\\
\end{equation}\\
for $i\geq 1$.
Then defining $N'=\max\{i\geq 1;u_i<1\}$ and using the facts that $u_{N'}<1$ and $u_1\geq \alpha(t)$ we obtain
\begin{equation}
\label{Blub2}
\frac{5}{6}\ln_2 t \geq \sum_{j=1}^{N'-1}\ln (1+\vphi_j(t)).
\end{equation}\\
Now using (\ref{Blub2}), (\ref{UF}) and the fact that $\gam(N'+1)\geq N$, we obtain that  for $\omega \in \limsup G_{n} $
\begin{equation}
\label{contrad}
\frac{5}{6}\ln_2 t_{n_l} \geq \sum_{j=1}^{\lfloor (A_1/2\gam) \ln_2 t_{n_l}\rfloor}\ln (1+\vphi_j(t_{n_l}))
\end{equation}\\
for all sufficiently large $l$.
\\
On the other hand, by Cram\'er theorem and the Borel-Catelli lemma, we can deduce that
\begin{equation}
\frac{1}{(A_1/2\gam) \ln_2 t_{n_l}} \sum_{j=1}^{\lfloor (A_1/2\gam) \ln_2 t_{n_l}\rfloor}\ln (1+\vphi_j(t_{n_l}))\longrightarrow  \IE[\ln (1+\vphi_{1}(t_1))]\nonumber\\ 
\end{equation}\\
$\IP$-a.s., as $l \rightarrow \infty$ with $\IE[\ln (1+\vphi_{1}(t_1))]$ a positive finite constant.  
Combining this last result with (\ref{contrad}) we obtain for $\IP$-a.a.\ $\omega \in \limsup_{n\rightarrow \infty}G_{n}$, 

\begin{equation}
\frac{5}{6}\ln_2 t_{n_l} \geq  \frac{A_1}{4\gam}\IE[\ln (1+\vphi_{1}(t_1))] \ln_2 t_{n_l}
\end{equation}\\ 
for all $l$ large enough.
\\
Taking $A_1>\frac{10\gam}{3\IE[\ln (1+\vphi_{1}(t_1))]}$ this last relation is clearly impossible,  and we showed that $\IP[\limsup_{n \rightarrow \infty}G_{n}]=0$.
\medskip

\noindent
Now let us define the sequence of intervals $I_n=[t_n,t_{n+1})$ for $n\geq 0$. It remains to show that we can find a new constant $A$ sufficiently large such that $\IP$-a.s.,  $N(\omega,t)\leq A \ln\ln t $ for $t $ sufficiently large. 
This is easy if we observe that for $t \in I_n$,  $N(\omega,t)\leq N(\omega,t_n)+N(\omega,t_{n+1})$ (this follows from the simple fact that that for $t \in I_n$, $[\alpha(t)\ln t,\ln t]\subset [\alpha(t_n)\ln t_n,\ln t_n]\cup[\alpha(t_{n+1})\ln t_{n+1},\ln t_{n+1}]$).  This implies together with $\IP[\limsup_{n \rightarrow \infty}G_{n}]=0$ that $\IP$-a.s., for every $t \in I_n$ $N(\omega,t)\leq 2A_1 \ln\ln t_{n+1}$ for all $n$ sufficiently large. With this last observation it is sufficient to take $A\geq3A_1$ to have that  (i) of Definition \ref{goodenviron} holds for $\IP$-a.a.\ $\omega$ for all $t$ sufficiently large.
\medskip

We continue the proof of Lemma \ref{l_good_omega} by showing that item (ii) of Definition \ref{goodenviron} holds for $\IP$-a.a.\ $\omega$ for all $t$ sufficiently large.
Again consider the sequence of intervals $I_n=[t_n,t_{n+1})$ for $n\geq 0$ and define the events
\begin{equation*}
D_n=\{\mbox{there exist $\gam$ $t_{n+1}$-stable wells in the interval $[0,\ln^2 t_n]$}\}
\end{equation*}
for $n\geq 0$.
\\
By the scaling property of the Brownian motion the probability of $D_n$ is independent of $n$. Furthermore observe that 
\begin{align*}
D_0&\supset \bigcap_{j=1}^{2\gam-1}\Big\{W\Big(\frac{j}{2\gam}\Big)-W\Big(\frac{j-1}{2\gam}\Big)<-3,  W\Big(\frac{j+1}{2\gam}\Big)-W\Big(\frac{j}{2\gam}\Big)>3\Big\},\nonumber\\
\end{align*}
for $\delta<2$.
As this last event has positive probability we obtain that there exists a positive constant $C$ such that  $\IP[D_n]\geq C$ for $n\geq 0$ and $\delta<2$.
Then consider the events
\begin{equation*}
E_n=\{\mbox{there exist $\gam$ $t_{n+1}$-stable wells in the interval $[0,M\ln n \ln^2 t_n]$}\}
\end{equation*}\\ 
for $M>(-\ln (1-C))^{-1}$ and $n\geq 1$.
By the Markov property and the fact that $\IP[D_n]\geq C$  we obtain that
\begin{equation*}
\IP[E_n]\geq 1-n^{-M\ln \frac{1}{1-C}}
\end{equation*}\\ 
for $n\geq 1$. 
Using the Borel-Cantelli lemma, we obtain that $\IP$-a.s.\ there exists $n_0=n_0(\omega)$ such that for every $n>n_0$ there exist $\gam$ $t_{n+1}$-stable wells in the interval $[0,M\ln n \ln^2 t_n]$.
Now observe that for $t\in I_n$ we have the two following inequalities
\begin{equation}
\label{buz1}
(1+\delta)^n \leq \ln t \leq (1+\delta)^{n+1}
\end{equation}
and
\begin{equation}
\label{buz2}
L_{\omega}(t_{n+1},M\ln n \ln^2 t_n )\leq L_{\omega}(t,M\ln n \ln^2 t_n ) \leq L_{\omega}(t_{n}, M\ln n \ln^2 t_n)
\end{equation}\\ 
where $L_{\omega}(x,y)$ is the number of $x$-stable wells in the interval $[0,y]$. Inequality (\ref{buz2}) follows from the fact that ${\mathcal S}_{t_{n+1}}\subset {\mathcal S}_t \subset {\mathcal S}_{t_n}$. 
Let us take $t \in I_n$ for $n>n_0$. Then we know that there exist $\gam$ $t_{n+1}$-stable wells in the interval $[0, M\ln n \ln^2 t_n]$. By inequality (\ref{buz2}) we know that there exist $\gam$ $t$-stable wells in the interval $[0, M\ln n \ln^2 t_n]$. And by inequality (\ref{buz1}) we know that there exists $\gam$ $t$-stable wells in the interval \[[0, M[\ln_3 t - \ln_2 (1+\delta)]\ln^2 t].\] Finally, taking $\delta=e-1$ we showed that (ii) of Definition \ref{goodenviron} holds for $\IP$-a.a.\ $\omega$ for all $t$ sufficiently large. The proof of (iii) is completely equivalent to the proof of (ii).
\medskip

We will now show that item (iv) of Definition \ref{goodenviron} holds for $\IP$-a.a.\ $\omega$ for all $t$ sufficiently large.
In fact, it is sufficient to show the following statement. Fix $t>e$ and consider the first $\gam$ $t$-stable points to the right of the origin. Let us call the last one ${\hat m}_{\gam}$.  We will show that $\IP$-a.s., for $t$ sufficiently large $L_{\omega}(t^{1-\eps}, {\hat m}_\gam) \leq \ln_2 t$.\\
Note that as $\eps(t) \rightarrow 0$ as $t \rightarrow \infty$ we have for $t$ sufficiently large that $L_{\omega}(t^{1-\eps}, {\hat m}_\gam)\leq L_{\omega}(\frac{t}{2},{\hat m}_\gam)$.
Then, consider the sequence of intervals $I_n$ with $\delta <2$. By the proof of (ii), we know that there exists a positive constant $M$ such that $\IP$-a.s.\ there exists $n_0=n_0(\omega)$ such that for all $n>n_0$ there exists $\gam$ $t_{n}$-stable wells in the interval $[0,M\ln n \ln^2 t_n]$. 

Now let $\psi_i$ be the width of the $i$-th  $\frac{t_{n}}{2}$-stable well divided by $(1+\delta)^{2n}$. By Lemma 2.1 of \cite{CP2} we know that the variables $(\psi_i)_{i\geq 2}$ are i.i.d.\ with exponential tail. It follows from Cram\'er theorem that there exists a positive finite constant $M_1$ such that 
\[
\IP\left[
\sum_{i=1}^{n} \psi_i < M \ln n\right] \leq  \exp\{-M_1 n\}\;,
\]
for all $n> 1$.
Therefore, by the Borel-Cantelli lemma we have that $\IP$-a.s.\ there exists $n_1=n_1(\omega)$ such that for every $n>n_1$ we have $L_{\omega}(\frac{t_n}{2},{\hat m}_\gam)\leq n$.  Finally, for $t \in I_n$, by (\ref{buz2}) we have that $L_{\omega}(\frac{t}{2},{\hat m}_\gam)\leq n$ and by (\ref{buz1}) $L_{\omega}(\frac{t}{2},{\hat m}_\gam)\leq \frac{\ln_2 t}{\ln (1+\delta)}$. Taking $\delta=e-1$ we showed that (iv) of Definition \ref{goodenviron} holds for $\IP$-a.a.\ $\omega$ for all $t$ sufficiently large. .
\medskip

We end the proof of Lemma \ref{l_good_omega} by showing that item (v) of Definition~\ref{goodenviron} holds for $\IP$-a.a. $\omega$ for all $t$ sufficiently large.
The law of the Iterated Logarithm for Brownian motion tells us that
\begin{equation}
\IP\Big[\limsup_{t \rightarrow \infty} \frac{|W(t)|}{\sqrt{2\sigma^{2} t\ln_2 t}}=1\Big]=1 \nonumber\\
\end{equation}
This implies that $\IP$-a.s., there exists $t_{1}(\omega)$ such that $|W(t)| \leq 2\sqrt{\sigma^{2}t\ln_2 t}$, $\forall t \geq t_{1}(\omega)$, which at its turn implies that there exists $t_{2}(\omega)>t_{1}(\omega)$ such that $\max_{x \in [0,t]} |W(x)| \leq 2\sqrt{\sigma^{2}t\ln_2 t}$, $\forall t \geq t_{2}(\omega)$.
And finally we deduce for $t$ sufficiently large,
\begin{equation*}
\max_{x \in [0,B\ln_3 t \ln^{2}t]} |W(x)| \leq \ln_3 t \ln t
\end{equation*}
which ends the proof of Lemma \ref{l_good_omega}. $\square$

\medskip

\noindent

%%%%%%%%%%%%%%%%%%%%%%%%%%%%%%%%%%%%%%%%
\section{Proofs of  Theorems \ref{theo1} and \ref{theo2} in the case $\gam=2$}
\label{SEC2}

%%%%%%%%%%%%%%%%%%%%%%%%%%%%%%%%%%%%
\subsection{Some auxiliary results}

We will need some results that we recall below.
\medskip

For any integers $a<x<b$, the probability for $\xi$ a random walk starting at $x$ to reach $b$ before $a$ is given by:
\begin{equation}
\label{exit}
\Po^x [ \tau_b(\xi) <  \tau_a(\xi)] = \frac{ \sum_{y=a+1}^x
  e^{V(y)-V(a)}}{ \sum_{y=a+1}^b  e^{V(y)-V(a)}},
\end{equation}
see e.g.\ Lemma~1 in \cite{Sinai}.
\medskip

We will also need the following upper bound on the probability of confinement which is derived from Proposition 4.1 of \cite{PGF}. Let $I=[a,b]$ with $0\leq a<b<\infty$ be an interval of $\Z_+$. Then define
\[
H_+(I)=\max_{x\in I}\Big(\max_{y\in[x,b]}V(y)-\min_{y\in[a,x)}V(y)\Big),
 \]
\[
H_-(I)=\max_{x\in I}\Big(\max_{y\in[a,x]}V(y)-\min_{y\in(x,b]}V(y)\Big)
 \]
and let $H(I)=H_+(I)\wedge H_-(I)$.
For $a<x<b$ we have
\begin{equation}
\label{LOWB}
\Po^x[\tau_{\{a,b\}}(\xi)>t] \leq \exp \Big\{-\frac{t}{C_1(b-a)^4e^{H(I)}}\Big\}
\end{equation}
with $C_1$ a positive constant.
\medskip

\noindent
The upper bound (\ref{LOWB}) can be unadapted, for our purpose, if we consider an interval $[a,b]=[m,m']$ where $m$ and $m'$ are two neighboring $t$-stable points (i.e. $(m,m')\cap  {\mathcal S}_t=\emptyset$). That is why we also need the following upper bound (see Lemma 3.1 of \cite{CP}).
Let $t>1$, and $I^+:=[ h, m'] $ and $I^-:=[m, h ] $ with $h=\argmax_{x\in (m,m')}W(x)$. Let
\[\Delta_1=m'-m,
\]
\[
\gamma_1=\max_{x\in[m,m']}V(x)-\min_{x\in [m,m']} V(x).
\]
For  all $x\in [m,m']$, it holds on~$\Gamma_t$
that
\begin{align}
\label{LOWB2}
\Po^x [ \tau_{\{m', m\}}(\xi) > t/k]
 &\leq  \exp\Big\{ - t^{\frac12 (1-\e(I^+)\ln^{-1}t)}
\Big(C_2(\Delta_1\ln^{2K_0}t)^{-1}\nonumber\\
&~~~~{}- C_3 e^{\gamma_1/2}\exp\Big\{
-\frac{\lambda(I^+)e^{\e(I^+)}t^{\frac12 (1-\e(I^+)\ln^{-1}t)}}
{2k} \Big\} \Big) \Big\}\nonumber\\
&~~{}+ \exp\Big\{ - t^{\frac12 (1-\e(I^-)\ln^{-1}t)}
\Big(C_2(\Delta_1\ln^{2K_0}t)^{-1}\nonumber \\
&~~~~{}- C_3 e^{\gamma_1/2}\exp\Big\{
-\frac{\lambda(I^-)e^{\e(I^-)}t^{\frac12 (1-\e(I^-)\ln^{-1}t)}}
{2k} \Big\} \Big) \Big\},
\end{align}
where $C_2$ and $C_3$ are positive constants, $\lambda(I^+)$ (respectively $\lambda(I^-)$) is the spectral gap of the reflected random walk on the interval $I^+$ (respectively~$I^-$) and $k$ is such that $\ln k=o(\ln t)$ as $t \rightarrow \infty$.
\medskip

We eventually need to estimate the cost of escaping a well to the right. Let again $I=[a,b]$, $m=\argmin_{x\in [a,b]}W(x)$ and suppose that $W(b)=\max_{x\in(m,b]}W(x)$. We will use the following estimate (see e.g.\ Lemma 3.4 in \cite{CP}), for any $s>0$,
\begin{equation}
\label{CER}
\Po^m[\tau_b(\xi)<s] \leq C_4se^{-V(b)+V(m)}
\end{equation}
with $C_4$ a positive constant.

%%%%%%%%%%%%%%%%%%%%%%%%%%%%%
\subsection{Upper bound for $P_{\omega}[T_2>t]$}
\label{Subsecup2}
To bound the probability distribution tail $P_{\omega}[T_2>t]$ from above the main idea is to define a sequence of increasing stopping times $(\sigma_i)_{i\geq 0}$ such that on each interval of the form $[\sigma_i, \sigma_{i+1}]$ we can find a simple strategy for the random walks not to meet. We formalize this argument as follows.
\medskip

\noindent
First, for technical reasons let us find an upper bound for $P_{\omega}[T_2>t^{1+\delta}]$ for $\delta>0$ arbitrary instead of $P_{\omega}[T_2>t]$.
Fix $t>1$ and define for $1\leq i \leq N$, the following stopping times
\begin{equation*}
\sigma_{i} =\inf\{s>0;\xi_2(s)=h_1(i)\} 
\end{equation*}
and the following events
\begin{equation*}
B_{i}=\{ T_2 \notin[\sigma_{i-1},\sigma_{i}]\} 
\end{equation*} 
for $1\leq i \leq N-1$ (with the convention $\sigma_{0}=0$)
and
\begin{equation*}
B_{N}=\{T_2 \notin[\sigma_{N-1},t^{1+\delta}]\}\cap \Big\{\sigma_{N-1}\leq \frac{t^{1+\delta}}{2}\Big\}. \nonumber\\ 
\end{equation*}
\medskip

To find the upper bound for $P_{\omega}[T_2>t^{1+\delta}]$, the following decomposition is the key of our analysis
\begin{align}
\label{d1}
\Po[T_2>t^{1+\delta}]
&=
\Po\Big[T_2>t^{1+\delta},\sigma_{N-1}\leq \frac{t^{1+\delta}}{2}\Big]+\Po\Big[T_2>t^{1+\delta},\sigma_{N-1}>\frac{t^{1+\delta}}{2}\Big].  \nonumber\\ 
\end{align}
As by (\ref{ELEVPROP}), $\e[0,h_1(N-1)]<\ln t$, by a similar argument as that we will use in subsection \ref{Riton} for the term $\Po[B_{1},\sigma_{1}\geq t^{b_1}]$, we can show that 
\begin{equation}
\label{DF}
\Po\Big[T_2>t^{1+\delta},\sigma_{N-1}>\frac{t^{1+\delta}}{2}\Big]\leq o(\exp(-\ln^8 t))
\end{equation}
as $t \rightarrow \infty$ and $\omega \in \Gamma_t\cap \Lambda_t$.
\medskip

\noindent
For the upper bound of the term $\Po[T_2>t,\sigma_{N-1}\leq t^{1+\delta}/2]$ we start by noting that
\begin{eqnarray}
\label{produit}
\lefteqn{\Po\Big[T_2>t^{1+\delta},\sigma_{N-1}\leq \frac{t^{1+\delta}}{2}\Big]}\phantom{******}\nonumber\\
&\leq& 
\Po[B_{1}\cap \ldots \cap B_{N}]
\nonumber\\
&=& \Po[B_{1}]\Po[B_{2}\mid B_{1}]\ldots \Po[B_{N}\mid B_{1}B_{2}\ldots B_{N-1}]. 
\end{eqnarray}
In the next subsections we will find upper bounds for the terms of the product of the right-hand side of (\ref{produit}). 
%%%%%%%%%%%%%%%%%%%%%%%%%%%%%
\subsubsection{Upper bound for $\Po[B_{1}]$}
\label{Riton}
Fix $t>1$ and let $\eps'=\eps'(t)=(28+8K_0)\ln_2t \ln^{-1}t$ where $K_0$ is from $(\ref{blizko})$. We bound from above the event $B_1$ in the following way,
\begin{equation}
\label{decomp1}
\Po[B_{1}]\leq \Po[B_{1},\sigma_{1}\geq t^{\alpha+\eps'}]+\Po[\sigma_{1}< t^{\alpha+\eps'}].
\end{equation}
From now on, for the sake of brevity we will denote $b_1=\alpha+\eps'$.
In the next two paragraphs we treat both terms of the right-hand side of (\ref{decomp1}).
%%%%%%%%%%%%%%%%%%%%%%%%%
\medskip

\noindent
{\bf Upper bound for $\Po[\sigma_{1}< t^{b_1}]$}
\medskip

\noindent
Using the Markov property, we obtain
\begin{align*}
\Po[\sigma_{1}< t^{b_1}]
&= \Po[\tau_{h_{1}(1)}(\xi_2)<t^{b_1}]\nonumber\\
&= \int_{0}^{t^{b_1}} P^{m_{1}(1)}_{\omega}[\tau_{h_{1}(1)}(\xi_2)<t^{b_1}-s] d\Po[\tau_{m_{1}(1)}(\xi_2)\leq s]
\nonumber\\
&\leq \Po^{m_{1}(1)}[\tau_{h_{1}(1)}(\xi_2)< t^{b_1}].
\end{align*}
Therefore, applying (\ref{CER}) to the last term, we obtain for $t$ sufficiently large and $\omega\in \Gamma_{t}\cap \Lambda_t$ 
\begin{equation}
\label{TP1}
\Po[\sigma_{1}< t^{b_1}]\leq C_4t^{\eps'} t^{-(r_{1}(1)-\alpha)}.
\end{equation}
{\bf Upper bound for $\Po[B_{1},\sigma_{1}\geq t^{b_1}]$}
\medskip

\noindent
We will show that this term is negligible in comparison with (\ref{TP1}) as $t \rightarrow \infty$ on the set of environments which belong to $\Gamma_{t} \cap \Lambda_{t}$. To show this, we will couple the two random walks with the random walks restricted to the interval $[0,h_{1}(1)]$.
So, we need to consider two other
processes, which are reflected versions of our
random walks in random environment: let
${\hat \xi}_1$ and ${\hat \xi}_2$ be the reflected RWRE on the interval $I_{1}=[0,h_{1}(1)]$.
The processes ${\hat \xi}_1$ and ${\hat \xi}_2$ have the same
jump rates as $\xi_1$ and $\xi_2$ on $[0,h_{1}(1)]$,
but jump from~$h_{1}$ to~$h_{1}-1$ at rate $\omega_{h_{1}(1)}^-$.
By construction we obtain
\begin{equation}
\label{COUPLE}
\Po[B_{1},\sigma_{1}\geq t^{b_1}]\leq \Po[{\hat \tau} > t^{b_1}],
\end{equation}
where ${\hat \tau}$ is the first meeting time of the two random walks ${\hat \xi}_1$ and ${\hat \xi}_2$.

The idea is now to use the spectral properties of the reflected random walks to find an upper bound for the right-hand side of (\ref{COUPLE}). We start by showing that, for $\omega \in \Gamma_{t} \cap \Lambda_{t}$, with probability at least $\ln^{-K}t$, with $K$ a positive constant, the two random walks restricted to the interval $I_{1}$ will meet in a time of order $t^{\alpha+\frac {\eps'}{2}}$. We will denote $b'_1=\alpha+\frac {\eps'}{2}$.
First observe that
\begin{equation}
\Po[{\hat \tau}\leq t^{b'_1}] \geq \Po\Big[{\hat \xi}_1(t^{b'_1})\in [m_{1}(1),h_{1}(1)],{\hat \xi}_2(t^{b'_1})\in [0,m_{1}(1)]\Big].\nonumber\\
\end{equation}
As the two random walks in fixed environment are independent we have
\begin{equation}
\Po[{\hat \tau}\leq t^{b'_1}] \geq \Po\Big[{\hat \xi}_1(t^{b'_1})\in [m_{1}(1),h_{1}(1)]\Big] \Po\Big[{\hat \xi}_2(t^{b'_1})\in [0,m_{1}(1)]\Big]\nonumber\\
\end{equation}
Then we can write
\begin{equation}
\Po\Big[{\hat \xi}_1(t^{b'_1})\in [m_{1}(1),h_{1}(1)]\Big]= \sum_{i=m_{1}(1)}^{h_{1}(1)} \Po\Big[{\hat \xi}_1(t^{b'_1}) = i \Big]\nonumber\\
\end{equation}
and apply Corollary 2.1.5 in \cite{SC} to obtain
\begin{align}
\label{SG1}
\lefteqn{\Po\Big[{\hat \xi}_1(t^{b'_1})\in [m_{1}(1),h_{1}(1)]\Big]}\phantom{******}\nonumber\\
&\geq \sum_{i=m_{1}(1)}^{h_{1}(1)} \mu^{I_{1}}(i) - \exp(-\lambda(I_{1})t^{b'_1}) \sum_{i=m_{1}(1)}^{h_{1}(1)} \Big( \frac {\mu^{I_{1}}(i)}{\mu^{I_{1}}(1)} \Big)^{\frac{1}{2}} 
\end{align}
where $\mu^{I_{1}}$ and $\lambda(I_{1})$ are respectively the invariant measure and the spectral gap of the reflected random walks ${\hat \xi}_1$ and ${\hat \xi}_2$.
At this point, let us define
\begin{equation*}
U_1:=\sum_{i=m_{1}(1)}^{h_{1}(1)} \mu^{I_{1}}(i)\phantom{*}\mbox{and}\phantom{*}U_2:=\sum_{i=m_{1}(1)}^{h_{1}(1)} \Big( \frac {\mu^{I_{1}}(i)}{\mu^{I_{1}}(1)} \Big)^{\frac{1}{2}}.
\end{equation*}
We can write a similar estimate for $\Po\Big[{\hat \xi}_2(t^{b'_1})\in [0,m_{1}(1)]\Big]$ that is
\begin{align}
\label{SG2}
\lefteqn{\Po\Big[{\hat \xi}_2(t^{b'_1})\in [0,m_{1}(1)]\Big]}\phantom{******}\nonumber\\
&\geq \sum_{i=0}^{m_{1}(1)} \mu^{I_{1}}(i) - \exp(-\lambda(I_{1})t^{b'_1}) \sum_{i=0}^{m_{1}(1)} \Big( \frac {\mu^{I_{1}}(i)}{\mu^{I_{1}}(2)} \Big)^{\frac{1}{2}}
\end{align}
and let us write
\begin{equation*}
V_1:= \sum_{i=0}^{m_{1}(1)} \mu^{I_{1}}(i)\phantom{*}\mbox{and}\phantom{*}V_2:=\sum_{i=0}^{m_{1}(1)} \Big( \frac {\mu^{I_{1}}(i)}{\mu^{I_{1}}(2)} \Big)^{\frac{1}{2}}.
\end{equation*}
Combining (\ref{SG1}) and (\ref{SG2}) we obtain that
\begin{equation}
\label{SG3}
\Po[{\hat \tau}\leq t^{b'_1}] \geq U_{1}V_{1}- (U_{1}V_{2}+U_{2}V_{1}) \exp(-\lambda(I_{1})t^{b'_1}).
\end{equation}
\medskip

Now, we treat the term $U_{1}V_{1}$.
Note that if $\omega \in \Gamma_{t}\cap \Lambda_t$ we obtain, using (14) of \cite{CP}, that
\begin{equation}
U_{1}=\sum_{i=0}^{m_{1}(1)} \mu^{I_{1}}(i)\geq \mu^{I_{1}}(m_{1})\geq \frac{K'_{1} e^{-V(m_{1}(1))}}{\sum_{y \in I_{1}}e^{-V(y)}} \geq \frac{K'_{1}}{(h_{1}(1)+1)}\ln^{-2K_{0}}t \nonumber\\ 
\end{equation}
where $K_{0}$ is the positive constant from (\ref{blizko}) and $K'_1$ is a positive constant. Applying the same treatment to $V_{1}$, we obtain that 
\begin{equation}
U_{1}V_{1} \geq \frac{K'_{2}}{(h_{1}(1)+1)^{2}} \ln^{-4K_{0}}t \nonumber\\ 
\end{equation}
where $K'_{2}$ is a positive constant.
On $\Lambda_{t}$, we have that $(h_{1}(1)+1)\leq A\ln_3 t \ln^2 t$, so we obtain 
\begin{equation}
\label{SG4}
U_{1}V_{1} \geq K'_{3}\ln_3^{-2} t \ln^{-4(1+K_{0})}t
\end{equation}
with $K'_{3}$ a positive constant.
\medskip

Concerning the second term of (\ref{SG3}), we will see that if $\omega \in \Gamma_{t} \cap \Lambda_{t}$, it is negligible in comparison with the first one.
We use the upper bound
\begin{align*}
U_{1}V_{2}+U_{2}V_{1} 
&\leq (h_{1}(1)+1) \Big[ \frac{1}{(\mu^{I_{1}}(2))^{\frac{1}{2}}} + \frac{1}{(\mu^{I_{1}}(1))^{\frac{1}{2}}} \Big].
\end{align*}
Using (14) of \cite{CP}, it is then elementary to show that if $\omega \in \Gamma_{t} \cap \Lambda_{t}$,
\begin{equation}
\label{SG5}
U_{1}V_{2}+U_{2}V_{1} \leq K'_{4}t^{\ln_3 t} \ln^{5}t
\end{equation}
with $K'_4$ a positive constant.
\medskip

On the other hand, by Proposition 3.1 of \cite{CP} and the fact that by (\ref{ELEVPROP}), $\e(I_1)<\alpha \ln t$, we obtain 
\begin{equation*}
\exp(-\lambda(I_{1})t^{b'_1})= o(\exp(-\ln^{3/2}t)).
\end{equation*}
This shows together with (\ref{SG5}) that the second term of the right-hand side of (\ref{SG3}) is negligible in comparison with the first one. 
Finally, from (\ref{SG4}) and (\ref{SG5}), we deduce that
\begin{equation*}
\Po[{\hat \tau}\leq t^{b'_1}] \geq \frac {K'_{3}}{2}\ln_3^{-2} t \ln^{-4(1+K_{0})}t
\end{equation*}
for $t$ large enough and $\omega \in \Gamma_{t} \cap \Lambda_{t}$.
\medskip

\noindent
Now, dividing the time interval $[0,t^{\alpha+\eps'}]$ into (roughly speaking) $t^{\frac{\eps'}{2}}$ intervals of size $t^{\alpha+\frac {\eps'}{2}}$ and applying the Markov property we get that
\begin{equation*}
\Po[{\hat \tau}> t^{b_1}] \leq \Big(1- K'_{4} \ln_3^{-2} t \ln^{-4(1+K_{0})}t\Big)^{\lfloor t^{\frac {\eps'}{2}}\rfloor}\leq \exp{(-K'_4\ln^9 t)}
\end{equation*}
for $t$ large enough and $\omega \in \Gamma_{t} \cap \Lambda_{t}$.
As a consequence, the term $\Po[B_{1},\sigma_{1}\geq t^{b_1}]$ is negligible in comparison with $\Po[\sigma_{1}< t^{b_1}]$.
\medskip

\medskip

\medskip

%%%%%%%%%%%%%%%%%
To sum up this subsection, we obtained
\begin{equation}
\label{PB1}
\Po[B_{1}] \leq 2C_4t^{-(r_{1}(1)-\alpha-\eps')}=t^{-(r_{1}(1)+o(1))} 
\end{equation}
as $t\rightarrow \infty$ and $\omega \in \Gamma_{t}\cap \Lambda_{t}$.

%%%%%%%%%%%%%%%%%%%%%%%%%%%%%%%%%%%%%%%%%%%%%%%%%%%%%%%%%%%%%%%%
\subsubsection{Upper bounds for $\Po[B_{i}\mid B_{1} \ldots B_{i-1}]$, where $2\leq i \leq N-1$}
\label{SUBsection}
If $h_1(i)=h_1(i-1)$ then we obviously have  $\Po[B_{i}\mid  B_{1} \ldots B_{i-1}]=1$.
If $h_1(i)>h_1(i-1)$ we do the following.
Let us write $b_i=a_{i-1}+\eps'$.
We will use the following decomposition 
\begin{align*}
\lefteqn{\Po[B_{i}\mid  B_{1} \ldots B_{i-1}]}\phantom{****}\nonumber\\
&\leq \Po[B_{i} \cap \{\sigma_{i} -\sigma_{i-1} \geq t^{b_i}\} \mid B_{1} \ldots B_{i-1}] \nonumber\\
&\phantom{**}+ \Po[\sigma_{i}-\sigma_{i-1} < t^{b_i} \mid B_{1} \ldots B_{i-1}].
\end{align*}
As the event $\{\sigma_{i}-\sigma_{i-1} < t^{b_i}\}$ is independent of $B_{1} \ldots B_{i-1}$ we obtain
\begin{align}
\label{B2|B1}
\Po[B_{i} \mid B_{i-1}] 
&\leq \Po[B_{i} \cap \{\sigma_{i} -\sigma_{i-1} \geq t^{b_i}\} \mid B_{1} \ldots B_{i-1}]\nonumber\\
&\phantom{**}+ \Po[\sigma_{i}-\sigma_{i-1} < t^{b_i}].
\end{align}  
{\bf Upper bounds for $\Po[\sigma_{i}-\sigma_{i-1} < t^{b_i}]$}
\medskip

\noindent
Using the Markov property and (\ref{CER}), we obtain for $\omega \in \Gamma_{t}\cap \Lambda_t$
\begin{align}
\label{2G}
\lefteqn{\Po[\sigma_{i}-\sigma_{i-1} < t^{b_i}]}\phantom{*****}\nonumber\\
&= \Po^{h_{1}(i-1)}[\tau_{h_{1}(i)}(\xi_2)<t^{b_i}]
\nonumber\\
&= \int_{0}^{t^{b_i}} \Po^{m_{2}(i-1)}[\tau_{h_{2}(i-1)}(\xi_2)<t^{b_i}-s] d\Po^{h_{1}(i-1)}[\tau_{m_{2}(i-1)}(\xi_2)\leq s]
\nonumber\\
&\leq  \Po^{m_{2}(i-1)}[\tau_{h_{2}(i-1)}(\xi_2)< t^{b_i}]\nonumber\\
&\leq  C_4t^{\eps'}t^{-(r_{2}(i-1)-a_{i-1})}.
\end{align}
{\bf Upper bounds for $\Po[B_{i} \cap \{\sigma_{i} -\sigma_{i-1} \geq t^{b_i}\} \mid B_{1} \ldots B_{i-1}]$}
\medskip

\noindent
As the event $B_{1} \ldots B_{i-1}$ belongs to ${\cal F}_{\sigma_{i-1}}$ the $\sigma$-field generated by $\xi_1$ and $\xi_2$ up to the stopping time $\sigma_{i-1}$, applying the Markov property, we obtain
\begin{align*}
\lefteqn{\Po[B_{i} \cap \{\sigma_{i} -\sigma_{i-1} \geq t^{b_i}\} \mid B_{1} \ldots B_{i-1}]}\phantom{*****}\nonumber\\
&= \sum_{x=0} ^{h_{1}(i-1)} \Po[B_{i} \cap \{\sigma_{i} -\sigma_{i-1} \geq t^{b_i}\} \mid \xi_1(\sigma_{i-1})=x]\nonumber\\
&\phantom{*********}\times \Po[\xi_1(\sigma_{i-1})=x\mid B_{1} \ldots B_{i-1}] \nonumber\\
&\leq \max_{x< h_1(i-1)} \Po[B_{i} \cap \{\sigma_{i} -\sigma_{i-1} \geq t^{b_i}\} \mid \xi_1(\sigma_{i-1})=x].\nonumber
\end{align*} 
The next step is to find an upper bound for \[\Po[B_{i} \cap \{\sigma_{i} -\sigma_{i-1} \geq t^{b_i}\} \mid \xi_1(\sigma_{i-1})=x]\] uniformly in $x$ for $x\in [0,h_1(i-1))$.\\

As we did in subsection \ref{Riton}, we can couple the two random walks with the random walks restricted to the interval $I_{i}=[0,h_{1}(i)]$. We immediately conclude that
\begin{equation*}
\Po[B_{i} \cap \{\sigma_{i} -\sigma_{i-1} \geq t^{b_i}\}\mid B_{1} \ldots B_{i-1}] \leq \exp{(-K'_4\ln^9 t)}, 
\end{equation*}
which is negligible in comparison with the right-hand side of (\ref{2G}).
\medskip

\medskip

\medskip
Thus, in the case $h_1(i)> h_1(i-1)$, we obtain
\begin{equation}
\label{PBi}
\Po[B_{i}\mid B_{1} \ldots B_{i-1}] \leq 2C_4t^{-(r_2(i-1)-a_{i-1}-\eps')}=t^{-(r_{2}(i-1)-a_{i-1}+o(1))} 
\end{equation}
as $t \rightarrow \infty$ and $\omega \in \Gamma_{t}\cap \Lambda_{t}$.

%%%%%%%%%%%%%%%%%%%%%%%%%%%%%%%%%%%%%
\subsubsection{Upper bound for $\Po[B_{N}\mid B_{1} \ldots B_{N-1}]$}
If $h_1(N)=h_1(N-1)$ then we obviously have  $\Po[B_{N}\mid B_{1} \ldots B_{N-1}]=1$.\\
 If $h_1(N)> h_1(N-1)$ we use the following decomposition.
\begin{align}
\label{Deck}
\lefteqn{\Po[B_{N} \mid B_{1}\ldots B_{N-1}]}\phantom{****} \nonumber\\
&\leq \Po[B_{N} \cap \{\sigma_{N} -\sigma_{N-1} \geq t^{a_{N-1}+\eps'}\} \mid B_{1}\ldots B_{N-1}]\nonumber\\
& \phantom{**}+ \Po[\sigma_{N}-\sigma_{N-1} < t^{a_{N-1}+\eps'}].
\end{align}
The upper estimate of the last term of (\ref{Deck}) is easily computed.
\begin{equation}
\label{TP3}
\Po[\sigma_{N}-\sigma_{N-1} < t^{a_{N-1}+\eps'}] \leq 2C_4t^{-(r_{2}(N-1)-a_{N-1}-\eps')} 
\end{equation}
for $t$ large enough and $\omega \in \Gamma_{t}\cap \Lambda_{t}$.
\medskip

It is not much more difficult to show that the first term of (\ref{Deck}) is negligible in comparison with the last one. To this end, just notice that
\begin{align*}
\lefteqn{\Po[B_{N} \cap \{\sigma_{N} -\sigma_{N-1} \geq t^{a_{N-1}+\eps'}\} \mid B_{1}\ldots B_{N-1}]}\phantom{**}\nonumber\\
&\leq
\Po[T_2 \notin [\sigma_{N-1}, \sigma_N \wedge t^{1+\delta}], \sigma_N \wedge t^{1+\delta}-\sigma_{N-1}\geq t^{a_{N-1}+\eps' }\mid B_{1}\ldots B_{N-1}]
\end{align*}
holds for sufficiently large $t$ since $a_{N-1}+\eps'<1+\delta/2$ for sufficiently large $t$.
This last term can be bounded from above by the same method we used for the term  $\Po[B_{1} \cap \sigma_{1} \geq t^{b_1}]$ of section \ref{Riton}. A similar estimate can be computed which shows that it is negligible in comparison with (\ref{TP3}). 
\medskip

\medskip

Finally in the case $h_1(N)> h_1(N-1)$ we obtained 
\begin{equation}
\label{PBN}
\Po[B_{N}\mid B_{1} \ldots B_{N-1}] \leq 2C_4t^{-(r_2(N-1)-a_{N-1}-\eps')}=t^{-(r_{2}(N-1)-a_{N-1}+o(1))} 
\end{equation}
as $t \rightarrow \infty$ and $\omega \in \Gamma_{t}\cap \Lambda_{t}$.

%%%%%%%%%%%%%%%%%%%%%%%%
\subsubsection{Conclusion}
\label{Concup2}
Using the results of the precedent sections, we now show that 
\begin{equation}
\prod_{i=1}^{N}\Po[B_{i}\mid B_{1} \ldots B_{i-1}] \leq t^{-(r_1(N)+o(1))} 
\end{equation}
as $t \rightarrow \infty$ and $\omega \in \Gamma_{t}\cap \Lambda_{t}$.\\
We will proceed by induction. For every $i\leq 1$, consider the product 
\[\Pi_i=\prod_{j=1}^{i}\Po[B_{j}\mid B_{1} \ldots B_{j-1}].\]
Then we make the following induction hypothesis 
\begin{equation*}
\Pi_i\leq t^{-(r_1(i)+o(1))} 
\end{equation*}
for some $i\geq 1$.
We will show that it implies 
\begin{equation*}
\Pi_{i+1}\leq t^{-(r_1(i+1)+o(1))}. 
\end{equation*}
Indeed we have,
\begin{equation*}
\Pi_{i+1}= \Pi_i \Po[B_{i+1}\mid B_{1} \ldots B_{i}]. 
\end{equation*}
If $h_1(i+1)=h_1(i)$ we have $\Po[B_{i+1}\mid B_{1} \ldots B_{i}]=1$ and $\Pi_{i+1}= \Pi_i$.
But in this case note that by construction we have $r_1(i+1)=r_1(i)$.
So, we obtain
\begin{equation*}
\Pi_{i+1}\leq t^{-(r_1(i+1)+o(1))}.
\end{equation*}
If $h_1(i+1)>h_1(i)$ we have by (\ref{PBi}) and (\ref{PBN}) \[\Po[B_{i+1}\mid B_{1} \ldots B_{i}]\leq t^{-(r_{2}(i)-a_{i}+o(1))}\]
which implies with the induction hypothesis
\begin{equation*}
\Pi_{i+1}\leq t^{-(r_1(i)-a_{i}+r_{2}(i)+o(1))}.
\end{equation*}
But in this case we have by construction either $r_1(i+1)=r_1(i)-a_i+r_2(i)$ or $r_1(i+1)=r_2(i)$ and $r_1(i)=a_i$. This leads to 
\begin{equation*}
\Pi_{i+1}\leq t^{-(r_1(i+1)+o(1))}.
\end{equation*}
As the induction hypothesis is verified for $i=1$ by (\ref{PB1}), we showed that 
\begin{equation}
\label{PBN1}
\prod_{i=1}^{N}\Po[B_{i}\mid B_{1} \ldots B_{i-1}] \leq t^{-(r_1(N)+o(1))} 
\end{equation}
as $t \rightarrow \infty$ and $\omega \in \Gamma_{t}\cap \Lambda_{t}$.

We are now able to deduce the first part of Theorem \ref{theo2}.
Indeed, from Lemma \ref{l_good_omega}, (\ref{DF}) and (\ref{PBN1}),  we have $\IP$-a.s.,
\begin{equation}
\Po[T_2>t^{1+\delta}]\leq t^{-(r_1(N)-o(1))}\nonumber\\
\end{equation}
as $t \rightarrow \infty$.
We obtain $\IP$-a.s.,
\begin{equation}
\label{RE}
\Po[T_2>t]\leq t^{-\Big(\frac{r_1(N)}{1+\delta}+o(1)\Big)}
\end{equation}
as $t \rightarrow \infty$.
We can immediately deduce the first part of Theorem \ref{theo2}. 
Indeed, as $r_1(N)\geq 1$, we obtain $\IP$-a.s.,
\begin{equation}
\Po[T_2>t]\leq t^{-\Big(\frac{1}{1+\delta}+o(1)\Big)}\nonumber\\ 
\end{equation}
as $t \rightarrow \infty$.
To conclude, fix $c<1$, as $\delta$ is arbitrary we can choose $\delta$ sufficiently small such that $\delta<\frac{1-c}{c}$. This shows that $\Eo[T_2^{c}]<\infty$, $\IP$-a.s.
\medskip

\medskip

\medskip

Now, in order to prove Theorem \ref{theo1}, we shall improve (\ref{RE}). That is why we shall consider convergence in $\IP$-probability instead of $\IP$-a.s.\ to obtain
\begin{equation}
\label{RE1}
\ln \Po[T_2>t]\leq -\zeta_2(t)\ln t + o(\ln t)
\end{equation}
as $t \rightarrow \infty$.
At this point, we can explain why we chose to bound $\Po[T_2>t^{1+\delta}]$ instead of $\Po[T_2>t]$. We know by definition of the elevation that $\e[0,h_1(N)]<\ln t$ but when considering convergence $\IP$-a.s., we do not control the elevation and it can be arbitrarily close to $\ln t$. In order to show (\ref{DF}) and to bound the term $\Po[B_N\mid B_1\dots B_{N-1}]$, we had to consider an order of time greater than $t$, that is $t^{1+\delta}$. Considering convergence in $\IP$-probability things become easier since from now on, we can control the elevation of the interval $[0,h_1(N)]$ as follows. 
Fix $t>0$ and $0<\rho<1$. Let $\Xi_t$ be the set of environments $\omega$ such that for every $\omega \in \Xi_t$ we have $\e[0,h_1(N)]<(1-\rho)\ln t$.
By Lemma 4.1 of \cite{CP}, for any $\epsilon>0$, we can choose $\rho>0$ small enough in such a way that $\IP[\Xi_t]> (1-\epsilon)\ln t$. Furhermore, $\IP[\Xi_t]$ does not depend on $t$. Since our goal is now to prove convergence in $\IP$-probability,  we can restrict ourselves to the set $\Xi_t\cap \Lambda_t \cap \Gamma_t$ and bound directly $\Po[T_2>t]$ instead of $\Po[T_2>t^{1+\delta}]$. Finally, to obtain (\ref{RE1}) we can repeat exactly the same computations we used to get (\ref{RE}).

%%%%%%%%%%%%%%%%%%%%%%%%%%%%%%%%%%%%%%%%
\subsection{Lower bound for $\Po[T_2>t]$}
In order to bound  $\Po[T_2>t]$ from below, the main idea  is to find the best strategy for the two random walks not to meet until time $t$. To this end we define the following events.
First, let us write $\tau^2_{i}(\xi_{2})$ for the time of second return to state $i$ for the random walk $\xi_2$. Let
\\
\begin{equation}
D_{0}=\{ \tau_{\{0,2\}}(\xi_1)>\tau_{m'_2(1)}(\xi_{2})\}\cap \{\xi_{2}(u)\leq \xi_{2}(s),0\leq u \leq s \leq \tau_{m'_2(1)}(\xi_{2}) \}. \nonumber\\
\end{equation} 
Then for $1\leq i < N$, 
\begin{align*}
D_{i}=&\{\tau_{h_{1}(i)}(\xi_{1})-\tau_{m'_{2}(i)}(\xi_{2})>t^{a_{i}}\}\cap \{\tau_{m'_2(i+1)}(\xi_{2})<\tau^2_{h_{1}(i)}(\xi_{2})\}\nonumber\\
&\phantom{******} \cap \{\tau_{m'_{2}(i+1)}(\xi_{2})-\tau_{m'_{2}(i)}(\xi_{2})\leq t^{a_{i}}\}, \nonumber
\end{align*}
and finally
\begin{equation*}
D_{N}=\{ \tau_{h_1(N)}(\xi_{1})-\tau_{m'_{2}(N)}(\xi_{2})>t\}\cap \{\tau_{h_1(N)}^{2}(\xi_{2})-\tau_{m'_{2}(N)}(\xi_{2})>t\}.
\end{equation*}
On the event $D_0$, we want the random walk $\xi_1$ to stay at position 1 until $\xi_2$ reaches the point $m'_2(1)$.  On the events $D_i$, for $1\leq i<N$, we want $\xi_1$ to reach for the first time the point $h_1(i)$ in a time greater than $t^{a_i}$ and we want $\xi_2$ to reach the point $m'_2(i+1)$ in a time less than $t^{a_i}$ before it reaches $h_1(i)$ for the second time. Finally, on $D_N$, we want the random walks to stay in their respective $t$-stable wells until time $t$.
\medskip

\noindent
Thus we obtain
\begin{align}
\label{decLow1}
\Po[T_2>t] 
&\geq \Po[D_{0}\cap \ldots \cap D_{N}]
\nonumber\\
&\geq  \Po[D_{0}]\Po[D_{1}\mid D_{0}]\ldots \Po[D_{N}\mid D_{0}D_{1}\ldots D_{N-1}].
\end{align}
We will now find a lower bound for each term of the right-hand side of (\ref{decLow1}).

%%%%%%%%%%%%%%%%%%%%%%%%%%%%%%%%%%% 
\subsubsection{Lower bound for $\Po[D_{0}]$}

First of all, note that by condition B, the random walks remain in mean a time smaller than $\kappa/2$ in every state.
Thus, we can obtain a lower estimate of the event $D_{0}$ writing
\begin{align*}
D_{0}\supset
& \{ \tau_{\{0,2\}}(\xi_1)>\tau_{m'_2(1)}(\xi_{2})\}\cap \{\xi_{2}(u)\leq \xi_{2}(s),0\leq u \leq s \leq \tau_{m'_2(1)}(\xi_{2}) \}\nonumber\\
&\phantom{*} \cap \Big \{ \tau_{m'_2(1)}(\xi_{2})\leq \frac{m'_{2}(1)\kappa}{2} \Big \}
\end{align*}
which implies by independence of the two  random walks
\begin{align*}
\Po[D_{0}]
&\geq \Po\Big[\xi_{2}(u)\leq \xi_{2}(s),0\leq u \leq s\leq \tau_{m'_2(1)}(\xi_{2}),\tau_{m'_2(1)}(\xi_{2})\leq \frac{m'_{2}(1)\kappa}{2}\Big]\nonumber\\
&\phantom{**} \times \Po\Big[\tau_{\{0,2\}}(\xi_1)>\frac{m'_{2}(1)\kappa}{2}\Big].
\end{align*}
Then, using condition B and the fact that jump Markov processes have independent and exponentially distributed jump time increments, we can write
\begin{align*}
\Po[D_{0}]
&\geq  \Po\Big[\tau_{\{0,2\}}(\xi_1)>\frac{m'_{2}(1)\kappa}{2}\Big]\Po[\xi_{2}(s)\leq \xi_{2}(t),0\leq s \leq t \leq \tau_{m'_2(1)}(\xi_{2})] \nonumber\\
&\phantom{**}\times \Big(P_{X}\Big[X\leq \frac{\kappa}{2}\Big]\Big)^{m'_2(1)}
\end{align*}
where $X$ is a random variable with exponential law of parameter $2/\kappa$ and $P_{X}$ is the law of $X$.
From the last expression we obtain 
\begin{equation}
\label{EV0}
\Po[D_{0}] \geq \Big(\frac{e^{-\kappa^2}(e-1)}{2e\kappa^2}\Big)^{\ln^{1/2}t}
\end{equation}
for $\omega \in \Gamma_t\cap \Lambda_t$.
Note that, as $\kappa>1$, we have that $\Big(\frac{e^{-\kappa^2}(e-1)}{2e\kappa^2}\Big)<1$.

%%%%%%%%%%%%%%%%%%%%%%%%
\subsubsection{Lower bound for $\Po[D_{1}\mid D_{0}]$}
\label{sub0low}

First, for the sake of brevity let us write $\beta_i= \tau_{m'_2(i)}(\xi_2)$ for $i\geq 1$. Using again the Markov property and the independence of the two random walks, we obtain
\begin{align}
\label{Low1}
\Po[D_{1}\mid D_{0}]
&= \Po^{m'_2(1)}[\tau_{m'_2(2)}(\xi_2)<\tau_{h_1(1)}(\xi_2),\tau_{m'_2(2)}(\xi_2)\leq t^{a_1}]\nonumber\\
& \phantom{**}\times  \Po^1[\tau_{h_1(1)}(\xi_1)>t^{a_1}].
\end{align}
We first show that the second term of the right-hand side of (\ref{Low1}) is greater than $1/2$ for $t$ sufficiently large and $\omega \in \Gamma_t\cap \Lambda_t$. Observe that by the Markov property we have
\begin{equation}
\label{Fluck}
\Po^1[\tau_{h_1(1)}(\xi_1)>t^{a_1}] \geq \Po^{m_1(1)}[\tau_{h_1(1)}(\xi_1)>t^{a_1}].
\end{equation}
We can apply (\ref{CER}) to the right-hand side of (\ref{Fluck}) to obtain that 
\begin{equation*}
\Po^{1}[\tau_{h_1(1)}(\xi_1)>t^{a_1}]\geq (1-C_4t^{a_1}e^{-V(h_1(1))+V(m_1(1))})^{+}
\end{equation*}
where $(\cdot)^+$ is the positive part.
One can see that for $\omega \in \Gamma_t\cap \Lambda_t$ and $t$ sufficiently large the right-hand side of the last inequality will be larger than $1/2$. By construction, we have $W(h_1(1))-W(m_1(1))\geq a_1\ln t$ so we obtain
\begin{equation}
\label{TY1}
C_4t^{a_1}e^{-V(h_1(1))+V(m_1(1))}\leq C_4\ln^{-2K_0}t\leq \frac{1}{2}
\end{equation}\\
for $t$ sufficiently large.

Now we will treat the first term of the right-hand side of (\ref{Low1}).
We have to consider two cases:
\begin{itemize}
\item $m'_2(1)$ is $t^{a_1}$-stable;
\item $m'_2(1)$ is not $t^{a_1}$-stable.
\end{itemize}
\medskip

\noindent
%First case: $m'_2(1)$ is $t^{a_1}$-stable.
{\textbf{Case 1:} $m'_2(1)$ is $t^{a_1}$-stable}\\
First using the Markov property let us write
\begin{align}
\label{EQp}
\lefteqn{
\Po^{m'_2(1)}[\tau_{m'_2(2)}(\xi_2)<\tau_{h_1(1)}(\xi_2),\tau_{m'_2(2)}(\xi_2)\leq t^{a_1}]}\phantom{*******}\nonumber\\
&\geq
\Po^{m'_2(1)}\Big[\tau_{h_2(1)}(\xi_2)<\tau_{h_1(1)}(\xi_2),\tau_{h_2(1)}(\xi_2)\leq \frac{t^{a_1}}{2}\Big]\nonumber\\
 &\phantom{**}\times
\Po^{h_2(1)}\Big[\tau_{m'_2(2)}(\xi_2)<\tau_{m'_2(1)}(\xi_2),\tau_{m'_2(2)}(\xi_2)\leq \frac{t^{a_1}}{2}\Big]\nonumber\\
&:=  R_1\times U_1.
\end{align}
Let us bound from below the first term of the right-hand side of (\ref{EQp}). 
\medskip

\begin{figure}[!htb]
\begin{center}
\includegraphics[scale= 0.7]{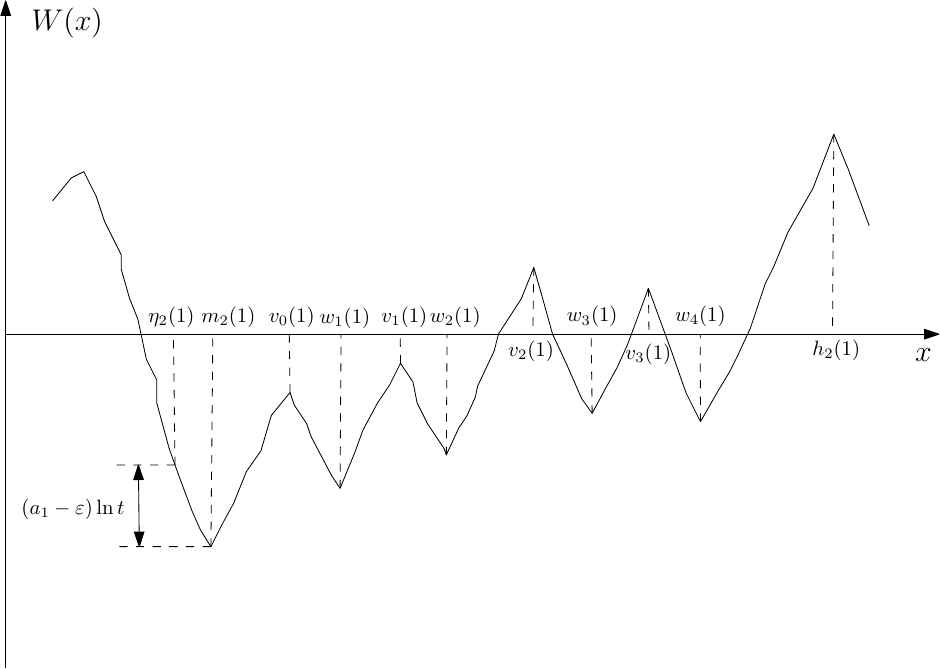}
\caption{On the definition of $\eta_2(1)$, $v_i(1)$ and $w_i(1)$. Case $J^1_2(1)=4$.}
\label{fig5}
\end{center}
\end{figure}

\noindent
% Lower bound for $T_1$
{\bf Lower bound for $R_1$}\\
In this case $m'_2(1)\equiv m_2(1)$. Fix $t>e$ and let $\eps=\eps(t)$ from Lemma \ref{l_good_omega}. Then define
 \[\eta_{2}(1)=\max \{x<m_{2}(1); W(x)-W(m_{2}(1))=(a_{1}-\eps)\ln t \}. \]
Futhermore, suppose that we have $J^1_2(1)=J^1_2(1)(\omega,t)$ $t^{a_1-\eps}$-stable wells in the interval $[m_2(1),h_2(1)]$. If $J^1_2(1)\geq 1$, then we define for $1 \leq i \leq J^1_2(1)$ the point $w_i(1)$ as the minimum of the $i$th $t^{a_1-\eps}$-stable well and  \[v_i(1)=\argmax_{w_i(1)<x<w_{i+1}(1)}W(x).\] We also define $v_0(1):= \arg \max_{m_2(1)<x<w_{1}(1)}W(x)$. See Figure \ref{fig5}.

\medskip

\noindent
If $J^1_2(1)=0$, we can directly find a lower bound for $R_1$.
We write
\begin{equation*}
R_1\geq \Po^{m_2(1)}[\tau_{h_2(1)}(\xi_2)<\tau_{\eta_2(1)}(\xi_2)]-\Po^{m_2(1)}\Big[\tau_{\{\eta_2(1),h_2(1)\}}(\xi_2)> \frac{t^{a_1}}{2}\Big].
\end{equation*}
Then using (\ref{exit}) for the first term and (\ref{LOWB}) for the second term (noting that $H([\eta_2(1),h_2(1)])\leq (a_1-\eps)\ln t+2K_0\ln_2t$), we obtain for $\omega \in \Gamma_t\cap \Lambda_t$ 
\begin{align}
\label{TY2}
R_1&\geq \frac{e^{-(W(h_2(1))-W(\eta_2(1)))}\ln^{-2K_0}t}{h_2(1)-\eta_2(1)}-\exp \Big( -\frac{t^{\eps}}{C_1\ln^{2K_0}t(h_2(1)-\eta_2(1))^4}\Big)\nonumber\\
&= \frac{e^{-(W(h_2(1))-W(\eta_2(1)))}\ln^{-2K_0}t}{h_2(1)-\eta_2(1)}-o(\exp(-t^{\frac{\eps}{2}}))
\end{align}
as $t \rightarrow \infty$.
If $J^1_2(1)\geq 1$ we have a bit more work to do.
By the Markov property we obtain 
\begin{align}
\label{Prod2}
R_1
&\geq \prod_{i=0}^{J^1_2(1)}\Po^{w_i(1)}\Big[\tau_{w_{i+1}(1)}(\xi_2)<\tau_{v_{i-1}(1)}(\xi_2),
 \tau_{w_{i+1}(1)}(\xi_2)\leq \frac{t^{a_1}}{2(J^1_2(1)+1)}\Big]
\end{align}
with the conventions $v_{-1}(1)=\eta_2(1)$, $w_0(1)=m_2(1)$ and $w_{J^1_2(1)+1}(1)=h_2(1)$.\\
Let us first bound from below the first term of the product. The other terms of the product will be treated in a similar way.
First note that by the Markov property we have
\begin{align}
\label{EQR}
\lefteqn{
\Po^{m_2(1)}\Big[\tau_{w_{1}(1)}(\xi_2)<\tau_{\eta_{2}(1)}(\xi_2),\tau_{w_{1}(1)}(\xi_2)\leq \frac{t^{a_1}}{2(J^1_2(1)+1)}\Big]}\phantom{****}\nonumber\\
&\geq
\Po^{m_2(1)}\Big[\tau_{v_{0}(1)}(\xi_2)<\tau_{\eta_{2}(1)}(\xi_2),\tau_{v_{0}(1)}(\xi_2)\leq \frac{t^{a_1}}{4(J^1_2(1)+1)}\Big]\nonumber\\
 &\phantom{**}\times
\Po^{v_0(1)}\Big[\tau_{w_{1}(1)}(\xi_2)<\tau_{m_{2}(1)}(\xi_2),\tau_{w_{1}(1)}(\xi_2)\leq \frac{t^{a_1}}{4(J^1_2(1)+1)}\Big]\nonumber\\
&:= R_2\times R_3.
\end{align}
We have
\begin{equation*}
R_2
 \geq \Po^{m_2(1)}[\tau_{v_{0}(1)}(\xi_2)<\tau_{\eta_{2}(1)}(\xi_2)]- \Po\Big[\tau_{\{v_{0}(1),\eta_2(1)\}}(\xi_2)> \frac{t^{a_1}}{4(J^1_2(1)+1)}\Big]
\end{equation*}
Then, using (\ref{exit}) for the first term and  (\ref{LOWB}) for the second term, we obtain for $\omega \in \Gamma_t$ 
\begin{align}
\label{TY3}
R_2
&\geq \frac{e^{-(W(v_0(1))-W(\eta_2(1)))}\ln^{-2K_0}t}{v_0(1)-\eta_2(1)}-\exp \Big( -\frac{t^{\eps}}{ K'_1J^1_2(1)\ln^{2K_0}t(v_0(1)-\eta_2(1))^4}\Big) \nonumber\\
&=\frac{e^{-(W(v_0(1))-W(\eta_2(1)))}\ln^{-2K_0}t}{v_0(1)-\eta_2(1)}-o(\exp(-t^{\frac{\eps}{2}}))
\end{align}
as $t \rightarrow \infty$.
\medskip

\noindent
For $R_3$ we write
\begin{align}
R_3
&\geq \Po^{v_0(1)}[\tau_{w_{1}(1)}(\xi_2)<\tau_{m_{2}(1)}(\xi_2)]\nonumber\\
&\phantom{************}- \Po^{v_0(1)}\Big[\tau_{\{w_{1}(1),m_2(1)\}}(\xi_2)> \frac{t^{a_1}}{4(J^1_2(1)+1)}\Big].
\end{align}
Analogously, using (\ref{exit}) for the first term and (\ref{LOWB2}) for the second term, we obtain for $\omega \in \Gamma_t\cap \Lambda_t$ 
\begin{equation}
\label{TY4}
R_3
\geq \frac{\ln^{-2K_0}t}{w_1(1)-m_2(1)}-o(\exp{(-t^{\frac{\eps}{3}})})
\end{equation}
as $t \rightarrow \infty$. To obtain $o(\exp{(-t^{\frac{\eps}{3}})})$ we used Proposition 3.1 of~\cite{CP} and the facts that on $\Gamma_t\cap \Lambda_t$, $\Delta_1\leq B\ln_3t\ln^2t$, $\gamma_1\leq 2\ln_3t\ln t$, $J^1_2(1)\leq \ln_2 t$ and $a_1-\frac{\e(I^+)\vee \e(I^-)}{\ln t}\geq \eps$.
\\

We will now bound from below the terms of the right-hand side product of (\ref{Prod2}) for $1\leq i \leq J^1_2(1)$.
First note that by the Markov property we have
\begin{align*}
\lefteqn{
\Po^{w_i(1)}\Big[\tau_{w_{i+1}(1)}(\xi_2)<\tau_{v_{i-1}(1)}(\xi_2),\tau_{w_{i+1}(1)}(\xi_2)\leq \frac{t^{a_1}}{2(J^1_2(1)+1)}\Big]}\phantom{***}\nonumber\\
&\geq
\Po^{w_i(1)}\Big[\tau_{v_{i}(1)}(\xi_2)<\tau_{v_{i-1}(1)}(\xi_2),\tau_{v_{i}(1)}(\xi_2)\leq \frac{t^{a_1}}{4(J^1_2(1)+1)}\Big]\nonumber\\
 &\phantom{**}\times
\Po^{v_i(1)}\Big[\tau_{w_{i+1}(1)}(\xi_2)<\tau_{w_{i}(1)}(\xi_2),\tau_{w_{i+1}(1)}(\xi_2)\leq \frac{t^{a_1}}{4(J^1_2(1)+1)}\Big].
\end{align*}
\\
Let us define 
\begin{equation*}
R^i_4:=\Po^{w_i(1)}\Big[\tau_{v_{i}(1)}(\xi_2)<\tau_{v_{i-1}(1)}(\xi_2),\tau_{v_{i}(1)}(\xi_2)\leq \frac{t^{a_1}}{4(J^1_2(1)+1)}\Big]
\end{equation*}
and
\begin{equation*}
R^i_5:=\Po^{v_i(1)}\Big[\tau_{w_{i+1}(1)}(\xi_2)<\tau_{w_{i}(1)}(\xi_2),\tau_{w_{i+1}(1)}(\xi_2)\leq \frac{t^{a_1}}{4(J^1_2(1)+1)}\Big]
\end{equation*}
for $1\leq i \leq J^1_2(1)$.
\\

Now, to bound $R^i_4$, we are going to distinguish two cases: $W(v_{i-1}(1))<W(v_i(1))$ and $W(v_{i-1}(1))>W(v_i(1))$.
In the first case, we can directly use a similar computation as that we used for $R_2$. In the second case note that $W(v_{i-1}(1))-W(v_i(1))<\eps \ln t$ since $(a_1-\eps)\ln t \leq W(v_i(1))-W(w_i(1))< a_1\ln t$.
With this observation, we can use a similar computation as that we used for $R_2$. So, suming up, we obtain in both cases
\begin{align}
\label{TY6}
R^i_4 &\geq \frac{e^{-(W(v_i(1))-W(v_{i-1}(1)))}t^{-\eps}\ln^{-2K_0}t}{v_i(1)-v_{i-1}(1)}\nonumber\\
&\phantom{************}-\exp \Big( -\frac{t^{\eps}}{K'_1\ln^{-2K_0}tJ^1_2(1)(v_i(1)-v_{i-1}(1))^4}\Big)\nonumber\\
&\geq \frac{e^{-(W(v_i(1))-W(v_{i-1}(1)))}t^{-\eps}\ln^{-2K_0}t}{v_i(1)-v_{i-1}(1)}-o(\exp(-t^{\frac{\eps}{2}})
\end{align}
as $t\rightarrow \infty$ and $\omega \in \Gamma_t\cap \Lambda_t$.\\
For $R^i_5$ we use the decompositon
\begin{align*}
R^i_5
&\geq \Po^{v_i(1)}[\tau_{w_{i+1}(1)}(\xi_2)<\tau_{w_{i}(1)}(\xi_2)]\nonumber\\
&\phantom{************}- \Po^{v_i(1)}\Big[\tau_{\{w_{i}(1),w_{i+1}(1)\}}(\xi_2)> \frac{t^{a_1}}{4(J^1_2(1)+1)}\Big].
\end{align*}
Then using (\ref{exit}) for the first term and  (\ref{LOWB2}) for the second term, we obtain for $\omega \in \Gamma_t\cap \Lambda_t$ 
\begin{equation}
\label{TY7}
R^i_5\geq \frac{\ln^{-2K_0}t}{w_{i+1}(1)-w_i(1)}-o(\exp{(-t^{\frac{\eps}{3}})}).
\end{equation}
as $t\rightarrow \infty$.
\medskip

Finally from  (\ref{TY2}), (\ref{TY3}), (\ref{TY4}), (\ref{TY6}) and (\ref{TY7})  we obtain for $\omega \in \Gamma_t\cap \Lambda_t$ and sufficiently large $t$
\begin{equation}
\label{Gro1}
R_1\geq \Big[\frac{t^{-\eps}\ln^{-2K_0}t}{(2h_2(1))^2}\Big]^{J^1_2(1)+1} t^{-(r_2(1)-a_1)}.
\end{equation}
\medskip

\noindent
Now let us go back to equation (\ref{EQp}).
\medskip

\begin{figure}[!htb]
\begin{center}
\includegraphics[scale= 0.7]{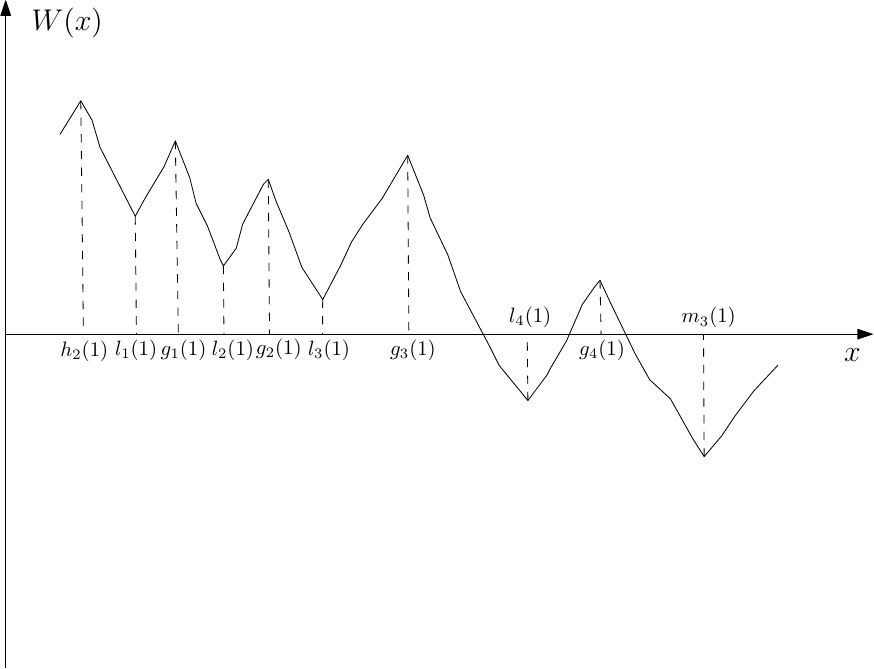}
\caption{On the definition of $l_i(1)$ and $g_i(1)$. Case $J^2_2(1)=4$.}
\label{fig6}
\end{center}
\end{figure}

\noindent
%Lower bound for $U_1$
{\bf Lower bound for $U_1$}\\
First observe that in this case, $m'_2(2)=m_3(1)$. Now, we divide the interval $[h_2(1),m_3(1)]$ into $t^{a_1-\eps}$-stable wells. We denote by $J^2_2(1)=J^2_2(1)(\omega,t)$ the number of $t^{a_1-\eps}$-stable wells in the interval $[h_2(1),m_3(1)]$. If $J^2_2(1)\geq 1$, then we define for $1 \leq i \leq J^2_2(1)$ the point $l_i(1)$ as the minimum of the $i$th $t^{a_1-\eps}$-stable well and  $g_i(1):=\arg \max_{l_i(1)<x<l_{i+1}(1)}W(x)$. We also define $g_{J^2_2(1)}^1=\arg \max_{l_{J^2_2(1)}(1)<x<m_3(1)}W(x)$. See Figure \ref{fig6}.

If $J^2_2(1)=0$ we can directly compute a lower bound for $U_1$.
\begin{align*}
U_1
&\geq
\Po^{h_2(1)}\Big[\tau_{m_3(1)}(\xi_2)<\tau_{w_{J^1_2(1)}(1)}(\xi_2),\tau_{\{w_{J^1_2(1)}(1),m_3(1)\}}(\xi_2)\leq \frac{t^{a_1}}{2}\Big]\nonumber\\
&\geq \Po^{h_2(1)}\Big[\tau_{m_3(1)}(\xi_2)<\tau_{w_{J^1_2(1)}(1)}(\xi_2)\Big]- \Po^{h_2(1)}\Big[\tau_{\{w_{J^1_2(1)}(1),m_3(1)\}}(\xi_2)> \frac{t^{a_1}}{2}\Big].
\end{align*}
Then using (\ref{exit}) for the first term and (\ref{LOWB2}) for the second term, we obtain for $\omega \in \Gamma_t\cap \Lambda_t$ 
\begin{equation}
\label{TY8}
U_1\geq \frac{\ln^{-2K_0}t}{m_3(1)-w_{J_2^1}(1)}-o(\exp(-t^{\frac{\eps}{3}})
\end{equation}
as $t \rightarrow \infty$.

If  $J^2_2(1)\geq 1$ we have a bit more work to do.
We will use the following decomposition
\begin{align}
\label{Low2}
U_1 &\geq \Po^{h_2(1)}\Big[\tau_{l_{1}(1)}(\xi_2)<\tau_{w_{J^1_2(1)}(1)}(\xi_2),\tau_{l_{1}(1)}(\xi_2)\leq \frac{t^{a_1}}{2(J^2_2(1)+1)}\Big]\nonumber\\
&\phantom{**}\times \prod_{i=1}^{J^2_2(1)}\Po^{l_i(1)}\Big[\tau_{l_{i+1}(1)}(\xi_2)<\tau_{g_{i-1}(1)}(\xi_2),\tau_{l_{i+1}(1)}(\xi_2)\leq \frac{t^{a_1}}{2(J^2_2(1)+1)}\Big].
\end{align}\\
Let us define \[U^i_2:=\Po^{l_i(1)}\Big[\tau_{l_{i+1}(1)}(\xi_2)<\tau_{g_{i-1}(1)}(\xi_2),\tau_{l_{i+1}(1)}(\xi_2)\leq \frac{t^{a_1}}{2(J^2_2(1)+1)}\Big]\] for $0\leq i \leq J^2_2(1)$.
Then, let us note that by the Markov property,
\begin{align}
U_2^i
&\geq \Po^{l_i(1)}\Big[\tau_{g_i(1)}(\xi_2)<\tau_{g_{i-1}(1)}(\xi_2),\tau_{g_i(1)}(\xi_2)\leq \frac{t^{a_1}}{4(J^2_2(1)+1)}\Big]\nonumber\\
 &\phantom{**} \times \Po^{g_i(1)}\Big[\tau_{l_{i+1}(1)}(\xi_2)<\tau_{l_i(1)}(\xi_2),\tau_{l_{i+1}(1)}(\xi_2)\leq \frac{t^{a_1}}{4(J^2_2(1)+1)}\Big]\nonumber\\
&:=U_3^i\times U_4^i.
\end{align}

Now, to bound the term $U_3^i$ we have to distinguish two cases.
\medskip

\noindent
{\sl Case 1:} $W(g_i(1))<W(g_{i-1}(1))$\\
In this case, by a similar computation as that we used for $R^i_4$ in the case $W(v_{i-1}(1)<W(v_i(1)$, we obtain for $\omega \in \Gamma_t\cap \Lambda_t$ 
\begin{equation}
\label{TY9}
U_3^i \geq \frac{\ln^{-2K_0}t}{g_i(1)-g_{i-1}(1)}-o(\exp(-t^{\frac{\eps}{2}})
\end{equation}
as $t \rightarrow \infty$.
\medskip

\noindent
{\sl Case 2:} $W(g_i(1))>W(g_{i-1}(1))$\\
By a similar computation as that we used for $R^i_4$ in the case $W(v_{i-1}(1))>W(v_i(1))$, we obtain for $\omega \in \Gamma_t\cap \Lambda_t$,
\begin{equation}
\label{TY10}
U_3^i \geq \frac{t^{-\eps}\ln^{-2K_0}t}{g_i(1)-g_{i-1}(1)}-o(\exp(-t^{\frac{\eps}{2}})
\end{equation}
$t \rightarrow \infty$.

For $U_4^i$ we use the decompositon
\begin{equation*}
U_4^i
\geq \Po^{g_i(1)}[\tau_{l_{i+1}(1)}(\xi_2)<\tau_{l_{i}(1)}(\xi_2)]- \Po^{g_i(1)}[\tau_{\{l_{i}(1),l_{i+1}(1)\}}(\xi_2)> \frac{t^{a_1}}{4(J'''_1+1)}].
\end{equation*}
Then using (\ref{exit}) for the first term and (\ref{LOWB2}) for the second term, we obtain for $\omega \in \Gamma_t\cap \Lambda_t$ 
\begin{equation}
\label{TY11}
U_4^i \geq \frac{\ln^{-2K_0}t}{l_{i+1}(1)-l_i(1)}-o(\exp(-t^{\frac{\eps}{3}})
\end{equation}
as $t \rightarrow \infty$.
Finally as the first term  of the right-hand side of (\ref{Low2}) is completely similar to the term $U_3^i$, from  (\ref{TY8}), (\ref{TY9}), (\ref{TY10}) and (\ref{TY11}) we obtain for $\omega \in \Gamma_t\cap \Lambda_t$ and sufficiently large~$t$
\begin{equation}
\label{Gro4}
U_1\geq \Big[\frac{t^{-\eps}\ln^{-2K_0}t}{(2m_3(1))^2}\Big]^{J^2_2(1)+1}.
\end{equation}
\medskip

\noindent
%Second case: $m'_2(1)$ is not $t^{a_1}$-stable.
{\textbf{Case 2:} $m'_2(1)$ is not $t^{a_1}$-stable}\\
In this case we have to distinguish two subcases.
\begin{itemize}
\item $a_1=l_1(1)$

To find a lower bound for $\Po^{m'_2(1)}[\tau_{m'_2(2)}(\xi_2)<\tau_{h_1(1)}(\xi_2),\tau_{m'_2(2)}(\xi_2)\leq t^{a_1}]$ in this case, we can use a similar computation as for the term $U_1$ above. We obtain that for $\omega \in \Gamma_t\cap \Lambda_t$ and sufficiently large $t$,
\begin{equation}
\label{Gro2}
\Po^{m'_2(1)}[\tau_{m'_2(2)}(\xi_2)<\tau_{h_1(1)}(\xi_2),\tau_{m'_2(2)}(\xi_2)\leq t^{a_1}]\geq
\Big[\frac{t^{-\eps}\ln^{-2K_0}t}{(2m_2(1))^2}\Big]^{J^2_1(1)}.
\end{equation}

\item $a_1= r_1(1)$

In this case, we use the decomposition (\ref{EQp}). The unique difference is that, using the Markov property we decompose the term $R_1$ as follows
\begin{align}
\label{EQt}
\lefteqn{
\Po^{m'_2(1)}\Big[\tau_{h_2(1)}(\xi_2)<\tau_{h_1(1)}(\xi_2),\tau_{h_2(1)}(\xi_2)\leq \frac{t^{a_1}}{2}\Big]}\phantom{******}\nonumber\\
&\geq
\Po^{m'_2(1)}\Big[\tau_{m_2(1)}(\xi_2)<\tau_{h_1(1)}(\xi_2),\tau_{m_2(1)}(\xi_2)\leq \frac{t^{a_1}}{4}\Big]\nonumber\\
&\phantom{**}\times
\Po^{m_2(1)}\Big[\tau_{h_2(1)}(\xi_2)<\tau_{h_1(1)}(\xi_2),\tau_{h_2(1)}(\xi_2)\leq \frac{t^{a_1}}{4}\Big].
\end{align}
To find a lower bound for the first term of the right-hand side of (\ref{EQt}) we can use a similar computation as for the term $U_1$. For the second term of the right-hand side of (\ref{EQt}) we can use a similar computation as for $R_1$. 
We obtain for $\omega \in \Gamma_t\cap \Lambda_t$ and sufficiently large $t$
\begin{align}
\label{Gro3}
\lefteqn{\Po^{m'_2(1)}[\tau_{m'_2(2)}(\xi_2)<\tau_{h_1(1)}(\xi_2),\tau_{m'_2(2)}(\xi_2)\leq t^{a_1}]}\nonumber\\
& \phantom{************}\geq \Big[\frac{t^{-\eps}\ln^{-4K_0}t}{(2m_3(1))^2}\Big]^{J^1_2(1)+J^2_1(1)+J^2_2(1)+2} t^{-(r_2(1)-a_1)}.
\end{align}
\end{itemize}

Finally, using (\ref{TY1}), (\ref{Gro1}), (\ref{Gro4}), (\ref{Gro2}) and (\ref{Gro3}) we obtain 
\begin{align}
\label{EV1}
 \Po[D_{1}\mid D_{0}]
&\geq \frac{1}{2}\Big[\frac{t^{-\eps}\ln^{-4K_0}t}{(2m'_2(2))^2}\Big]^{J^1_2(1)+J^2_1(1)+J^2_2(1)+2} t^{-(r_2(1)-a_1)^+}\nonumber\\
&= t^{-((r_2(1)-a_1)^++o(1))}
\end{align}
as $t\rightarrow \infty$ and $\omega \in \Gamma_t\cap \Lambda_t$.

%Lower bounds for $\Po[D_{i}\mid D_{0}\ldots D_{i-1}], 2\leq i < N$
\subsubsection{Lower bounds for $\Po[D_{i}\mid D_{0}\ldots D_{i-1}]$, where $2\leq i < N$}
\label{Hu}

Using the Markov property, we first obtain
\begin{align*}
\label{Lowref}
\lefteqn{\Po[D_{i}\mid D_{0}\ldots D_{i-1}]}\phantom{***}\nonumber\\
=&
\sum_{x=0}^{h_1(i-1)}\Po[D_{i}\mid \xi_{1}(\beta_i)=x,\xi_{2}(\beta_i)=m'_2(i)]\Po[\xi_{1}(\beta_i)=x\mid D_{0}\ldots D_{i-1}].
\end{align*}
Let us bound  $\Po[D_{i}\mid\xi_{1}(\beta_i)=x,\xi_{2}(\beta_i)=m'_2(i)]$  uniformly in $x$ for  $x<h_1(i-1)$. We have by the Markov property and the independence of the two random walks 
\begin{eqnarray}
\label{Low30}
\lefteqn{\Po[D_{i}\mid  \xi_{1}(\beta_i)=x,\xi_{2}(\beta_i)=m'_2(i)]}\nonumber\\
&=&
\Po^{x}[\tau_{h_1(i)}(\xi_1)>t^{a_i}]\Po^{m'_2(i)}[\tau_{m'_2(i+1)}(\xi_2)<\tau_{h_1(i)}(\xi_2),\tau_{m'_2(i+1)}(\xi_2)\leq t^{a_i}]\nonumber\\
&:=& V^i_1\times V^i_2.
\end{eqnarray}
Let us first treat the first term of the right-hand side of (\ref{Low30}).
%%%%%%%%%
\medskip

\noindent
{\bf Lower bounds for $V^i_1$}\\
If $h_1(i)>h_1(i-1)$, then by the Markov property and as $x<h_1(i-1)$ we obtain,
\begin{equation*}
V^i_1 \geq \Po^{m_1(i)}[\tau_{h_1(i)}(\xi_1)>t^{a_i}]
\end{equation*}
As we did in (\ref{Fluck}) we show that for $\omega \in \Gamma_t\cap \Lambda_t$ and $t$ large enough we have 
\begin{equation*}
\Po^{m_1(i)}[\tau_{h_1(i)}(\xi_1)>t^{a_i}]\geq \frac{1}{2}. 
\end{equation*}

If $h_1(i)=h_1(i-1)$, we use the following decomposition
\begin{align}
\label{blu}
V^i_1
&\geq \Po^{x}[\tau_{m_1(i)}(\xi_1)<\tau_{h_1(i)}(\xi_1),\tau_{m_1(i)}(\xi_1)\leq t^{a_i}]\nonumber\\
& \phantom{**}\times \Po^{m_1(i)}[\tau_{h_1(i)}(\xi_1)>t^{a_i}].
\end{align}
Again we have 
\begin{equation}
\label{chuchu}
\Po^{m_1(i)}[\tau_{h_1(i)}(\xi_1)>t^{a_i}]\geq \frac{1}{2}. 
\end{equation}
For the first term of the right-hand side of (\ref{blu}), we can perform exactly the same computation as in subsection \ref{sub0low} for the term $U_1$,  to find that, being $J^1_1(i)$ the number of $t^{a_i-\eps}$-stable wells in the interval $[m_1(i),h_1(i)]$,
\begin{equation}
\label{PU2}
\Po^{x}[\tau_{m_1(i)}(\xi_1)<\tau_{h_1(i)}(\xi_1),\tau_{m_1(i)}(\xi_1)\leq t^{a_i}]
\geq \Big[\frac{t^{-\eps}\ln^{-2K_0}t}{(2h_1(i))^2}\Big]^{J^1_1(i)+1}
\end{equation}
for $\omega \in \Gamma_t \cap \Lambda_t$ and sufficiently large $t$.
\medskip

\noindent
{\bf Lower bounds for $V^i_2$}\\
This term is completely similar to the term \[\Po^{m'_2(1)}[\tau_{m'_2(2)}(\xi_2)<\tau_{h_1(1)}(\xi_2),\tau_{m'_2(2)}(\xi_2)\leq t^{a_1}]\]
 of section \ref{sub0low}.  We follow step by step the method we applied in that section to obtain that
 
 \begin{equation}
\label{PU3}
V^i_2 \geq \Big[\frac{t^{-\eps}\ln^{-4K_0}t}{(2m'_2(i+1))^2}\Big]^{J^1_2(i)+J^2_1(i)+J^2_2(i)+2} t^{-(r_2(i)-a_i)^+}.
 \end{equation}

%%%%%%%%%%%%
\medskip

Finally using (\ref{chuchu}),  (\ref{PU2}) and (\ref{PU3}) we obtain 
\begin{align}
\label{EVi}
 \Po[D_{i}\mid D_{0}\dots D_{i-1}] 
&\geq \frac{1}{2}\Big[\frac{t^{-\eps}\ln^{-4K_0}t}{(2m'_2(i+1))^2}\Big]^{J^1_1(i)+J^1_2(i)+J^2_1(i)+J^2_2(i)+3} t^{-(r_2(i)-a_i)^+}\nonumber\\
&=t^{-((r_2(i)-a_i)^++o(1))}
\end{align}
as $t\rightarrow \infty$ and $\omega \in \Gamma_t\cap \Lambda_t$.

%Lower bounds for $\Po[D_{i}\mid D_{0}\ldots D_{i-1}], i=N,N+1$
\subsubsection{Lower bound for $\Po[D_{N}\mid D_{0}\ldots D_{N-1}]$}

As we have already done in subsection \ref{Hu} , we have by the Markov property 
\begin{align*}
\label{Lowref3}
\lefteqn{\Po[D_{N}\mid D_{1}\ldots D_{N-1}]}\phantom{*****}\nonumber\\
&=
\sum_{x=0}^{h_1(N-1)}\Po[D_{N}\mid \xi_1(\beta_{N})=x,\xi_2(\beta_{N})=m'_2(N)]\nonumber\\
& \phantom{******}\times \Po[\xi_1(\beta_{N})=x,\xi_2(\beta_{N})=m'_2(N)\mid D_{0}\ldots D_{N-1}].
\end{align*}
Then, we have by the Markov property and the independence of the two random walks
\begin{align*}
\lefteqn{\Po[D_{N}\mid \xi_1(\beta_{N})=x,\xi_2(\beta_{N})=m'_2(N)]}\phantom{*********}\nonumber\\
&=
\Po^{x}[\tau_{h_1(N)}(\xi_1)>t]\Po^{m'_2(N)}[\tau_{h_1(N)}(\xi_2)>t]
\end{align*}
with $x<h_1(N-1)$.
As we did in subsection \ref{Hu} for $V^i_1$, one can see that
\begin{equation*}
\label{PF3}
\Po^x[\tau_{h_1(N)}(\xi_1)>t]\geq \frac{1}{2}\Big[\frac{t^{-\eps}\ln^{-2K_0}t}{(2h_1(N))^2}\Big]^{J^1_{1}(N)}
\end{equation*}
for $t$ large enough and $\omega \in \Gamma_{t}\cap \Lambda_t$.
Using (\ref{CER}), we show that for $t$ large enough and $\omega \in \Gamma_{t}\cap \Lambda_t$,
\begin{equation*}
\Po^{m'_2(N)}[\tau_{h_1(N)}(\xi_2)>t]\geq \frac{1}{2}.
\end{equation*}
Therefore, we obtain
\begin{equation}
\label{EVN}
\Po[D_{N}\mid D_{0}\ldots D_{N-1}]
\geq \frac{1}{4}\Big[\frac{t^{-\eps}\ln^{-2K_0}t}{(2m'_2(N))^2}\Big]^{J^1_1(N)}
=t^{o(1)}
\end{equation}
as $t\rightarrow \infty$ and $\omega \in \Gamma_t\cap \Lambda_t$.

%Conclusion
\subsubsection{Conclusion}
\label{conclow2}
Using the results of the precedent subsections, we now aim at showing that 
\begin{equation}
\label{ETL}
\prod_{i=1}^{N}\Po[D_{i}\mid D_{1} \ldots D_{i-1}] \geq t^{-(r_1(N)+o(1))} 
\end{equation}
as $t \rightarrow \infty$ and $\omega \in \Gamma_{t}\cap \Lambda_{t}$.\\
We will proceed by induction. For every $i\geq 1$ consider the product 
\[\Pi_i=\prod_{j=1}^{i}\Po[D_{j}\mid D_{1} \ldots D_{j-1}].\]
Then we make the following induction hypothesis 
\begin{equation*}
\Pi_i\geq t^{-(r_1(i+1)+o(1))} 
\end{equation*}
for some $1\leq i \leq N-2$ ( we suppose $N>2$ otherwise (\ref{ETL}) is trivial).
We will show that it implies 
\begin{equation*}
\Pi_{i+1}\geq t^{-(r_1(i+2)+o(1))}. 
\end{equation*}
Indeed we have,
\begin{equation*}
\Pi_{i+1}= \Pi_i \Po[D_{i+1}\mid D_{1} \ldots D_{i}]. 
\end{equation*}
If $m'_2(i+1)$ is not $t^{a_{i+1}}$-stable then by (\ref{EVi}), we have $\Po[D_{i+1}\mid D_{1} \ldots D_{i}]\geq t^{o(1)}$ (since in this case $(r_2(i)-a_i)^+=0$) and by the induction hypothesis
\begin{equation*}
\Pi_{i+1}\geq t^{-(r_1(i+1)+o(1))}. 
\end{equation*}
But in this case note that by construction we have $r_1(i+2)=r_1(i+1)$.
So we obtain
\begin{equation*}
\Pi_{i+1}\geq t^{-(r_1(i+2)+o(1))}.
\end{equation*}
If $m'_2(i+1)$ is $t^{a_{i+1}}$-stable we have by (\ref{EVi}) \[\Po[D_{i+1}\mid D_{1} \ldots D_{i}]\geq t^{-(r_{2}(i+1)-a_{i+1}+o(1))}\]
which implies with the induction hypothesis
\begin{equation*}
\Pi_{i+1}\geq t^{-(r_1(i+1)-a_{i+1}+r_{2}(i+1)+o(1))}.
\end{equation*}
But in this case we have by construction either $r_1(i+2)=r_1(i+1)-a_{i+1}+r_2(i+1)$ or $r_1(i+2)=r_2(i+1)$ and $r_1(i+1)=a_{i+1}$. This leads to 
\begin{equation*}
\Pi_{i+1}\geq t^{-(r_1(i+2)+o(1))}.
\end{equation*}
As the induction hypothesis is verified for $i=1$ by (\ref{EV1}), we showed that for all $1\leq i<N$ we have
\begin{equation*}
\prod_{j=1}^{i}\Po[D_{j}\mid D_{1} \ldots D_{j-1}] \geq t^{-(r_1(j+1)+o(1))} 
\end{equation*}
as $t \rightarrow \infty$ and $\omega \in \Gamma_{t}\cap \Lambda_{t}$.\\
Now taking $i=N-1$ and using (\ref{EVN}), we obtain
\begin{equation}
\label{EVN1}
\prod_{j=1}^{N}\Po[D_{j}\mid D_{1} \ldots D_{j-1}] \geq t^{-(r_1(N)+o(1))} 
\end{equation}
as $t \rightarrow \infty$ and $\omega \in \Gamma_{t}\cap \Lambda_{t}$.\\
By Lemma \ref{l_good_omega}, (\ref{EV0}) and (\ref{EVN1}),  we obtain that $\IP$-a.s. 
\begin{equation}
\label{LAST0}
\Po[T_2>t]\geq t^{-(r_1(N)+o(1))}=t^{-(\zeta_2(t)+o(1))}
\end{equation}
as $t\rightarrow \infty$.
Together with (\ref{RE1}), the inequality (\ref{LAST0}) proves Theorem \ref{theo1}.
\medskip

To conclude, we deduce the second part of Theorem \ref{theo2}. Note first that by Theorem \ref{theo3}, we have for $\theta>0$ and every $t\geq e$
\begin{equation}
\label{LAST1}
\IP[\zeta_{2}(t)<1+\theta]=\IP[\zeta_{2}(e)<1+\theta]>0
\end{equation}
the last inequality is obtained using the fact that, as $W$ is a Brownian motion, $\zeta_2(e)$ is absolutely continuous in relation to the Lebesgue measure with density non almost surely zero in $[1,1+\theta]$.
From (\ref{LAST1}) and as the sequence $(\zeta_{2}(n))_{n\geq 3}$ is ergodic, there exists an increasing subsequence $(n_i)_i$ such that $1\leq \zeta_2(n_i)<1+\theta$, $\IP$-a.s.\ for all $i\geq 1$.
Now by Markov's inequality, we get 
\begin{equation}
\label{LAST2}
\Eo[T_2^{c}]\geq \Po[T_2>n_i]n_{i}^{c}
\end{equation}
for every $n_i$.
Thus from (\ref{LAST0}), (\ref{LAST2}) and the fact that $1\leq \zeta_2(n_i)<1+\theta$ for the subsequence $(n_i)_i$ we get that 
\begin{equation}
\Eo[T_2^{c}]\geq n_{i}^{(c-(1+\theta)+o(1))}\nonumber\\
\end{equation}
for $i$ sufficiently large. If $c>1$, taking $\theta<c-1$ and letting $i\rightarrow \infty$ we obtain $\Eo[T_2^{c}]=+\infty$, $\IP$-a.s.

%%%%%%%%%%%%%%%%%%%%%%%%%%%%%%%%%%%%%%%%%%%%%%%%%%%%%%%%%%%%%%%%%%%%%%%%%%%%%%%%%%%%%%%%%%%%%%%%%%%%%%%%%%%%%%%%%%%%%%%%%%%%%%%%%%%%%%

\section{Proofs of  Theorems \ref{theo1} and \ref{theo2} in the general case}
\label{SECK}

\subsection{Upper bound for $\Po[T_\gam>t]$}
First, as we did in the case $\gam=2$, let us find an upper bound for $\Po[T_\gam>t^{1+\delta}]$ for $\delta>0$ arbitrary instead of $\Po[T_\gam>t]$.
Then we write
\begin{equation*}
 \Po[T_\gam>t^{1+\delta}]= \Po[T_{\{\gam-1,\gam\}}>t^{1+\delta}\mid T_{\gam-1}>t^{1+\delta}]\Po[T_{\gam-1}>t^{1+\delta}]
\end{equation*}
where $T_{\{\gam-1,\gam\}}=\inf\{s>0;\xi_{\gam-1}(s)=\xi_\gam(s)\}$.
Let us suppose that for some $\gam>2$,
\begin{equation}
\label{INDHYP}
\Po[T_{\gam-1}>t^{1+\delta}]\leq t^{-(\sum_{i=1}^{\gam-2}(\gam-i)r_i(N)+o(1))}
\end{equation}
as $t\rightarrow \infty$ and $\omega \in \Gamma_t\cap \Lambda_t$.
The goal of this section is to show that this implies 
\begin{equation}
\label{IND}
\Po[T_\gam>t^{1+\delta}]\leq t^{-(\sum_{i=1}^{\gam-1}(\gam-i)r_i(N)+o(1))}
\end{equation}
as $t\rightarrow \infty$ and $\omega \in \Gamma_t\cap \Lambda_t$.
As by subsection \ref{Subsecup2}, (\ref{INDHYP}) is true for $\gam=3$, we will conclude by induction that 
(\ref{IND}) is true for all $\gam\geq 2$.

Now, we have to bound from above the term $\Po[T_{\{\gam-1,\gam\}}>t^{1+\delta}\mid T_{\gam-1}>t^{1+\delta}]$ for $\gam \geq 3$ to obtain an upper bound for $\Po[T_\gam>t^{1+\delta}]$. For the sake of brevity, let us write 
$ \Po^*[\cdot]:=\Po[\phantom{*}\cdot \mid T_{\gam-1}>t^{1+\delta}]$.

Fix $t>1$ and define, for $1\leq i \leq N$, the following stopping times
\begin{equation}
\label{STP}
\sigma_{i} =\inf\{s>0;\xi_\gam(s)= h_{\gam-1}(i)\}, 
\end{equation}
the following events
\begin{equation}
\label{EVT}
B_{i}=\{ T_{\{\gam-1,\gam\}} \notin[\sigma_{i-1},\sigma_{i}]\} 
\end{equation} 
for $1\leq i \leq N-1$, with the convention $\sigma_{0}=0$
and
\begin{equation}
B_{N}=\{T_{\{\gam-1,\gam\}} \notin[\sigma_{N-1},t^{1+\delta}]\}\cap \Big\{\sigma_{N-1}<\frac{t^{1+\delta}}{2}\Big\}. \nonumber\\ 
\end{equation}
To find the upper bound for $\Po^*[T_{\{\gam-1,\gam\}}>t^{1+\delta}]$, we will use the following decomposition
\begin{align}
\label{d1K}
\Po^*[T_{\{\gam-1,\gam\}}>t^{1+\delta}]
&=
\Po^*\Big[T_{\{\gam-1,\gam\}}>t^{1+\delta},\sigma_{N-1}< \frac{t^{1+\delta}}{2}\Big]\nonumber\\
&\phantom{**}+\Po^*\Big[T_{\{\gam-1,\gam\}}> t^{1+\delta},\sigma_{N-1}\geq \frac{t^{1+\delta}}{2}\Big].
\end{align}\\
As the second term of the right-hand side of (\ref{d1K}) is completely similar to the term $\Po^*[B_{1}, \sigma_1\geq t^{\alpha+\eps}]$ of section \ref{UPPERKparB_1}, we can show that it is negligible in comparison with the first one.
For the upper bound of the first term of (\ref{d1K}) we start by noting that
\begin{eqnarray}
\label{produitK}
\lefteqn{\Po^*\Big[T_{\{\gam-1,\gam\}}>t^{1+\delta},\sigma_{N-1}< \frac{t^{1+\delta}}{2}\Big]}\phantom{*********}\nonumber\\
&\leq& 
\Po^*[B_{1}\cap \ldots \cap B_{N}]
\nonumber\\
&=& \Po^*[B_{1}]\Po^*[B_{2}\mid B_{1}]\ldots \Po^*[B_{N}\mid B_{1}B_{2}\ldots B_{N-1}].
\end{eqnarray}
In the next subsections we will find upper bounds for the terms of the product of the right-hand side of (\ref{produitK}). 

\subsubsection{Upper bound for $\Po^*[B_{1}]$}
\label{UPPERKparB_1}
We start by writing
\begin{equation*}
\Po^*[B_{1}]\leq \Po^*[B_{1}, \sigma_1\geq t^{\alpha+\eps}]+\Po^*[\sigma_1< t^{\alpha+\eps}]
\end{equation*}
where $\eps=\eps(t)$ is from Lemma \ref{l_good_omega}.
\medskip

\noindent
{\bf Upper bound for $\Po^*[\sigma_1< t^{\alpha+\eps}]$}
\medskip

\noindent
By independence of the random walks and the Markov property, observe that 
\begin{equation*}
\Po^*[\sigma_1< t^{\alpha+\eps}]\leq \prod_{j=1}^{\gam-1}\Po^{h_{j-1}(1)}[\tau_{h_j(1)}(\xi_{\gam})<t^{\alpha+\eps}]
\end{equation*}
with the convention $h_0(1)=\gam$.
Now using (\ref{CER}) to bound from above each term of the product, we obtain
\begin{equation}
\Po^*[\sigma_1< t^{\alpha+\eps}]\leq C_4^{\gam-1}\prod_{j=1}^{\gam-1}t^{-(r_j(1)-(\alpha+\eps))}
\end{equation}
for $\omega \in \Gamma_t\cap \Lambda_t$.
\medskip

\noindent
{\bf Upper bound for $\Po^*[B_{1}, \sigma_1\geq t^{\alpha+\eps}]$}
\medskip

\noindent
We will show that the term $\Po^*[B_{1}, \sigma_1\geq t^{\alpha+\eps}]$ is negligible in comparison with $\Po^*[\sigma_1< t^{\alpha+\eps}]$. Let us first treat two important cases.
Denote by $M(1)$ the set of the first $\gam-1$ $t^{\alpha}$-stable points. That is
\[
M(1)=\{m_1(1), m_2(1),\dots,m_{\gam-1}(1)\}.
 \]
Then, take $\eps'(t)= (28+8K_0) \ln_2 t\ln^{-1} t$,
we will first treat the following case.
\medskip

\noindent
{\sl Case 1:} $a_1\geq \alpha+2\eps'$
\medskip

\noindent
We define the following event
\[
 C=\{\mbox{There exists $s \in [0, t^{\alpha+\eps}]$ such that $\xi_j \in M(1)$, for $1 \leq j \leq \gam$}\}.
\]
Observe that 
\begin{equation}
\Po^*[B_{1}, \sigma_1\geq t^{\alpha+\eps}]\leq \frac{\Po[C^c]}{\Po[T_{\gam}>t^{1+\delta}]}
\end{equation}
for $t$ sufficiently large.\\
In section \ref{secLOWKpart} we will obtain that
\begin{equation}
\label{ESTRELA}
\Po[T_{\gam}>t^{1+\delta}]\geq t^{2\gam^2\ln_3t}
\end{equation}
for $\delta<1$, so let us focus on the term $\Po[C^c]$. Divide the time interval $[0, t^{\alpha+\eps}]$ into $N_0=\lfloor \frac{t^{\alpha+\eps}}{t^{\alpha+\eps'}}\rfloor$ intervals of size $t^{\alpha+\eps'}$ and observe that 
\begin{equation}
\Po[C^c]\leq \Po[\cap_{j=0}^{N_0}\cup_{i=1}^{\gam}\{\xi_i((j+1/2)t^{\alpha+\eps'})\notin M(1)\}]
\end{equation}
By the Markov property and the independence of the random walks we obtain,
\begin{equation}
\Po[C^c]\leq \Big(1-\Big(\min_{x\in [0,h_{\gam-1}(1))}\Po^x[\xi_1(t^{\alpha+\eps'}/2)\in M(1)]\Big)^\gam\Big)^{N_0}
\end{equation}
Now we need to bound from below the quantity $\Po^x[\xi_1(t^{\alpha+\eps'}/2)\in M(1)]$ uniformly in $x$.
Suppose that $x\in [h_{j-1}(1), h_j(1))$ for some $2\leq j \leq \gam-2$. Without loss of generality we can suppose that $x\in [m_j(1),h_j(1))$. Conditioning on the $\sigma$-field ${\cal F}_{\{m_j(1),m_{j+1}(1)\}}$
generated by $\xi_1^x$ up to
the stopping time $ \tau_{\{m_j(1), m_{j+1}(1)\}}$, we get by the
Markov property
\begin{align}
\label{negK}
\lefteqn{\Po^x[\xi_1(t^{\alpha+\eps'}/2)\in M(1)]}\phantom{****}\nonumber\\
&\geq \inf_{s\in [t^{\alpha+\eps'}/4,t^{\alpha+\eps'}/2]} \Po^{mj(1)}[\xi_1(s)=m_j(1)]\nonumber\\
&\phantom{***}\times \Po^x[\tau_{m_{j+1}(1)}(\xi_1)>\tau_{m_{j}(1)}(\xi_1),\tau_{\{m_{j+1}(1),m_{j}(1)\}}(\xi_1)\leq t^{\alpha+\eps'}/4].
\end{align}
For the first term of the right-hand side of (\ref{negK}) we can apply Lemma 3.5 of \cite{CP} with the fact that $a_1\geq \alpha+2\eps'$ to obtain
\begin{equation}
\inf_{s\in [t^{\alpha+\eps'}/4,t^{\alpha+\eps'}/2]} \Po^{mj(1)}[\xi_1(s)=m_j(1)]\geq \frac{1}{2h_{\gam-1}(1)\ln^{2K_0}t} 
\end{equation}
for suficiently large $t$ and $\omega \in \Gamma_t\cap \Lambda_t$.
For the second term of the right-hand side of (\ref{negK}), we first use the decomposition 
\begin{align*}
\lefteqn{\Po^x[\tau_{m_{j+1}(1)}(\xi_1)>\tau_{m_{j}(1)}(\xi_1),\tau_{\{m_{j+1}(1),m_{j}(1)\}}(\xi_1)\leq t^{\alpha+\eps'}/4]}\phantom{*****}\nonumber\\
&\geq \Po^x[\tau_{m_{j+1}(1)}(\xi_1)>\tau_{m_{j}(1)}(\xi_1)]-\Po^x[\tau_{\{m_{j+1}(1),m_{j}(1)\}}(\xi_1)\leq t^{\alpha+\eps'}/4].
\end{align*}
Then using (\ref{exit}) for the first term and (\ref{LOWB2}) for the second term we obtain
\begin{equation*}
\Po^x[\tau_{m_{j+1}(1)}(\xi_1)>\tau_{m_{j}(1)}(\xi_1),\tau_{\{m_{j+1}(1),m_{j}(1)\}}(\xi_1)\leq t^{\alpha+\eps'}/4]\geq \frac{1}{2h_{\gam-1}(1)} 
\end{equation*}
for $t$ sufficiently large and $\omega\in \Gamma_t\cap \Lambda_t$.\\
Finally we obtain
\begin{equation*}
\Po[C^c] \leq \Big(1-\Big(\frac{1}{4h_{\gam-1}^2\ln^{2K_0}t}\Big)^\gam\Big)^{N_0}= o(\exp(-t^{\frac{\eps}{2}}))
\end{equation*}
as $t\rightarrow \infty$ and $\omega\in \Gamma_t\cap \Lambda_t$.
Together with (\ref{ESTRELA}), this last result implies that $\Po^*[\sigma_1< t^{\alpha+\eps}]$ is negligible in comparison with $\Po^*[\sigma_1< t^{\alpha+\eps}]$.
\medskip

Another case which is easy to treat is the case $\max_{1\leq j<\gam-1}\{r_j(1)\vee l_j(1)\}\leq \alpha+2\eps'$.
\medskip

\noindent
{\sl Case 2:} $\max_{1\leq j<\gam-1}\{r_j(1)\vee l_j(1)\}\leq \alpha+2\eps'$
\medskip

\noindent
In this case, we can see that the elevation of the interval $[0, h_{\gam-1}(1)]$ is smaller than $\alpha+2(\gam-2)\eps'$. Then we can apply the spectral gap technique of subsection \ref{Riton} for the term $\Po[B_{1},\sigma_{1}\geq t^{b_1}]$ to show that
\begin{equation*} 
\Po^*[B_{1}, \sigma_1\geq t^{\alpha+\eps}]=o(\exp(-\frac{1}{2}\exp(\ln^{1/12}t)))  
\end{equation*}
as $t \rightarrow \infty$ and $\omega \in \Gamma_t\cap \Lambda_t$.
\medskip

\begin{figure}[!htb]
\begin{center}
\includegraphics[scale= 0.7]{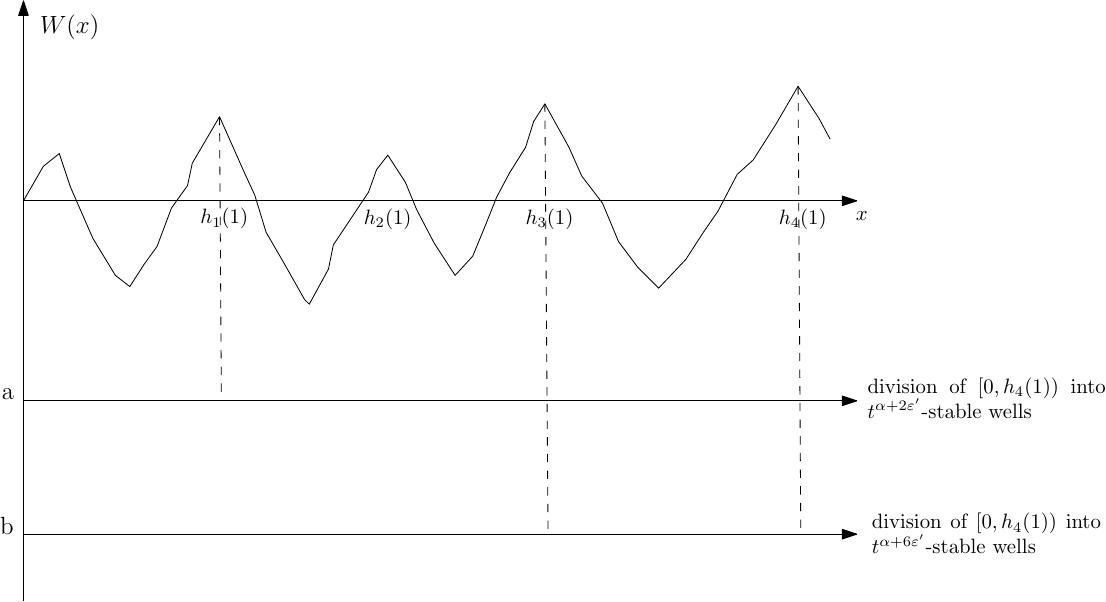}
\caption{Construction of new partitions}
\label{figrenorm}
\end{center}
\end{figure}

{\sl Case 1} and {\sl Case 2} are in fact two limit cases. Let us explain it formally. We want to obtain situations in which we can control the elevation and the height of each well in the interval $[0, h_{\gam-1}(1))$. In the first case, we know that the elevation of each well is smaller than $\alpha\ln t$, so if two random walks stay in a given well a time of order $\alpha+\eps'$ they will meet with high probability (of order $1-\exp(-\ln^9t))$. To guarantee that two particles can stay a sufficiently long time in a well we impose its height to be greater than $(\alpha+2\eps')\ln t$.
In case 2, the situation is different, there is no sufficiently high well that can capture a random walk for a suficiently long time, as a consequence the elevation of the interval $[0, h_{\gam-1}(1))$  is controled and as the random walks will stay in the interval $[0, h_{\gam-1}(1))$ a time much larger than the order of its elevation, we can apply the spectral gap argument of subsection \ref{Riton} to see that the random walks will meet with high probability.
In the intermediate cases, observe that it is possible to find a partition of the interval $[0, h_{\gam-1}(1))$ in a new serie of wells, in a finite number of steps, so that we can go back to one of the two cases above. More precisely, as the number of wells in $[0, h_{\gam-1}(1))$ is finite, we will have two situations
\begin{itemize}
 \item  there exists $4\leq j\leq 2\gam-1$ such that we can define a partition of $[0, h_{\gam-1}(1))$ in $t^{\alpha+j\eps'}$-stable wells such that the elevation of each well is smaller than $(\alpha+(j-2)\eps')\ln t$.
\item  there exists $4\leq j\leq 2\gam-1$ such that we can define a partition of $[0, h_{\gam-1}(1))$ into wells such that the height of each well is smaller than $(\alpha+ j\eps')\ln t$.
\end{itemize}
In order to let things clearer, on Figure \ref{figrenorm} we give an example of such a partition in the case $\gam=5$ and $r_1(1)\wedge l_1(1)\wedge r_3(1)\wedge l_3(1) \geq \alpha+2\eps'$ and $r_2(1)\wedge l_2(1)<\alpha+2\eps'$. On line a, we have a division of the interval $[0, h_{4}(1))$ into $t^{\alpha+2\eps'}$-stable wells. The elevation of the well $(h_1(1), h_3(1))$ is greater than $\alpha\ln t$ but smaller $(\alpha+2\eps')\ln t$, so we have to see if the $t^{\alpha+2\eps'}$-stable wells are also $t^{\alpha+4\eps'}$-stable wells. If this is the case we stop here since we go back to case 1. If this is not the case we take the smaller well, on the picture, $[0,h_1(1))$ and include it in the bigger well $[0, h_3(1))$. Now as the elevation of the well $[0, h_3(1))$ is smaller than $(\alpha+4\eps')\ln t$ we check if both wells $[0, h_3(1))$ and $[h_3(1), h_4(1))$ are $t^{(\alpha+6\eps')}$-stable (Notice that for the well $[h_3(1), h_4(1))$ in fact we just need that $W(h_3(1))-W(m_4(1))\geq (\alpha+6\eps')\ln t$). If this is the case we stop since we go back to case 1. If this not we continue the procedure to finally go back to case 2. 
\medskip

To sum up, in all the cases we obtain,
\begin{equation}
\label{UPBK1}
\Po^*[B_{1}]\leq t^{-(\sum_{j=1}^{\gam-1}r_{j}(1)+o(1))}
\end{equation}
as $t\rightarrow \infty$ and $\omega \in \Gamma_t \cap \Lambda_t$.

%%%%%%%%%%%%%%%%%%%%%%%%%%%%%%%%%%%%%%%%%%%%%%%%%%
\subsubsection{Upper bound for $\Po^*[B_i\mid B_1\dots B_{i-1}]$, $2\leq i <N$}
\label{secUPBKi}
If $h_{\gam-1}(i)=h_{\gam-1}(i-1)$, we obviously have $\Po^*[B_i\mid B_1\dots B_{i-1}]=1$.
If $h_{\gam-1}(i)>h_{\gam-1}(i-1)$, we start by writing
\begin{align*}
\Po^*[B_{i}\mid B_1\dots B_{i-1}]
&\leq \Po^*[B_{i}, \sigma_i-\sigma_{i-1}\geq t^{a_{i-1}+\eps}\mid B_1\dots B_{i-1}]\nonumber\\
&\phantom{**}+\Po[\sigma_i-\sigma_{i-1}< t^{a_{i-1}+\eps}].
\end{align*}

\noindent
{\bf Upper bounds for $\Po[\sigma_i-\sigma_{i-1}< t^{a_{i-1}+\eps}]$}
\medskip

By the Markov property and (\ref{CER}) we obtain,
\begin{equation*}
\Po[\sigma_i-\sigma_{i-1}< t^{a_{i-1}+\eps}]\leq C_4t^{\eps'}t^{-(r_{\gam}(i-1)-a_{i-1})}.
\end{equation*}
\medskip
for $\omega \in \Gamma_t \cap \Lambda_t$.
\medskip

\noindent
{\bf Upper bounds for $\Po^*[B_{i}, \sigma_i-\sigma_{i-1}\geq t^{a_{i-1}+\eps}\mid B_1\dots B_{i-1}]$}
\medskip

We will show that this term is negligible in comparison with $\Po[\sigma_i-\sigma_{i-1}< t^{a_{i-1}+\eps}]$. Observe that
\begin{align}
\label{NEGK2}
\lefteqn{\Po^*[B_{i}, \sigma_i-\sigma_{i-1}\geq t^{a_{i-1}+\eps}\mid B_1\dots B_{i-1}]}\phantom{********}\nonumber\\
&\leq \frac{\Po[B_{i},T_{\gam-1}>t^{1+\delta}, \sigma_i-\sigma_{i-1}\geq t^{a_{i-1}+\eps}\mid B_1\dots B_{i-1}]}{\Po[T_\gam>t^{1+\delta}]}.
\end{align}
We just need to focus on the numerator of (\ref{NEGK2}).
For the sake of brevity we introduce the vector $\vec{x}=(x_1,\dots,x_{\gam-1})$ with all its components smaller than $h_{\gam-1}(i-1)$.
Observe that
\begin{align}
\label{NEGK4}
\lefteqn{\Po[B_{i},T_{\gam-1}>t^{1+\delta}, \sigma_i-\sigma_{i-1}\geq t^{a_{i-1}+\eps}\mid B_1\dots B_{i-1}]}\phantom{*****}\nonumber\\
&\leq \max_{\vec{x}}\Po^{\vec{x}}[B_{i},T_{\gam-1}\notin [\sigma_{i-1}, \sigma_i], \sigma_i-\sigma_{i-1}\geq t^{a_{i-1}+\eps}]\nonumber\\
&\phantom{**}+\Po^{h_{\gam-1}(i-1)}[\sigma_i>t^{1+\delta}].
\end{align}
Using (\ref{LOWB}), we can see that the second term of the right-hand side of (\ref{NEGK4}) is stretched exponential. For the first term, we can apply a similar argument as for the term $\Po^*[B_{1}, \sigma_1\geq t^{\alpha+\eps}]$ in section \ref{UPPERKparB_1} to show that it is $o(\exp(-(1/2)\exp(\ln^{1/12}t)))$ as $t \rightarrow \infty$ and $\omega \in \Gamma_t\cap \Lambda_t$.
\medskip

Therefore, in the case $h_{\gam-1}(i)>h_{\gam-1}(i-1)$, we obtain
\begin{equation}
\label{UPBKi}
\Po^*[B_{i}\mid B_1\dots B_{i-1}]\leq 2C_4t^{\eps}t^{-(r_{\gam}(i-1)-a_{i-1})}=t^{-(r_{\gam}(i-1)-a_{i-1}+o(1))}
\end{equation}
as $t\rightarrow \infty$ and $\omega \in \Gamma_t \cap \Lambda_t$.

%%%%%%%%%%%%%%%%%%%%%%%%%%%%%%%%%%%%%%%%%
\subsubsection{Upper bound for $\Po^*[B_N\mid B_1\dots B_{N-1}]$}
\label{secUPBKN}
If $h_{\gam-1}(N)=h_{\gam-1}(N-1)$, we have $\Po^*[B_i\mid B_1\dots B_{i-1}]=1$.
If $h_{\gam-1}(N)>h_{\gam-1}(N-1)$, we do the following.
We denote by $\vec{x}=(x_1,\dots,x_{\gam-1})$ a vector with all its components smaller than $h_{\gam-1}(N-1)$. 
We use the following decomposition
\begin{align}
\lefteqn{\Po^*[B_N\mid B_1\dots B_{N-1}]}\phantom{*****}\nonumber\\
&\leq \max_{\vec{x}}\Po^{\vec{x}}[T_\gam\notin[\sigma_{N-1},\sigma_N\wedge t^{1+\delta}],\sigma_N\wedge t^{1+\delta}-\sigma_{N-1}\geq t^{a_{N-1}+\eps}]\nonumber\\
&\phantom{**}+\Po[\sigma_N-\sigma_{N-1}<t^{a_{N-1}+\eps}]
\end{align}
where the last inequality holds for sufficiently large $t$ since $a_{N-1}+\eps< 1+\delta/2$ for sufficiently large $t$.
\medskip

\noindent
{\bf Upper bound for $\Po[\sigma_N-\sigma_{N-1}<t^{a_{N-1}+\eps}]$}
\medskip

By the Markov property and (\ref{CER}) we obtain,
\begin{equation}
\label{atal}
\Po[\sigma_N-\sigma_{N-1}<t^{a_{N-1}+\eps}]\leq C_4t^{-(r_{\gam}(N-1)-(a_{N-1}+\eps))}.
\end{equation}
for sufficiently large $t$ and $\omega \in \Gamma_t \cap \Lambda_t$.
\medskip

\noindent
{\bf Upper bound for $\Po^{\vec{x}}[T_\gam\notin[\sigma_{N-1},\sigma_N\wedge t^{1+\delta}],\sigma_N\wedge t^{1+\delta}-\sigma_{N-1}\geq t^{a_{N-1}+\eps}]$}
\medskip
Again to bound this term from below, we can apply the same serie of argument as for the term $\Po^*[B_{i}, \sigma_i-\sigma_{i-1}\geq t^{a_{i-1}+\eps}\mid B_1\dots B_{i-1}]$ of subsection \ref{secUPBKi} to show that it is negligible in comparison with (\ref{atal}).
\medskip

Hence, in the case $h_{\gam-1}(N)>h_{\gam-1}(N-1)$, we obtain
\begin{equation}
\label{UPBKN}
\Po^*[B_{N}\mid B_1\dots B_{N-1}]\leq 2C_4t^{\eps}t^{-(r_{\gam}(N-1)-a_{N-1})}=t^{-(r_{\gam}(N-1)-a_{N-1}+o(1))}
\end{equation}
as $t\rightarrow \infty$ and $\omega \in \Gamma_t \cap \Lambda_t$.

\subsubsection{Conclusion}

By (\ref{UPBK1}), (\ref{UPBKi}), (\ref{UPBKN}) and the construction of section \ref{goodenv} we can deduce by the same type of induction as that we used in subsection \ref{Concup2} that for $\gam\geq 3$,
\begin{equation*}
\Po^*[T_{\{\gam-1,\gam\}}>t]\leq t^{-\Big(\frac{1}{1+\delta}\sum_{j=1}^{\gam-1}r_j(N)+o(1)\Big)}
\end{equation*}
as $t\rightarrow \infty$ and $\omega \in \Gamma_t \cap \Lambda_t$.
This leads to 
\begin{equation*}
\Po[T_\gam>t]\leq t^{-\Big(\frac{\zeta_\gam(t)}{1+\delta}+o(1)\Big)}
\end{equation*}
as $t\rightarrow \infty$ and $\omega \in \Gamma_t \cap \Lambda_t$.
We can immediately deduce the first part of Theorem \ref{theo2}. 
Indeed, as $r_j(N)\geq 1$, for $1\leq j\leq \gam-1$, we obtain $\IP$-a.s.,
\begin{equation*}
\Po[T_\gam>t]\leq t^{-\Big(\frac{\gam(\gam-1)}{2(1+\delta)}+o(1)\Big)}\nonumber\\ 
\end{equation*}
as $t \rightarrow \infty$.
To conclude, fix $c<\frac{\gam(\gam-1)}{2}$, as $\delta$ is arbitrary we can choose $\delta$ sufficiently small such that $\delta<\frac{\gam(\gam-1)-2c}{2c}$. This shows that $\Eo[T_\gam^{c}]<\infty$, $\IP$-a.s.
\medskip

Now, in order to prove Theorem \ref{theo1}, we shall improve (\ref{RE}). That is why we shall consider convergence in $\IP$-probability instead of $\IP$-a.s.\ to obtain
\begin{equation}
\label{REK}
\ln \Po[T_\gam>t]\leq -\zeta_\gam(t)\ln t + o(\ln t)
\end{equation}
as $t \rightarrow \infty$.
By a similar argument as that we used in subsection \ref{Concup2} for the case $\gam=2$ we can obtain (\ref{REK}) in $\IP$-probability.

\subsection{Lower bound for $\Po[T_\gam>t]$}
\label{secLOWKpart}
Let $\tau_{i}^{2}(\xi_j)$ the time of second return to state i for the $j$th random walk.
\medskip

We define the following events
\begin{align*}
D^1_0
&=\cap_{i=1}^{\gam-1} \{ \tau_{\{i-1,i+1\}}(\xi_{i})>\tau_{m'_{\gam}(1)}(\xi_{\gam})\}\nonumber\\
&\phantom{**}\cap \{\xi_{\gam}(s)\leq \xi_{\gam}(t),0\leq s < t \leq \tau_{m'_{\gam}(1)}(\xi_\gam) \} 
\end{align*} 
and for  $2 \leq j \leq \gam-1$,
\begin{align*}
D^j_{0}
&= \cap_{i=1}^{\gam-j} \{ \tau_{\{i-1,i+1\}}(\xi_{i})>\tau_{m_{\gam-j+1}(1)}(\xi_{\gam-j+1})\}\nonumber\\
&\phantom{**}\cap_{i=\gam-j+2}^{\gam-1}  \{ \tau^{2}_{m_{i}(1)-1}(\xi_{i})\wedge \tau_{m_{i}(1)+1}(\xi_{i})>\tau_{m_{\gam-j+1}(1)}(\xi_{\gam-j+1})\}\nonumber\\
&\phantom{**}\cap  \{ \tau^{2}_{m'_{\gam}(1)-1}(\xi_{\gam})\wedge \tau_{m'_{\gam}(1)+1}(\xi_{\gam})>\tau_{m_{\gam-j+1}(1)}(\xi_{N-j+1})\}\\
&\phantom{**}\cap \{\xi_{\gam-j+1}(s)\leq \xi_{\gam-j+1}(t),0\leq s < t \leq \tau_{m_{\gam-j+1}(1)}(\xi_{\gam-j+1}) \}
\end{align*}
with the convention $\cap_{i=\gam}^{\gam-1}\{\dots \}=\Omega$ the sure event.
We define $D_0=\cap_{j=1}^{\gam-1}D_0^j$.
For $1\leq n \leq N$,
define 
\[
\beta(n)=\max_{2\leq j\leq \gam}\{\tau_{m_j(n)}(\xi_j)\}.
\]
Next, for $1\leq n < N$,
if $a_n=l_{\gam-1}(n)$ then define
\begin{align*}
D_{n}
&=\{ \tau^2_{h_{\gam-1}(n)}(\xi_{\gam})>\tau_{m'_\gam(n+1)}(\xi_{\gam}),\tau_{m'_\gam(n+1)}(\xi_{\gam})-\beta(n)\leq t^{a_n}\}\nonumber\\
&\phantom{**} \cap_{j=1}^{\gam-1} \{(\tau^2_{h_{j-1}(n)}(\xi_j)\wedge \tau_{h_{j}(n)}(\xi_j)) -\beta(n)>t^{a_{n}}\} \nonumber
\end{align*} 
with the convention $\tau^2_{h_0(n)}(\xi_1)=\tau_{h_1(n)}(\xi_1)$.\\
If $a_n= r_i(n)$ for some $1\leq i \leq \gam-1$ or $l_{i}(n)$ for some $1\leq i \leq \gam-2$  then define
\begin{align*}
D_{n}
&=\cap_{j=i}^{\gam}\{ \tau^2_{h_{j-1}(n)}(\xi_{j})>\tau_{m_{j+1}(n)}(\xi_{j}),\theta_jt^{a_n}<\tau_{m_{j+1}(n)}(\xi_{j})-\beta(n)\leq \theta'_j t^{a_n}\}\nonumber\\
&\phantom{**} \cap_{j=1}^{i-1} \{(\tau^2_{h_{j-1}(n)}(\xi_j)\wedge \tau_{h_{j}(n)}(\xi_j))-\beta(n)>t^{a_{n}}\} \nonumber\\
&\phantom{**} \cap_{j=i}^{\gam} \{\tau_{h_{j+1}(n)}(\xi_j)-\beta(n)>t^{a_{n}}\} 
\end{align*} 
with the conventions $\tau^2_{h_0(n)}(\xi_1)=\tau_{h_1(n)}(\xi_1)$, $\cap_{j=1}^{0}\{\dots \}=\Omega$, \[h_{\gam+1}(n)=\min_{x>m_{\gam+1}(n)}\{W(x)-W(m_{\gam+1}(n))=a_n\ln t\}\]
and 
\begin{equation*}
\theta_j=\frac{\gam-j}{\gam-i+1},\phantom{**}
\theta'_j=\frac{\gam-j+1}{\gam-i+1}.
\end{equation*}
Finally define
\begin{align*}
D_{N}
&= \{ \tau^2_{h_{\gam-1}(N)}(\xi_{\gam})-\beta(N))>t\} \nonumber\\
&\phantom{**} \cap_{j=1}^{\gam-1} \{(\tau^2_{h_{j}(N)}(\xi_j)\wedge \tau_{h_{j+1}(N)}(\xi_j))-\beta(N)>t\} \nonumber\\
&\phantom{**}\cap \{\tau_{h_1(N)}(\xi_1)-\beta(N)>t\}.\nonumber\\
\end{align*} 

On the event $D_0$ we place all the $\gam$ random walks at the points $m_j(1)$, for $1\leq j\leq \gam-1$ and $m'_\gam(1)$. Then, on the event $D_n$ for $1\leq n \leq N-1$, suppose that the smallest well is $I(n)=[h_{i-1}(n),h_{i}(n)]$, that is to say, $r_i(n)\wedge l_{i-1}(n)=a_n$. By (\ref{LOWB}), we know that the random walk in the well $I(n)$ could escape from it in a time of order $t^{a_n+\eps}$ with high probability. In order for the random walk which is in this well not to meet the random walks which are in the neighboring wells, we will oblige this random walk and those that are in the following wells to move to the right to occupy the first $\gam-i+1$ $t^{a_n}$-stable wells to the right of the well $I(n)$. Futhermore, we want this to occur in a time of order less than $t^{a_n}$, avoiding the random walks meeting by controling the intervals of time each random walk moves to its well to the next well. The random walks to the left of the well $I(n)$ stay in their respective wells. Finally on $D_N$ we oblige all the random walks to stay in their respective $t$-stable wells until time $t$.
\medskip

Observe that
\begin{equation}
\label{decompLow}
\Po[T_\gam>t] \geq  \Po[D_{0}]\Po[D_{1}\mid D_{0}]\ldots \Po[D_{N}\mid D_{0}D_{1}\ldots D_{N-1}].
\end{equation}
We will now find a lower bound for every term of the right-hand side of (\ref{decompLow}).

\subsubsection{Lower bound for $\Po[D_{0}]$}

First we will use the decomposition
\begin{align}
\label{decompLow1}
\Po[D_0] &= \Po\Big[D_0= \cap_{j=1}^{\gam-1} D^{j}_0\Big]\nonumber\\
&= \Po\Big[D^1_{0}\Big]\Po\Big[D^2_{0}\mid D^1_{0}\Big]\ldots \Po\Big[D^{\gam-1}_{0}\mid D^1_{0}D^2_{0}\ldots D^{\gam-2}_{0}\Big].
\end{align}
As we did for the case of two particles, it is not hard to find that 
\begin{equation}
\label{POS2}
\Po\Big[D^1_{0}\Big]\geq \Big(\frac{e^{-(\gam-1)\kappa^2}(e-1)}{2e\kappa^2}\Big)^{m'_{\gam}(1)}
\end{equation}
and for $2\leq j\leq \gam-1$,
\begin{equation}
\label{POS3}
\Po\Big[D^j_{0}\mid D^1_0 \dots D^{j-1}_{0}\Big]\geq \Big(\frac{e^{-(\gam-1)\kappa^2}(e-1)}{2e\kappa^2}\Big)^{m_{\gam-j+1}}.
\end{equation}
Finally using (\ref{decompLow1}), (\ref{POS2}) and (\ref{POS3}) we can obtain a lower bound for $\Po[D_{0}]$
\begin{align}
\label{LBK0}
\Po[D_{0}]
&\geq  \Big(\frac{e^{-(\gam-1)\kappa^2}(e-1)}{2e\kappa^2}\Big)^{m'_\gam(1)+\sum_{j=2}^{\gam-1}m_j(1)}\nonumber\\
&\geq  \Big(\frac{e^{-(\gam-1)\kappa^2}(e-1)}{2e\kappa^2}\Big)^{\gam\ln^{1/2}t}
\end{align}\\
for $\omega \in \Gamma_t \cap \Lambda_t$.

\subsubsection{Lower bounds for $\Po[D_{n}\mid D_{0}\dots D_{n-1}]$, $1\leq n<N$}
\label{secD1}
We start by noting that by definition of $D_n$ and by the Markov property we have
\begin{align}
\lefteqn{\Po[D_{n}\mid D_{0}\dots D_{n-1}]}\phantom{***}\nonumber\\
&=\sum_{x_1\in I_1(n)}\sum_{x_2\in I_2(n)}\dots \sum_{x_\gam\in I_\gam(n)} \Po[D_n\mid \xi_1(\beta_n)=x_k,\dots,\xi_\gam(\beta_n)=x_\gam]\nonumber\\
&\phantom{*****} \times \Po[\xi_1(\beta_n)=x_k,\dots,\xi_\gam(\beta_n)=x_\gam\mid D_{0}\dots D_{n-1}]. 
\end{align}
where if $n>1$, $I_1(n)=[0,h_1(n))$, $I_j(n)=(h_{j-1}(n),h_j(n))$ for $2\leq j< \gam$
and $I_\gam(n)=(h_{\gam-1}(n),h_{\gam+1}(n-1))$ and
$I_1(1)=\{1\}$, $I_j(1)=m_j(1)$ for $2\leq j< \gam$ and $I_\gam(1)=m'_\gam(1)$.
We will now bound uniformly in $x_1,\dots, x_\gam$, the quantity 
\[\Po[D_n\mid \xi_1(\beta_n)=x_k,\dots,\xi_\gam(\beta_n)=x_\gam].\]

\noindent
{{\bf Case 1:} $a_n=l_{\gam-1}(n)$}
\medskip

\noindent
We have by the Markov property and the independence of the $\gam$ random walks
\begin{align}
\lefteqn{\Po[D_n\mid \xi_1(\beta_n)=x_k,\dots,\xi_\gam(\beta_n)=x_\gam]}\phantom{******} \nonumber\\
&=\Po^{x_\gam}\Big[\tau_{h_{\gam-1}(n)}(\xi_\gam)>\tau_{m'_{\gam}(n+1)}(\xi_\gam),\tau_{m'_{\gam}(n+1)}(\xi_\gam)\leq t^{a_n}\Big]\nonumber\\
&\phantom{**}\times \prod_{j=2}^{\gam-1} \Po^{x_j}[\tau_{h_{j-1}(n)}(\xi_j)\wedge \tau_{h_{j}(n)}(\xi_j)>t^{a_n}]\nonumber\\ 
&\phantom{**}  \times \Po^{x_1}[\tau_{h_1(n)}(\xi_1)>t^{a_n}].
\end{align}
Let us denote
\begin{equation*}
R_1=\Po^{x_\gam}\Big[\tau_{h_{\gam-1}(n)}(\xi_\gam)>\tau_{m'_{\gam}(n+1)}(\xi_\gam),\tau_{m'_{\gam}(n+1)}(\xi_\gam)\leq t^{a_n}\Big],
\end{equation*}
\begin{equation*}
R^j_2=\Po^{x_j}[\tau_{h_{j-1}(n)}(\xi_j)\wedge \tau_{h_{j}(n)}(\xi_j)>t^{a_n}],\phantom{*}\mbox{$2\leq j \leq \gam-1$,}
\end{equation*}
and
\begin{equation*}
R_3=\Po^{x_1}[\tau_{h_1(n)}(\xi_1)>t^{a_n}].
\end{equation*}
By a similar method as that we used in subsection \ref{sub0low} for the term
\[\Po^{m'_2(1)}[\tau_{m'_2(2)}(\xi_2)<\tau_{h_1(1)}(\xi_2),\tau_{m'_2(2)}(\xi_2)\leq t^{a_1}],\]
 we can deduce that for $\omega \in \Gamma_t\cap \Lambda_t$, 
\begin{equation*}
R_1\geq t^{-((r_\gam(n)-a_n)^++o(1))}
\end{equation*}
as $t\rightarrow \infty$.\\
The terms $R^j_2$ $2\leq j \leq \gam-1$ and $R_3$ are quite similar to the term $V_1^i$ of subsection \ref{Hu}, so we can show that $\omega \in \Gamma_t\cap \Lambda_t$,
\begin{equation*}
R^j_2\leq t^{o(1)}, \phantom{**} R_3\leq t^{o(1)}
\end{equation*}
as $t\rightarrow \infty$.

To sum up we obtain in this case
\begin{equation}
\label{LBKn1}
\Po[D_n\mid D_0\dots D_{n-1}]\geq t^{-((r_\gam(n)-a_n)^++o(1))}
\end{equation}
as $t\rightarrow \infty$ and $\omega \in \Gamma_t\cap \Lambda_t$.
\medskip

\noindent
{{\bf Case 2:} $a_n= r_i(n)$ for some $1\leq i \leq \gam-1$ or $l_{i}(n)$ for some $1\leq i \leq \gam-2$}
\medskip

\noindent
As an example let us treat the case $i>1$.
Using the Markov property and the independence of the $\gam$ random walks we obtain the following decomposition

\begin{align}
\label{EQ1J}
\lefteqn{\Po[D_n\mid \xi_1(\beta_n)=x_k,\dots,\xi_\gam(\beta_n)=x_\gam]}\phantom{*******} \nonumber\\
&=  \prod_{j=i}^{\gam} \Po^{x_j}\Big[\tau_{h_{j-1}(n)}(\xi_j)>\tau_{m_{j+1}(n)}(\xi_j),\nonumber\\
&\phantom{*******}\theta_jt^{a_n}<\tau_{m_{j+1}(n)}(\xi_j)\leq \theta'_jt^{a_n},\tau_{h_{j+1}(n)}(\xi_j)>t^{a_{n}}\Big]\nonumber\\
&\phantom{**}\times \prod_{j=2}^{i-1} \Po^{x_j}[\tau_{h_{j-1}(n)}(\xi_j)\wedge \tau_{h_{j}(n)}(\xi_j)>t^{a_n}]\nonumber\\ 
&\phantom{**}  \times \Po^{x_1}[\tau_{h_1(n)}(\xi_1)>t^{a_n}].
\end{align}
Let us denote
\begin{align*}
U^j_1
&=\Po^{x_j}\Big[\tau_{h_{j-1}(n)}(\xi_j)>\tau_{m_{j+1}(n)}(\xi_j),\nonumber\\
&\phantom{*****}\theta_jt^{a_n}<\tau_{m_{j+1}(n)}(\xi_j)\leq \theta'_jt^{a_n},\tau_{h_{j+1}(n)}(\xi_j)>t^{a_{n}}\Big]
\end{align*}
for $i\leq j \leq \gam$.
As the other terms of the right-hand side of (\ref{EQ1J}) are completely similar to the terms $R^j_2$ and $R_3$ of the previous paragraph, the new dificulty here comes from the terms $U^j_1$. We will now bound from below these terms.

Conditioning on the the $\sigma$-field  ${\cal F}_{m_{j+1}(n)}$  generated by $\xi_j$ up to the stopping time $\tau_{m_{j+1}(n)}$, we get by the Markov property that
\begin{align}
\label{cond1}
U^j_1
&\geq \Po^{x_j}\Big[\tau_{h_{j-1}(n)}(\xi_j)>\tau_{m_{j+1}(n)}(\xi_j),\theta_jt^{a_n}<\tau_{m_{j+1}(n)}(\xi_j)\leq \theta'_jt^{a_n}\Big]\nonumber\\
&\phantom{**} \times \Po^{m_{j+1}(n)}[\tau_{h_{j+1}(n)}(\xi_j)>t^{a_{n}}]\nonumber\\
&:=V^j_1\times V^j_2.
\end{align}
Applying (\ref{CER}), we can see that $V^j_2\geq t^{o(1)}$ as $t\rightarrow \infty$ and $\omega \in \Gamma_t\cap \Lambda_t$. Now it remains to bound from below $V_1^j$. First of all, let us denote 
\begin{equation*}
A_j=\{\tau_{h_{j-1}(n)}(\xi_j)>\tau_{m_{j+1}(n)}(\xi_j)\}
\end{equation*}
and
\begin{equation*}
B_j=\{\theta_j t^{a_n}<\tau_{m_{j+1}(n)}(\xi_j)\leq \theta'_j t^{a_n}\}.
\end{equation*}
Conditioning on the the $\sigma$-field  ${\cal F}_{m_{j}(n)}$  generated by $\xi_j$ up to the stopping time $\tau_{m_{j}(n)}$, we get by the Markov property that
\begin{align}
\label{EQi}
V^j_1
&\geq \Po^{m_j(n)}[\tau_{h_{j-1}(n)}(\xi_j)>\tau_{m_{j+1}(n)}(\xi_j),\theta_j t^{a_n}<\tau_{m_{j+1}(n)}(\xi_j)\leq (\theta'_j-\lambda)t^{a_n}]\nonumber\\
&\phantom{**}\times \Po^{x_j}[\tau_{m_j(n)}(\xi_j)< \tau_{\{h_{j-1}(n),h_j(n)\}}(\xi_j), \tau_{m_j(n)}(\xi_j)\leq \lambda t^{a_n}]
\end{align}
with $\lambda$ a positive constant smaller than $\frac{1}{\gam+1}$.
The second term of the right-hand side of (\ref{EQi}) is similar to the term (\ref{PU2}), so we can deduce that 
\begin{equation}
 \Po^{x_j}[\tau_{m_j(n)}(\xi_j)< \tau_{\{h_{j-1}(n),h_j(n)\}}(\xi_j), \tau_{m_j(n)}(\xi_j)\leq \lambda t^{a_n}]\geq t^{o(1)}
\end{equation}
as $t\rightarrow \infty$ and $\omega \in \Gamma_t\cap \Lambda_t$. 
Now we treat the first term of the right-hand side of (\ref{EQi}).
Writing
\begin{equation*}
E_j=\{\tau_{h_{j-1}(n)}(\xi_j)>\tau_{m_{j+1}(n)}(\xi_j)\}
\end{equation*}
\begin{equation*}
F_j=\{\tau_{m_{j+1}(n)}(\xi_j)\leq (\theta'_j-\lambda)t^{a_n}\},
\end{equation*}
we obtain
\begin{align}
\label{EQ2}
\lefteqn{\Po^{m_j(2)}[E_j,F_j,\tau_{m_{j+1}(n)}(\xi_j)>\theta_j t^{a_n}]}\nonumber\\ 
&= \Po^{m_j(n)}[E_j,F_j \mid \tau_{\{h_{j-1}(n),h_j(n)\}}(\xi_j)>\theta_j t^{a_n}] \Po^{m_j(n)}[\tau_{\{h_{j-1}(n),h_j(n)\}}(\xi_j)>\theta_j t^{a_n}].
\end{align}
Using (\ref{CER}), we can show that the last term of the right-hand side of (\ref{EQ2}) is greater than $1/2$ for $t$ large enough and $\omega \in \Gamma_t\cap \Lambda_t$. We will now  focus on the first term. First observe that
\begin{align}
 \lefteqn{\Po^{m_j(n)}[E_j,F_j \mid \tau_{\{h_{j-1}(n),h_j(n)\}}(\xi_j)>\theta_j t^{a_n}]}\phantom{**} \nonumber\\
&\geq \Po^{m_j(n)}[E_j,F_j,\xi_j(\alpha_j t^{a_n})\in [\eta_j(n), \gamma_j(n)]\mid \tau_{\{h_{j-1}(n),h_j(n)\}}(\xi_j)>\theta_j t^{a_n}]
\end{align}
where the interval $[\eta_j(n), \gamma_j(n)]$ is such that 
\begin{equation*}
\eta_j(n)= \max_{x<m_j(n)}\{W(x)-W(m_j(n))=\eps \ln t \}
\end{equation*}
and
\begin{equation*}
\gamma_j(n)= \min_{x>m_j(n)}\{W(x)-W(m_j(n))=\eps \ln t \}
\end{equation*}
and $\eps=\eps(t)$ is from Lemma \ref{l_good_omega} (observe that for $t$ sufficiently large these points are well defined).

Then using the Markov property we obtain
\begin{align*}
\lefteqn{ \Po^{m_j(n)}[E_j,F_j,\xi_j(\theta_j t^{a_n})\in [\eta_j(n), \gamma_j(n)]\mid \tau_{\{h_{j-1}(n),h_j(n)\}}(\xi_j)>\theta_j t^{a_n}]}\phantom{********}\nonumber\\
 &= \sum_{x=\eta_j(n)}^{\gamma_j(n)}  \Po^{m_j(n)}[E_j,F_j\mid \xi_j(\alpha_j t^{a_n})=x]\nonumber\\
&\phantom{******}\times \Po^{x_j}[\xi_j(\alpha_j t^{a_n})=x \mid \tau_{\{h_{j-1}(n),h_j(n)\}}(\xi_j)>\theta_j t^{a_n}].
\end{align*}
Now, let us find a uniform lower bound of $\Po^{m_j(2)}[E_j,F_j\mid \xi_j(\theta_j t^{a_n})=x]$ in $x$.
By the Markov property and the definitions of the events $E_j$ and $F_j$ we have
\begin{align*}
\lefteqn{\Po^{m_j(n)}[E_j,F_j\mid \xi_j(\alpha_j t^{a_2})=x]}\phantom{******}\nonumber\\
&=\Po^{x}[\tau_{h_{j-1}(n)}(\xi_j)>\tau_{m{j+1}(n)}(\xi_j),\tau_{m{j+1}(n)}(\xi_j)\leq (\theta'_j-\theta_j-\lambda)t^{a_n}].
\end{align*}
This last term is similar to the term $R_1$ of this subsection, so we obtain
\begin{equation}
\label{ESTIM1}
 \Po^{x}[\tau_{h_{j-1}(n)}(\xi_j)>\tau_{m{j+1}(n)}(\xi_j),\tau_{m{j+1}(n)}(\xi_j)\leq (\beta_j-\alpha_j-\lambda)t^{a_n}] \geq
  t^{-(r_{j}(n)-a_{n}+o(1))}.
 \end{equation}
as $t\rightarrow \infty$ and $\omega \in \Gamma_t\cap \Lambda_t$. 
On the other hand, using (\ref{LOWB}) and reversibility, we obtain, for $t$ large enough and $\omega \in \Gamma_t \cap \Lambda_t$
\begin{equation}
\label{Marcel}
\Po^{m_j(n)}[\xi_j(\alpha_j t^{a_n}) \in [\eta_j(n),\gamma_j(n)]\mid \tau_{\{h_{j-1}(n),h_j(n)\}}(\xi_j)>\alpha_j t^{a_n}]\geq \frac{1}{2}.
\end{equation}
From (\ref{cond1}), (\ref{ESTIM1}) and (\ref{Marcel}) we obtain for $i\leq j \leq \gam$, 
\begin{equation}
U^j_1 \geq  t^{-(r_{j}(n)-a_{n}+o(1))}\nonumber\\
 \end{equation}
as $t \rightarrow \infty$ and $\omega \in \Gamma_t\cap \Lambda_t$.
\medskip

Finally, in the case $a_n= r_i(n)$ or $l_{i-1}(n)$ for some $1\leq i \leq \gam-1$, we obtain that
\begin{equation}
\label{LBKn2}
 \Po[D_n\mid D_0\dots D_{n-1}]\geq t^{-(\sum_{j=i}^{\gam}r_{j+1}(n)-(\gam-i+1)a_n+o(1))}
\end{equation}
as $t \rightarrow \infty$ and $\omega \in \Gamma_t\cap \Lambda_t$.

\subsubsection{Lower bound for $\Po[D_{N}\mid D_{0}\dots D_{N-1}]$}
We have by the Markov property 
\begin{align*}
\lefteqn{\Po[D_{N}\mid D_{0}\dots D_{N-1}]}\phantom{**}\nonumber\\
&=\sum_{x_1\in I_1(N)}\sum_{x_2\in I_2(N)}\dots \sum_{x_\gam\in I_\gam(N)} \Po[D_2\mid \xi_1(\beta_N)=x_k,\dots,\xi_\gam(\beta_N)=x_\gam]\nonumber\\
&\phantom{*****} \times \Po[\xi_1(\beta_N)=x_k,\dots,\xi_\gam(\beta_N)=x_\gam\mid D_{0}\dots D_{N-1}]. 
\end{align*}
where $I_1=[0,h_1(N))$, $I_j=(h_{j-1}(N),h_j(N))$ for $2\leq j< \gam$
and $I_\gam=(h_{\gam-1}(N),h_{\gam+1}(N-1))$.
We will now bound uniformly in $x_1,\dots, x_\gam$, the quantity 
\[\Po[D_N \mid \xi_1(\beta_N)=x_k,\dots,\xi_\gam(\beta_N)=x_\gam].\]
We have by the Markov property and the independence of the random walks
\begin{align*}
\lefteqn{\Po[D_N \mid \xi_1(\beta_N)=x_k,\dots,\xi_\gam(\beta_N)=x_\gam]}\phantom{***********} \nonumber\\
&= \prod_{j=2}^{\gam-1} \Po^{x_j}[\tau_{\{h_{j}(N),h_{j+1}(N)\}}(\xi_j)>t]\nonumber\\
&\phantom{**}\times \Po^{x_\gam}[\tau_{h_{\gam-1}(N)}(\xi_\gam)>t]
\Po^{x_1}[\tau_{h_1(N)}(\xi_1)>t].
\end{align*}
All these terms are similar to the terms $R_2^j$ and $R_3$ of subsection \ref{secD1}, so we obtain
\begin{equation}
\label{LBKN}
\Po[D_{N}\mid D_{0}\dots D_{N-1}]\geq t^{o(1)} 
\end{equation}
as $t \rightarrow \infty$ and $\omega \in \Gamma_t\cap \Lambda_t$.

\subsubsection{Conclusion}

In this section we will show by induction that the total cost of our strategy immediately after step $n-1$ is bounded from below by  $t^{-(\sum_{j=1}^{\gam-1}(\gam-j)r_{j}(n)+o(1))}$. That is to say, if we denote by
\begin{equation*}
\Pi_{n-1}=\Po[D_{0}]\Po[D_{1}\mid D_{0}]\ldots \Po[D_{n-1}\mid D_{0}D_{1}\ldots D_{n-2}]
\end{equation*}
for $1\leq n \leq N$, we have
\begin{equation*}
 \Pi_{n-1} \geq t^{-(\sum_{j=1}^{\gam-1}(\gam-j)r_{j}(n)+o(1))}
\end{equation*}
as $t \rightarrow \infty$ and $\omega \in \Gamma_t\cap \Lambda_t$.\\
So suppose that 
\begin{equation*}
 \Pi_{n-1} \geq t^{-(\sum_{j=1}^{\gam-1}(\gam-j)r_{j}(n)+o(1))}
\end{equation*}
is true for some $1\leq n <N$. We will show that it implies that 
\begin{equation*}
 \Pi_{n} \geq t^{-(\sum_{j=1}^{\gam-1}(\gam-j)r_{j}(n+1)+o(1))}.
\end{equation*}
\noindent
{{\bf Case 1:} $a_n=l_{\gam-1}(n)$}
\medskip

\noindent
In this case we have by (\ref{LBKn1}),
\begin{equation*}
\Po[D_n\mid D_0\dots D_{n-1}]\geq  t^{-((r_\gam(n)-a_n)^++o(1))}
\end{equation*}

In the case $(r_\gam(n)-a_n)^+>0$, let us compute the total cost using the construction of section \ref{goodenv}.
\begin{align*}
\lefteqn{\sum_{j=1}^{\gam-1}(\gam-j)r_{j}(n)+r_\gam(n)-a_n}\phantom{*******}\nonumber\\
&=\sum_{j=1}^{\gam-2}(\gam-j)r_{j}(n)+ r_{\gam-1}(n)+r_\gam(n)-a_n \nonumber\\
&=\sum_{j=1}^{\gam-2}(\gam-j)r_{j}(n+1)+ r_{\gam-1}(n+1)\nonumber\\
&=\sum_{j=1}^{\gam-1}(\gam-j)r_{j}(n+1)
\end{align*}
where we use the fact that, in this case, $r_{\gam-1}(n+1)=r_{\gam-1}(n)+r_{\gam}(n)-a_n$.
In the case $(r_\gam(n)-a_n)^+=0$, we have directly
\begin{equation*}
 \sum_{j=1}^{\gam-1}(\gam-j)r_{j}(n)=\sum_{j=1}^{\gam-1}(\gam-j)r_{j}(n+1).
\end{equation*}
since $r_{j}(n)=r_{j}(n+1)$ for $1\leq j\leq \gam-1$.\\
Therefore, in both cases $(r_\gam(n)-a_n)^+>0$ and $(r_\gam(n)-a_n)^+=0$ we obtain
\begin{equation*}
 \Pi_{n} \geq t^{-(\sum_{j=1}^{\gam-1}(\gam-j)r_{j}(n+1)+o(1))}.
\end{equation*}

\noindent
{{\bf Case 2:} $a_n= r_i(n)$ for some $1\leq i \leq \gam-1$ or $l_{i}(n)$ for some $1\leq i \leq \gam-2$}
\medskip

\noindent
By (\ref{LBKn2}) we have
\begin{equation}
 \Po[D_n\mid D_0\dots D_{n-1}]\geq t^{-(\sum_{j=i}^{\gam}r_{j+1}(n)-(\gam-i+1)a_n+o(1))}
\end{equation} 
As an example, let us treat the case $a_n=r_i(n)$.
Let us compute the total cost using the construction of section \ref{goodenv}.
\begin{align*}
\lefteqn{\sum_{j=1}^{\gam-1}(\gam-j)r_{j}(n)+\sum_{j=i}^{\gam} r_j(n)-(\gam-i+1)a_n}\phantom{**}\nonumber\\
&= r_\gam(n)+\sum_{j=i+1}^{\gam-1}(\gam-j+1)r_{j}(n)+(\gam-i+1)r_i(n)\nonumber\\
&\phantom{**}+\sum_{j=1}^{i-2}(\gam-j)r_{j}(n)+(\gam-i+1)r_{i-1}(n)-(\gam-i+1)a_n\nonumber\\
&=\sum_{j=i}^{\gam-1}(\gam-j)r_{j}(n+1)+\sum_{j=1}^{i-2}(\gam-j)r_{j}(n+1)+(\gam-i+1)r_{i-1}(n+1)\nonumber\\
&=\sum_{j=1}^{\gam-1}(\gam-j)r_{j}(n+1)
\end{align*}
where we used the fact that for $j<i$, $r_j(n+1)=r_j(n)$ and for $j\geq i$, $r_j(n+1)=r_{j+1}(n)$.
Hence, we obtain 
\begin{equation*}
 \Pi_{n} \geq t^{-(\sum_{j=1}^{\gam-1}(\gam-j)r_{j}(n+1)+o(1))}.
\end{equation*}

By (\ref{LBK0}) and Definition \ref{goodenviron} item (iii), we have
\begin{equation*}
\Pi_0\geq t^{-(\sum_{j=1}^{\gam-1}(\gam-j)r_{j}(1)+o(1))}
\end{equation*}
as $t \rightarrow \infty$ and $\omega \in \Gamma_t\cap \Lambda_t$.
Thus, the induction is right for every $1\leq 1 <N$, in particular we obtain
\begin{equation*}
 \Pi_{N-1} \geq t^{-(\sum_{j=1}^{\gam-1}(\gam-j)r_{j}(N)+o(1))}
\end{equation*}
as $t \rightarrow \infty$ and $\omega \in \Gamma_t\cap \Lambda_t$.
Finally by (\ref{LBKN}), we have
\begin{equation*}
 \Pi_{N} \geq t^{-(\sum_{j=1}^{\gam-1}(\gam-j)r_{j}(N)+o(1))}
\end{equation*}
as $t \rightarrow \infty$ and $\omega \in \Gamma_t\cap \Lambda_t$.
Which implies
\begin{equation*}
 \Po[T_\gam>t]\geq t^{-(\sum_{j=1}^{\gam-1}(\gam-j)r_{j+1}(N)+o(1))}
\end{equation*}
or by definition of the process $\zeta_\gam$,
\begin{equation*}
 \Po[T_\gam>t]\geq t^{-(\zeta_\gam(t)+o(1))}
\end{equation*}
as $t$ tends to infinity and $\omega \in \Gamma_t\cap \Lambda_t$.
Together with (\ref{REK}), this last inequality concludes the proof of Theorem \ref{theo1} in the general case. To obtain the second part of Theorem \ref{theo2} we use a similar argument as in subsection \ref{conclow2} for the case $\gam=2$.

%%%%%%%%%%%%%%%%%%%%%%%%%%%%%%%%%%%%%%%%%%%%%%%%%%%%%%%%%%%%%%%%%%%%%

\section{Proof of Theorem \ref{theo3}}
\label{SECtheo3}
We will start this last section by treating the case $\gam=2$. From this special case, we will deduce the proof of the general case.

\begin{figure}[!htb]
\begin{center}
\includegraphics[scale= 0.5]{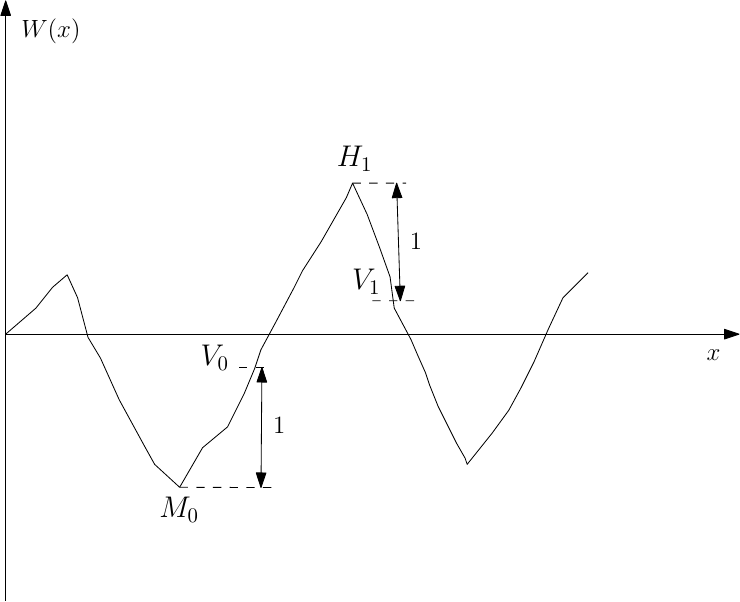}
\caption{On the definition of $V_0$, $M_0$, $V_1$ and $H_1$.}
\label{figul}
\end{center}
\end{figure}

Consider the Brownian motions $W$ and $W'$ given by $W'(\cdot)=\lambda W(\cdot/\lambda^2)$ for $\lambda>0$. Denoting respestively by ${\mathcal S}_t(W)$ and ${\mathcal S}_t(W')$ the sets of the $t$-stable points of $W$ and $W'$, it is elementary to observe that 
\[{\mathcal S}_t(W')=\lambda^2{\mathcal S}_t(W)
\]
and as ${\mathcal H}_t(W)={\mathcal S}_t(-W)$ we also have
\[
{\mathcal H}_t(W')=\lambda^2{\mathcal H}_t(W).
\]
Then using these two properties and expression (\ref{Procze}) of $\zeta_2(t)$ we obtain
\[
\zeta_2(t,W)= \zeta_2(t^{\lambda},W').
\]
As $W$ and $W'$ are identically distributed we deduce $\zeta_2(t)\eqlaw \zeta_2(t^{\lambda})$.
With this scaling property, we observe that the distribution of $\zeta_2(t)$, for fixed $t$, does not depend on $t$. Let us now compute the distribution of $\zeta_2(e)$.
To this end, define the following random variables
\begin{align*}
V_0  &=\inf\{t>0 : W(t)-\inf_{0\leq s\leq t}W(s)=1\},\\
M_0  &= \inf\{t\in [0,V_0] : W(t)=\inf_{0\leq s\leq V_0}W(s)\},\\
V_1  &= \inf\{t>V_0 : \sup_{V_0\leq s\leq t}W(s)-W(t)=1\},\\
H_1  &= \inf\{t\in [V_0,V_1] :
    W(t)=\sup_{V_0\leq s\leq V_1}W(s)\}.
\end{align*}
See Figure \ref{figul}.

By definition of $\zeta_2(e)$, one can immediately see that $\zeta_2(e)=W(H_1)-W(M_0)$.
Then we have
\begin{align*}
\IP[\zeta_2(e)-1>u] &= \IP[W(H_1)-W(M_0)-1>u]\\
&= \IP[W(H_1)-W(V_0)+W(V_0)-W(M_0)-1>u]\\
&= \IP[W(H_1)-W(V_0)>u]
\end{align*}
By Lemma 6.1 (ii) of \cite{CP} we know that the random variable $W(H_1)-W(V_0)$ is exponentially distributed with mean 1. So, this ends the proof of Theorem~\ref{theo3} in the case $\gam=2$. 

In the general case, observe that by the same argument we used for the case $\gam=2$, we show that the distribution of $\zeta_\gam(t)$ does not depend on $t$. Finally to obtain the distribution of $\zeta_\gam(e)$ observe that by the fact that $W$ has independent increments, $\zeta_\gam(e)-\frac{\gam(\gam-1)}{2}$ is the sum of $\gam-1$ independent exponential random variables with parameters $1/(\gam-i)$, $1\leq i \leq \gam-1$. Finally after some elementary computations we get expression (\ref{densityK}). $\square$

\section*{Final comments}

In this last part, we discuss two other natural questions closely related to ours. Suppose that we have $\gam \geq 2$ random walks starting at positions $x_1<x_2<\dots<x_{\gam}$.

The first question concerns the coalescing time of the $\gam$ random walks. Indeed, suppose that whenever two random walks meet they stay together forever. Thus, we can define the coalescing time $T'_{\gam}$ as the moment when all the random walks coalesce. Observe that this time is equal to the time when the first and the $\gam$th random walks meet. Thus, we can deduce that $T'_{\gam}$ behaves like $T_2$ for all $\gam \geq 2$.

The second question concerns the meeting time $T'_{\gam}$ of all the $\gam$ random walks at the same time, that is
\[T''_{\gam}=\inf\{s>0; \xi_1(s)=\xi_2(s)=\dots=\xi_{\gam}(s)\}.
\]
Unlike the case of standard random walks, for which $T''_{\gam}$ is almost surely finite only for $\gam \leq 3$, in the case of random walks in random environment, it is possible to prove, analogously to what is shown in the Appendix of~\cite{Z}, that  for all $\gam\geq 2$, the $\gamma$ random walks will all meet in the origin $\IP$-a.s., and from this follows the finiteness of $T''_{\gam}$. 

In fact, it is possible to show that $T''_{\gam}$ also has the same behavior as the meeting time $T_2$ of two random walks. First, it is possible to show that when the $\gam$ random walks are in the same $t$-stable well they meet at the same time with high probability. Then, in the light of the proofs of Theorems \ref{theo1} and \ref{theo2}, for the $\gam$ random walks not to meet at the same time until time $t$, the strategy which has the "cheapest cost" is to keep the first $\gam-1$ random walks in the first $t$-stable well and to send the $\gam$th random walk into the second $t$-stable well. This is exactly what we have done in the proof of Theorem \ref{theo2} for the case of two random walks. 

To conclude, Theorems \ref{theo1} and \ref{theo2} for $T_2$ should also be true for  $T'_{\gam}$ and $T''_{\gam}$.

\section*{Acknowledgements}
I wish to thank my Ph.D. advisor, Serguei Popov, whose insight was invaluable, and CNPq for financial support.
{

}

\end{document}